\documentclass[11pt]{article}

\usepackage{mathrsfs}
\usepackage{amssymb,amsmath,amsfonts,amsthm}
\usepackage{amsfonts}
\usepackage{graphicx,subfigure,epstopdf}
\usepackage[usenames]{color}
\usepackage{url}
\usepackage{algorithm,algorithmic}
\usepackage{dsfont,verbatim}
\usepackage[colorlinks,linktocpage,linkcolor=blue]{hyperref}
\usepackage[margin=1.0in]{geometry}
\usepackage{enumerate,enumitem}
\usepackage[normalem]{ulem}
\allowdisplaybreaks

\definecolor{darkred}{rgb}{.7,0,0}

\definecolor{green}{rgb}{0,0.7,0}

\newcommand{\nn}{\nonumber}

\newtheoremstyle{thmm}{1.5ex plus 1ex minus .2ex}{1.5ex plus 1ex minus.2ex}{\rmfamily}{}{\bfseries}{}{1em}{} \theoremstyle{thmm}
\newtheorem{theorem}{Theorem}[section]
\newtheorem{lemma}{Lemma}[section]

\newtheorem{corollary}{Corollary}[section]

\newtheorem{remark}{Remark}[section]
\renewenvironment{proof}[1][Proof]{\noindent\textit{#1. }
}{\hfill$\square$}

\renewcommand{\theequation}{\thesection.\arabic{equation}}

\newcommand{\vertiii}[1]
{{\left\vert\kern-0.25ex\left
\vert\kern-0.25ex\left\vert #1
    \right\vert\kern-0.25ex\right
\vert\kern-0.25ex\right\vert}}

\def\d{{\rm d}}

\def\R{\mathbb{R}}
\def\C{\mathbb{C}}

\title{\Large\bf Analyticity, maximal regularity and maximum-norm stability\\ 
of semi-discrete finite element solutions of parabolic equations\\ 
in nonconvex polyhedra 
} 
\author{Buyang Li\thanks{Department of Applied Mathematics, The Hong Kong Polytechnic University, 
Hung Hom, Hong Kong.
Email: bygli@polyu.edu.hk}}
\date{}

\begin{document}
\maketitle
\vspace{-10pt}

\begin{abstract}
In general polygons and polyhedra, possibly nonconvex, the analyticity of the finite element heat semigroup in the $L^q$ norm, $1\leq q\leq\infty$, and the maximal $L^p$-regularity of semi-discrete finite element solutions of parabolic equations are proved. By using these results, the problem of maximum-norm stability of the finite element parabolic projection is reduced to the maximum-norm stability of the Ritz projection, which currently is known to hold for general polygonal domains and convex polyhedral domains. \\

\end{abstract}

\section{Introduction}
\setcounter{equation}{0}

Let $\Omega$ be a polygonal or polyhedral domain in $\R^N$, $N=2,3$, and consider the heat equation 
\begin{align} 
\label{PDE1} 
\left\{\begin{array}{ll}
\displaystyle
\frac{\partial u(t,x)}{\partial t}- \Delta u(t,x) = f(t,x)
&\mbox{for}\,\,(t,x)\in \R_+\times \Omega,\\[10pt]
u(t,x)=0 
&\mbox{for}\,\,(t,x)\in \R_+\times\partial\Omega ,\\[7pt]
u(0,x)=u_0(x)
&\mbox{for}\,\,x\in \Omega . 
\end{array}\right. 
\end{align} 
In the case $f=0$, it is well known that the solution of \eqref{PDE1} is given by $u(t,x)=(e^{t\Delta} u_0)(x)$ , where $E(t)=e^{t\Delta}$ extends to a bounded analytic semigroup on $C_0(\overline\Omega)$ and $L^q(\Omega)$ for arbitrary $1\le q<\infty$ (cf. \cite{Ouhabaz1995}), satisfying the following estimates: 
\begin{align} 
\label{Semigroup} 
\begin{aligned}
&\sup_{t>0}\big(\|E(t)v\|_{L^q(\Omega)}+ t\|\partial_tE(t)v\|_{L^q(\Omega)}\big) \le C\|v\|_{L^{q}(\Omega)}  \, , &&v\in L^q(\Omega),
\,\,\, 1\le q<\infty , \\
&\sup_{t>0}\big(\|E(t)v\|_{C_0(\overline\Omega)}+ t\|\partial_tE(t)v\|_{C_0(\overline\Omega)}\big) \le C\|v\|_{C_0(\overline\Omega)}  ,
&& v\in C_0(\overline\Omega) .
\end{aligned}
\end{align} 
In the case $u_0=0$, the solution of \eqref{PDE1}  possesses the maximal $L^p$-regularity in the space $L^q(\Omega)$, namely, for all $f\in L^p(\R_+;L^q(\Omega))$,  
\begin{align}\label{MaxLpReg}
\|\partial_tu\|_{L^p(\R_+;L^q(\Omega))}
+\|\Delta u\|_{L^p(\R_+;L^q(\Omega))}\leq C_{p,q}
\|f\|_{L^p(\R_+;L^q(\Omega))}, 
\quad \mbox{if}\,\,\, u_0=0 ,\,\, 1<p,q<\infty .
\end{align} 
Such maximal $L^p$-regularity as \eqref{MaxLpReg} has important applications in the analysis of nonlinear partial differential equations (PDEs) \cite{Amann1995,ClementPruss1992,Lions1996}, and has been widely studied in the literature; see \cite{KunstmannWeis2004} and the references therein. 

This paper is concerned with the discrete analogues of \eqref{Semigroup}-\eqref{MaxLpReg}, namely, 
\begin{align}
&\sup_{t>0}\big(\|E_h(t)v_h\|_{L^q(\Omega)}+ t\|\partial_tE_h(t)v_h\|_{L^q(\Omega)}\big) \le C\|v_h\|_{L^{q}(\Omega)}  ,
&& \forall\, v_h\in S_h,\,\,\, 1\le q\le\infty , 
\label{Semigroup-FEM}\\[5pt]
&\|\partial_tu_h\|_{L^p(\R_+;L^q(\Omega))}
+\|\Delta_hu_h\|_{L^p(\R_+;L^q(\Omega))}\leq C_{p,q}
\|f\|_{L^p(\R_+;L^q(\Omega))}, 
&&
\mbox{if}\,\,u_{h,0}=0 ,\,\, 
1<p,q<\infty,
\label{MaxLpReg-FEM}
\end{align} 
where $E_h(t)=e^{t\Delta_h}$ is the semigroup generated by the discrete Laplacian operator $\Delta_h$ (on a finite element subspace $S_h\subset H^1_0(\Omega)$ with mesh size $h$), 
defined by  
\begin{align}\label{Discrete-Laplacian}
(\Delta_h\phi_h,\varphi_h)=-(\nabla \phi_h,\nabla\varphi_h),\quad \forall\, \phi_h, \varphi_h\in S_h, 
\end{align}
and $u_h$ is the finite element solution of \eqref{PDE1}, i.e. 
\begin{align}\label{FEMEq0}
\left\{\begin{array}{ll}
(\partial_t u_h,v_h)+(\nabla u_h ,\nabla v_h)=(f,v_h),
&\forall\, v_h\in S_h ,\,\,
\forall\, t\in(0,T),\\[5pt]
u_h(0)=u_{h,0} . 
\end{array}\right.
\end{align} 
The constants $C$ and $C_{p,q}$ in \eqref{Semigroup-FEM}-\eqref{MaxLpReg-FEM} should be independent of $f$ and $h$. For the maximal $L^p$-regularity \eqref{MaxLpReg-FEM} we require $u_{h,0}=0$ (as the continuous problem), while for error estimate we choose $u_{h,0}$ to be the $L^2$ projection of $u_0$ (see Corollary \ref{max-norm-stability}).

The discrete analyticity \eqref{Semigroup-FEM} and the discrete maximal $L^p$-regularity \eqref{MaxLpReg-FEM} are important mathematical tools for numerical analysis of parabolic equations. For example, \eqref{Semigroup-FEM} can be used to derive error estimates for both semi-discrete and fully discrete finite element methods \cite{Hansbo2002,OdenReddy2012,Palencia1996,Thomee2006}, and \eqref{MaxLpReg-FEM} has been used to study the convergence rates of finite element solutions of semilinear parabolic equations \cite{Geissert2007} as well as nonlinear parabolic equations with nonsmooth diffusion coefficients \cite{LiSun2015-regularity}. The time-discrete extension of the maximal $L^p$-regularity \eqref{MaxLpReg} has been used to study the stability and convergence of time discretization methods for nonlinear parabolic equations with general (possibly degenerate) nonlinearities \cite{AkrivisLi2017,AkrivisLiLubich2016,KunstmanLiLubich2016}. 

Being the foundation for many existing numerical analyses, the discrete analyticity \eqref{Semigroup-FEM} and the discrete maximal regularity \eqref{MaxLpReg-FEM} have been studied by many authors in the literature. In the case $q=2$, \eqref{Semigroup-FEM} holds trivially \cite[Lemma 3.2]{Thomee2006} and \eqref{MaxLpReg-FEM} is an immediate consequence of \eqref{Semigroup-FEM} due to the Hilbert space structure of $L^2(\Omega)$ (cf. \cite{KaltonoLancien2000}). The discrete analyticity \eqref{Semigroup-FEM} for $q\in[1,\infty]\backslash 2$ is a simple consequence of the result in the end-point case $q=\infty$ (via complex interpolation and duality), which was proved in \cite{SchatzThomeeWahlbin1980} for $N=2$ and $r=1$ and was proved in \cite{NitscheWheeler1982} for $N=1,2,3$ and $r\ge 4$, where $r$ is the degree of finite elements. The general case $N\ge 2$ and $r\ge 1$ was proved in \cite{SchatzThomeeWahlbin1998} and \cite{ThomeeWahlbin2000} for the Neumann and Dirichlet boundary conditions, respectively, and was extended to parabolic equations with nonsmooth diffusion coefficients in \cite{Li2015}. The analyses presented in these works were all restricted to smooth domains. The discrete analyticity \eqref{Semigroup-FEM} was proved in \cite{Rannacher1991} for convex polygons in the case $N=2$ and $r=1$ with a logarithmic factor $|\ln h|^{\frac{3}{2}}$, and was proved in \cite{LiSun2016} for convex polyhedra in the case $N=2,3$ and $r\ge 1$. 
In the presence of an extra logarithmic factor $|\ln h|$, the discrete analyticity \eqref{Semigroup-FEM} can be extended to general two-dimensional polygons (cf. \cite[Theorems 6.1 and 6.3]{Thomee2006}). However, the sharp estimate (without logarithmic factor) of \eqref{Semigroup-FEM} remains open in nonconvex polygons and polyhedra. 

Similarly, \eqref{MaxLpReg-FEM} has been proved in smooth domains and convex polygons/polyhedra \cite{Geissert2006,LiSun2016}. 
The extension of \eqref{MaxLpReg-FEM} to the fully discrete finite element methods has been considered in \cite[section 6]{KovacsLiLubich2016} and \cite{LeykekhmanVexler2016-NM,LiSun2017}, which rely on the semi-discrete results. 
For the lumped mass method, both \eqref{Semigroup-FEM} and \eqref{MaxLpReg-FEM} have been proved in general polygons by using the maximum principle \cite{Crouzeix2003,KemmochiSaito}. 
However, for the finite element method, sharp estimates of \eqref{Semigroup-FEM} and \eqref{MaxLpReg-FEM} remain open in nonconvex polygons and polyhedra. In particular, the techniques used in the existing works rely on the $H^2$-regularity of elliptic equations, which only holds in smooth or convex domains.

It is worth to mention that the proof of the discrete maximal $L^p$-regularity \eqref{MaxLpReg-FEM} is closely related to the proof of the maximum-norm stability (best approximation property) of finite element solutions of parabolic equations, namely, 
\begin{align}\label{maximum-norm}
\|u-u_h\|_{L^\infty(0,T;L^\infty(\Omega))}
\le C\ln(2+1/h)\inf_{\chi_h}\|u-\chi_h\|_{L^\infty(0,T;L^\infty(\Omega))},
\end{align} 
where the infimum extends over all $\chi_h\in L^\infty(0,T;S_h)$, and the logarithmic factor ``$\ln(2+1/h)$" in \eqref{maximum-norm} is sharp for piecewise linear finite elements (possibly removable for higher order finite elements).  
Such a priori $L^\infty$-norm best approximation property has been proved in smooth domains  \cite{Leykekhman2004,Li2015,NitscheWheeler1982,SchatzThomeeWahlbin1980,SchatzThomeeWahlbin1998,ThomeeWahlbin2000} and convex polygons (2D) \cite{Rannacher1991}, but remains open in convex polyhedra and nonconvex polygons/polyhedra, though the maximum-norm a posteriori error estimates for finite element solutions 
of parabolic equations have been derived in general polyhedra \cite{DLM2009}. The a priori $L^\infty$-norm best approximation property has been proved in the fully discrete settings with discontinuous Galerkin time-stepping methods \cite{LeykekhmanVexler2016-2} (the result does not cover the semi-discrete case due to a logarithmic dependence on the time-step size). Related maximum-norm stability for finite element solutions of elliptic equations can be found in \cite{DLSW2012,GLRS2009,LeykekhmanVexler2016,RannacherScott1982,Schatz1980}.

In this paper, we prove \eqref{Semigroup-FEM}-\eqref{MaxLpReg-FEM} in general polygons and polyhedra, possibly nonconvex (cf. Theorem \ref{MainTHM1}), and we reduce \eqref{maximum-norm} to the maximum-norm stability of the Ritz projection (cf. Corollary \ref{MainCor1}). In particular, \eqref{maximum-norm} is proved completely in nonconvex polygons and convex polyhedra (cf. Corollary \ref{max-norm-stability}). 
The proof of these results relies on a dyadic decomposition of the domain $(0,1)\times\Omega=\cup_{*,j}Q_j$ together with some local $L^2H^{1+\alpha}(Q_j)$ and $L^\infty H^{1+\alpha}(Q_j)$ estimates of the Green's function (Lemma \ref{GFEst1}) and a local energy error estimate for finite element solutions of parabolic equations (Lemma \ref{LocEEst}). In contrast to the existing work (cf. \cite{LiSun2016,SchatzThomeeWahlbin1998}), the local energy error estimate used here does not require any superapproximation property of the Ritz projection (which only holds in convex domains). These results help to prove the key lemma (Lemma \ref{LemGm2}) for the proof of our main results. 
The maximal $L^p$-regularity \eqref{MaxLpReg-FEM} is first proved for $p=q$ and then extended to $p\neq q$ by using the singular integral operator approach (Sections \ref{secp=q}--\ref{secpneqq}).

\section{Main results and their consequences}
\setcounter{equation}{0}


Let $L^q=L^q(\Omega)$. Let $\Gamma_h(t,x,x_0) $ be the kernel of the operator $E_h(t)$, i.e. 
\begin{align}
(E_h(t)v_h) (x_0)=\int_\Omega  \Gamma_h(t,x,x_0) v_h(x)\d x ,\quad\forall\, v_h\in S_h ,
\end{align}
and define $|E_h(t)|$ to be the linear operator on $L^q$ with the kernel $|\Gamma_h(t,x,x_0)|$, namely, 
\begin{align}
(|E_h(t)|v) (x_0):=\int_\Omega |\Gamma_h(t,x,x_0)| v(x)\d x ,\quad\forall\, v\in L^q .
\end{align}

The main result of this paper is the following theorem.

\begin{theorem}\label{MainTHM1}
{\it Let $\Omega$ be a polygon in $\R^2$ or a polyhedron in $\R^3$ {\rm(}possibly nonconvex{\rm)}, and let $S_h$, $0<h<h_0$, be a family of finite element subspaces of $H^1_0(\Omega)$ consisting of piecewise polynomials of degree $r\ge 1$ subject to a quasi-uniform triangulation of the domain $\Omega$ (with mesh size $h$). Then we have the following analytic semigroup estimate 
and maximal function estimate: 
\begin{align}
&\sup_{t>0}\,(\|E_h(t)v_h\|_{L^{q}}
+t\|\partial_tE_h(t)v_h\|_{L^{q}})
\leq C\|v_h\|_{L^{q}} , 
&&\forall\,\, v_h\in 
S_h ,
\,\,\,\forall\,\, 1\leq q\leq\infty , \label{analyticity} \\ 
&\left\|\sup_{t>0}|E_h(t)|\, |v| \right\|_{L^{q}}
\leq C_q\|v\|_{L^{q}} , 
&& \forall\,v\in L^q, \,\,\, \forall\, 1<q\le \infty . \label{MEgodic2}
\end{align}

Further, if $u_{h,0}=0$ and $f\in L^p(0,T;L^q)$, then the finite element solution given by \eqref{FEMEq0} possesses  
the following maximal $L^p$-regularity\,{\rm:} 
\begin{align}
&\|\partial_tu_h\|_{L^p(0,T;L^q)}
+\|\Delta_hu_h\|_{L^p(0,T;L^q)} 
\leq \max(p\, ,\, (p-1)^{-1}) C_q\|f \|_{L^p(0,T;L^q)} ,
&& \forall\,\, 1< p,q<\infty, 
\label{MaxLp1}\\
&\|\partial_tu_h\|_{L^\infty(0,T;L^{q})}
+\|\Delta_hu_h\|_{L^\infty(0,T;L^{q})}
\leq C\ell_h\|f\|_{L^\infty(0,T;L^{q})} ,
&&\forall\,\, 1\leq q\leq\infty,
\label{MaxLinfty}  
\end{align} 
where $\ell_h:=\log(2+1/h)$. 

The constant $C$ in \eqref{analyticity} and \eqref{MaxLinfty} is independent of $f$, $h$, $p$, $q$ and $T$, and the constant $C_q$ in \eqref{MEgodic2} and \eqref{MaxLp1} is independent of $f$, $h$, $p$ and $T$. 
}
\end{theorem}

\begin{remark}
By the theory of analytic semigroups \cite[page 254]{Yosida1980}, the inequality \eqref{analyticity} implies  
the existence of a positive constant
$\theta\in(0,\pi/2)$, independent of $h$ and $q$, 
such that the semigroup $\{E_h(t)\}_{t>0}$ extends to be a
bounded analytic semigroup 
$\{E_h(z)\}_{z\in\Sigma_{\theta}}$ 
in the sector $\Sigma_{\theta}:=\{z\in\C\backslash\{0\}:\, |{\rm arg}(z)|<\theta\}$, 
i.e. 
\begin{align}
&E_h(z_1+z_2)=E_h(z_1)E_h(z_2),\quad
\forall\,\, z_1,z_2\in \Sigma_{\theta},\\[5pt]
&\!\! \sup_{z\in \Sigma_{\theta}}\,(\|E_h(z)v_h\|_{L^{q}}
+|z|\|\partial_zE_h(z)v_h\|_{L^{q}})
\leq C\|v_h\|_{L^{q}} ,\quad \forall\,\, v_h\in S_h ,
\,\,\,\forall\,\, 1\leq q\leq\infty . 
\end{align}
\end{remark}


An immediate consequence of Theorem \ref{MainTHM1} is the following $R$-boundedness result for the discrete heat semigroup and discrete resolvent operator, which has important application in deriving the time-discrete maximal $\ell^p$-regularity of the fully discrete finite element solutions discretized with backward Euler, Crank-Nicolson, second-order BDF and A-stable Runge-Kutta schemes (cf. \cite[Section 6]{KovacsLiLubich2016}). 
\begin{corollary}[$R$-boundedness of the discrete resolvent]\label{MainCor0}
{\it 
Under the assumptions of Theorem \ref{MainTHM1}, for any $1<q<\infty$ there exists $\theta_q>0$ {\rm(}independent of $h${\rm)} 
such that 
\begin{enumerate}[label={\rm(\arabic*)},ref=\arabic*]\itemsep=5pt
\item The semigroup of operators
$\{E_h(z): z\in\Sigma_{\theta_q}\}$ is $R$-bounded
in ${\cal L}(L^q,L^q)$ (the space of bounded linear operators on $L^q$), and 
the $R$-bound is independent of $h$.

\item The collection of finite element resolvent operators $\{z(z-\Delta_h)^{-1}:z\in\Sigma_{\frac{\pi}{2}+\theta_q}\}$
is $R$-bounded in ${\cal L}(L^q,L^q)$, and 
the $R$-bound is independent of $h$.  
\end{enumerate}
}
\end{corollary}
\noindent{\it Proof.}$\,\,$
It is easy to see that the maximal semigroup estimate \eqref{MEgodic2} implies the maximal ergodic estimate
\begin{align}
&\bigg\|\sup_{t>0}\frac{1}{t}\int_0^t |E_h(s)| |v| \d s \bigg\|_{L^{q}}
\leq C_q\|v\|_{L^{q}} ,
&& \forall\,\, 1<q\le\infty . \label{MEgodic3}
\end{align}
According to \cite[Lemma 4.c]{Weis2001-1}, for $q\in(1,2]$ the above maximal ergodic estimate implies the $R$-boundedness of the semigroup of operators $\{E_h(z)\}_{z\in\Sigma_{\theta_q}}$ in ${\cal L}(L^q,L^q)$ with $\theta_q=(\theta -\epsilon)q/2$, where $\epsilon$ can be arbitrarily small. For $q\in[2,\infty)$, a duality argument shows that the semigroup $\{E_h(z)\}_{z\in\Sigma_{\theta_q}}$ is $R$-bounded in ${\cal L}(L^q,L^q)$ with angle $\theta_{q}=(\theta-\epsilon)q'/2$ (cf. \cite[Proof of Lemma 4.d]{Weis2001-1}).  

The second statement in Corollary \ref{MainCor0} 
is actually a consequence of the first statement (cf. \cite[Theorem 4.2]{Weis2001-2}).\qed
\medskip 

Recall that the $L^2$ projection $P_h:L^2(\Omega)\rightarrow S_h$ and Ritz projection $R_h:H^1_0(\Omega)\rightarrow S_h$ onto the finite element spaces are defined by 
\begin{align}
&(P_h\phi,\varphi_h)=(\phi,\varphi_h), &&\forall\,\phi\in L^2(\Omega),\,\,\, \forall\,\varphi_h\in S_h, \\
&(\nabla R_h\phi,\varphi_h)=(\nabla \phi,\varphi_h), &&\forall\,\phi\in H^1_0(\Omega),\,\, \forall\,\varphi_h\in S_h .
\end{align}
In particular, the $L^2$ projection actually can be extended to $ L^{q}(\Omega)$, $1\le q\le \infty$, satisfying the following estimate: 
\begin{align}\label{Ph-Lq-Est}
&\|P_h\phi\|_{L^{q}}\le C\|\phi\|_{L^{q}},\quad\forall\, \phi\in L^{q}(\Omega),\,\,\, 1\le q\le \infty ,
\end{align}
where the constant $C$ is independent of the mesh size $h$. The estimate above is a consequence of \cite[Lemma 6.1]{Thomee2006} and the self-adjointness of $P_h$; also see \cite{DouglasDupontWahlbin1974}, \cite[Lemma 7.2]{Wahlbin1991} and the properties of the finite element spaces stated in Section \ref{Sec2-2}.

Besides Corollary \ref{MainCor0}, the maximal $L^p$-regularity results \eqref{MaxLp1}-\eqref{MaxLinfty} also imply the following sharp $L^p(0,T;L^q)$ error estimates for finite element solutions of parabolic equations. 
 
\begin{corollary}\label{MainCor1}
{\it Let 
$u$ and $u_h$ be the solutions of \eqref{PDE1} and \eqref{FEMEq0}, respectively. 
Then, under the assumptions of Theorem \ref{MainTHM1}, we have 
\begin{align*}
&\|u_h-u\|_{L^p(0,T;L^q)}
\leq C_{p,q}(
\|u-R_hu\|_{L^p(0,T;L^q)}+ \|P_hu(0)-u_{h}(0)\|_{L^q}) ,\\
&\|u_h-u\|_{L^\infty(0,T;L^\infty)}
\leq C\ell_h(
\|u-R_hu\|_{L^\infty(0,T;L^\infty)}
+ \|P_hu(0)-u_{h}(0)\|_{L^\infty}) ,
\end{align*}
for $1<p,q<\infty$, 
where $P_h$ and $R_h$ denote the $L^2$-projection and Ritz projection onto the finite element space $S_h$, respectively, and the constants $C_{p,q}$ and $C$ are independent of $u$ and $T$.
}
\end{corollary}
\noindent{\it Proof.}$\,\,\,$
Let $\phi_h:=P_hu-u_h- e^{-t}(P_hu(0)-u_{h}(0))$.
Then $\phi_h$ satisfies the following operator equation:
\begin{align*}
\left\{
\begin{aligned}
&\partial_t\phi_h-\Delta_h\phi_h
=\Delta_h\big(R_hu-P_hu+e^{-t}(P_hu(0)-u_{h}(0)) \big) 
+ e^{-t}(P_hu(0)-u_{h}(0)) , \\
&\phi_h(0)=0 .
\end{aligned}
\right.
\end{align*}
Multiplying the last equation by $\Delta_h^{-1}$, we obtain
\begin{align}\label{EqDinvphi}
\left\{
\begin{aligned}
&\partial_t(\Delta_h^{-1}\phi_h)-\Delta_h(\Delta_h^{-1}\phi_h)
=R_hu-P_hu+(e^{-t}+ e^{-t}\Delta_h^{-1}) (P_hu(0)-u_{h}(0)) ,\\[5pt]
&\Delta_h^{-1}\phi_h(0)=0. 
\end{aligned}
\right.
\end{align}
By applying \eqref{MaxLinfty} to the equation above (with $q=\infty$), we have 
\begin{align}\label{phi_hLinftyLinfty}
\|\phi_h\|_{L^\infty(0,T;L^\infty)}
&=\|\Delta_h(\Delta_h^{-1}\phi_h)\|_{L^\infty(0,T;L^\infty)}  \nonumber \\
&\leq 
C\ell_h \|R_hu-P_hu+(e^{-t}+ e^{-t}\Delta_h^{-1}) (P_hu(0)-u_{h}(0))\|_{L^\infty(0,T;L^\infty)}  \nonumber \\
&\leq C\ell_h(\|R_hu-P_hu\|_{L^\infty(0,T;L^\infty)}
+ \|P_hu(0)-u_{h}(0)\|_{L^\infty}) ,
\end{align}
where we have used the following $L^\infty$ estimate of  finite element solutions of the Poisson equation (a proof is given in Appendix C)
\begin{align}\label{Linfty-Deltah-1}
\|\Delta_h^{-1}(P_hu(0)-u_{h}(0))\|_{L^\infty}
\le C\|P_hu(0)-u_{h}(0)\|_{L^\infty} .
\end{align}
By using the $L^\infty$ stability of the $L^2$ projection (i.e., using \eqref{Ph-Lq-Est} with $q=\infty$), \eqref{phi_hLinftyLinfty} further reduces to 
\begin{align}\label{phi_hLinftyLinfty2}
\|\phi_h\|_{L^\infty(0,T;L^\infty)}
&\leq C\ell_h(\|P_h(R_hu-u)\|_{L^\infty(0,T;L^\infty)}
+ \|P_hu(0)-u_{h}(0)\|_{L^\infty}) \nonumber \\
&\leq C\ell_h(\|R_hu-u\|_{L^\infty(0,T;L^\infty)}
+ \|P_hu(0)-u_{h}(0)\|_{L^\infty}) .
\end{align} 
This proves the second statement of Corollary \ref{MainCor1}. 
The first statement of Corollary \ref{MainCor1} can be proved similarly by applying \eqref{MaxLp1} to \eqref{EqDinvphi}. 
\qed

One of the advantages of Corollary \ref{MainCor1} is that it reduces the $L^\infty$ stability of finite element solutions of parabolic equations to the $L^\infty$ stability of the Ritz projection, which immediately implies the following $L^\infty$ stability results in nonconvex polygons and convex polyhedra.
\begin{corollary}\label{max-norm-stability}
{\it 
Under the assumptions of Theorem \ref{MainTHM1}, 
if $\Omega$ is a polygon in $\R^2$ (possibly nonconvex) or a convex polyhedra in $\R^3$, 
and $u_{h,0}=P_hu_0$ or $u_{h,0}=R_hu_0$, then the solutions of \eqref{PDE1} and \eqref{FEMEq0} satisfy
\begin{align}\label{best-apprx}
\|u-u_h\|_{L^\infty(0,T;L^\infty)}
\leq C\ell_h^2\inf_{\chi_h}\|u-\chi_h\|_{L^\infty(0,T;L^\infty)}  ,
\end{align} 
where the constant $C$ is independent of $h$ and $T$, and the infimum extends over all $\chi_h\in L^\infty(0,T;S_h)$. 
}
\end{corollary}
\noindent{\it Proof.}$\,\,$
In a two-dimensional polygon (possibly nonconvex) or a convex polyhedra, both the $L^2$-projection $P_h$ and the Ritz projection $R_h$ have been proved to be stable in the maximum norm (cf. \cite[Lemma 6.1]{Thomee2006} and \cite{LeykekhmanVexler2016,Schatz1980}), i.e.  
$$
\|u-P_hu\|_{L^\infty}
+\|u-R_hu\|_{L^\infty} \leq C\ell_h\inf_{\chi_h\in S_h}\|u-\chi_h\|_{L^\infty} .
$$
Hence, Corollary \ref{MainCor1} and the inequality above imply \eqref{best-apprx}. \qed\medskip

In the next section, we introduce the notations to be used in this paper. 
The proof of Theorem \ref{MainTHM1} is presented in Section \ref{Proof-Theorem}.

\section{Notations}\label{Sec:notation}
\setcounter{equation}{0}

\subsection{Function spaces}

We use the conventional notations of Sobolev spaces $W^{s,q}(\Omega)$, $s\ge 0$ and $1\leq q\leq\infty$ (cf. \cite{AdamsFournier2003}), 
with the abbreviations $L^q=W^{0,q}(\Omega)$, $W^{s,q}=W^{s,q}(\Omega)$ and $H^s:=W^{s,2}(\Omega)$. 
The notation $H^{-s}(\Omega)$ denotes the dual space of $H^s_0(\Omega)$, the closure of $C^\infty_0(\Omega)$ in $H^s(\Omega)$. 

For any given function $f:(0,T)\rightarrow W^{s,q}$
we define the Bochner norm 
\begin{align}
&\|f\|_{L^p(0,T;W^{s,q})} = 
\big\| \|f(\cdot)\|_{W^{s,q}}\big\|_{L^p(0,T)} ,\quad\forall\,\, 1\leq p,q\leq \infty,\,\, s\in\R   . 
\end{align} 
For any subdomain $D\subset \Omega$, we define 
\begin{align}\label{Def-HsD}
\|f\|_{W^{s,q}(D)}:=
\inf_{\widetilde f|_D=f}\|\widetilde f\|_{W^{s,q}(\Omega)} 
 ,\quad\forall\,\, 1\leq q\leq \infty ,\,\, s\in\R ,
\end{align} 
where the infimum extends over all possible
$\widetilde f$ defined on $\Omega$ such that
$\widetilde f=f$ in $D$. 
Similarly, for any subdomain $Q\subset {\cal Q}=(0,1)\times\Omega$, 
we define 
\begin{align}\label{DefLpX}
\|f\|_{L^pW^{s,q}(Q)}:=
\inf_{\widetilde f|_Q=f}\|\widetilde f\|_{L^p(0,T;W^{s,q})} 
 ,\quad\forall\,\, 1\leq p,q\leq \infty ,\,\,  s\in\R  ,
\end{align} 
where the infimum extends over all possible
$\widetilde f$ defined on ${\cal Q}$ such that
$\widetilde f=f$ in $Q$. 

We use the abbreviations 
\begin{align}\label{inner-products}
(\phi,\varphi):=\int_\Omega \phi(x)\varphi(x)\d x,\qquad
[u,v]:=\int_0^T\int_{\Omega} u(t,x)v(t,x)\d x \, \d t ,
\end{align} 
and denote $w(t)=w(t,\cdot)$ 
for any function $w$ defined on ${\cal Q}$. The notation $1_{0<t<T}$ will denote the characteristic function of the time interval $(0,T)$, i.e. $1_{0<t<T}(t)=1$ if $t\in(0,T)$ while $1_{0<t<T}(t)=0$ if $t\notin(0,T)$.

\subsection{Properties of the finite element space}
\label{Sec2-2}

For any subdomain $D\subset\Omega$, we denote by $S_h(D)$ the space of functions of $S_h$ restricted to the domain $D$, and denote by $S_h^0(D)$ the subspace of $S_h(D)$ consisting of functions which equal zero outside $D$. For any given subset $D\subset\Omega$, we denote $B_d (D) =\{x\in\Omega: {\rm dist}(x,D)\leq d\}$ for $d>0$. On a quasi-uniform triangulation of the domain $\Omega$, there exist positive constants $K $ and $\kappa$ such that the triangulation and the corresponding finite element space $S_h$ possess the following properties ($K$ and $\kappa$ are independent of the subset $D$ and $h$).
\medskip

\begin{enumerate}[label={\bf (P\arabic*)},ref=\arabic*]\itemsep=5pt

\item {\bf Quasi-uniformity:}

For all triangles (or tetrahedron) $\tau_l^h$ in the partition,
the diameter $h_l$ of $\tau_l^h$ and the radius $\rho_l$
of its inscribed ball satisfy
$$
K^{-1}h\leq \rho_l \leq h_l\leq Kh .
$$

\item {\bf Inverse inequality:}

If $D$ is a union of elements in the partition, then
\begin{align*}
\|\chi_h\|_{W^{l,p}(D)}
\leq K h^{-(l-k)-(N/q-N/p)}\|\chi_h\|_{W^{k,q}(D)} ,
\quad\forall\,\,\chi_h\in S_h,
\end{align*}
for $0\leq k\leq l\leq 1$ and 
$1\leq q\leq p\leq\infty$.  

\item {\bf Local approximation and superapproximation:} 

There exists an operator $I_h:H^1_0(\Omega)\rightarrow S_h$ with the following properties 
(cf. Appendix B):

\begin{enumerate}[label=(\arabic*),ref=\arabic*]\itemsep=5pt
\item 
\begin{align*} 
&\|v-I_hv\|_{L^2} +h \|\nabla(v- I_hv)\|_{L^2}  \leq Kh^{1+\alpha} \|v\|_{H^{1+\alpha}} ,  && \forall\,\,v\in H^{1+\alpha}(\Omega )\cap H^1_0(\Omega) ,\,\,\, \alpha\in[0,1] .
\end{align*}

\item 
If $d\geq 2h$ then the value of $I_h v$ in $D$ depends only on the value of $v$ in $B_d(D)$. If $d\ge 2h$ and supp$(v)\subset \overline D$, then $I_hv\in S_h^0(B_{d}(D))$.


\item 
If $d\geq 2h$, $\omega=0$ outside $D$
and $|\partial^\beta\omega|\leq Cd^{-|\beta|}$
for all multi-index $\beta$,
then 
%
\begin{align*}
&\psi_h\in S_h(B_{d}(D))\implies I_h(\omega\psi_h)\in S_h^0(B_{d}(D)),\\
&\|\omega\psi_h-I_h(\omega\psi_h)\|_{L^2}
+h\|\omega\psi_h-I_h(\omega\psi_h)\|_{H^1}
\leq K  h d^{-1}
\|\psi_h\|_{L^2(B_{d}(D))} .
\end{align*}

\item  If $d\ge 2h$ and $\omega\equiv 1$ on $B_{d}(D)$, then $I_h(\omega\psi_h)=\psi_h$ on $D$. 

\end{enumerate}

\end{enumerate}

The properties (P1)-(P3) hold for any quasi-uniform triangulation with the standard finite element spaces consisting of globally continuous piecewise polynomials of degree $r\ge 1$ (cf. \cite[Appendix]{SchatzWahlbin1995}), and have been used in 
many works in studying the discrete maximal $L^p$-regularity and maximum-norm stability of finite element solutions of parabolic equations; see \cite{Geissert2006,Leykekhman2004,Li2015,LiSun2016,SchatzThomeeWahlbin1998,SchatzWahlbin1995,ThomeeWahlbin2000}.  
Property (P3)-(1) and Definition \eqref{Def-HsD} imply the following estimate for $\alpha\in[0,1]$: 
\begin{align}\label{H1+slocal}
&\|v-I_hv\|_{L^2(D)} +h\|v-I_hv\|_{H^1(D)} \leq Kh^{1+\alpha} \|v\|_{H^{1+\alpha}(B_d(D))} ,  && \forall\,\,v\in H^{1+\alpha}(B_d(D))\cap H^1_0(\Omega).
\end{align} 

\subsection{Green's functions}\label{sec:Green}
For any $x_0\in \tau_{l}^h$ (where $\tau_{l}^h$ is a triangle or a tetrahedron in the triangulation of $\Omega$), there exists a function $\widetilde\delta_{x_0}\in C^3(\overline\Omega)$ with support in $\tau_{l}^h$ such that
\begin{align*}
\chi_h(x_0)=\int_{\Omega}\chi_h \widetilde\delta_{x_0}\d x,
\quad\forall\,\chi_h\in S_h ,
\end{align*}
and
\begin{align}
&\|\widetilde\delta_{x_0}\|_{W^{l,p}}
\leq K h^{-l-N(1-1/p)}
\quad\mbox{for}\,\,\,1\leq p\leq\infty,
\,\,\, l=0,1,2,3 , \label{reg-Delta-est}\\
&\sup_{y\in\Omega} \int_\Omega |\widetilde\delta_{y}(x)|\d x+
\sup_{x\in\Omega}\int_\Omega |\widetilde\delta_{y}(x)|\d y \le C .
\label{reg-Delta-est-2}
\end{align}
The construction of $\widetilde\delta_{x_0}$ can be found in \cite[Lemma 2.2]{ThomeeWahlbin2000}. 

Let $\delta_{x_0}$ denote the Dirac Delta function centered at $x_0$. In other words, $\int_\Omega\delta_{x_0}(y)\varphi (y)\d y=\varphi(x_0)$ for arbitrary $\varphi\in C(\overline\Omega)$.  
Then the discrete Delta function 
$$
\delta_{h,x_0}:=P_h  \delta_{x_0}=P_h \widetilde\delta_{x_0} 
$$ 
decays exponentially away from $x_0$ (cf. \cite[Lemma 7.2]{Wahlbin1991}): 
\begin{align}\label{Detal-pointwise}
& |\delta_{h,x_0}(x)|=|P_h \widetilde\delta_{x_0}(x)|
\leq Kh^{-N} e^{-\frac{|x-x_0| }{K h }} ,
\quad \forall\, x,x_0\in\Omega  .
\end{align}

Let $G(t,x,x_0)$ denote the Green's function
of the parabolic equation, i.e. $G=G(\cdot,\cdot\, ,x_0)$ is the solution
of 
\begin{align}\label{GFdef}
\left\{\begin{array}{ll}
\partial_tG(\cdot,\cdot\, ,x_0)-\Delta G(\cdot,\cdot\, ,x_0)=0
&\mbox{in}\,\,\, (0,T)\times \Omega,\\
G(\cdot,\cdot\, ,x_0) = 0
&\mbox{on}\,\,\, (0,T)\times \partial\Omega ,\\
G(0,\cdot,x_0)= \delta_{x_0}
&\mbox{in}\,\,\,\Omega .
\end{array}\right.
\end{align}
The Green's function $G(t,x,y)$ is symmetric with respect to $x$ and $y$. It has an analytic extension to the right half-plane, satisfying the following Gaussian estimate
(cf. \cite[p.\ 103]{Davis1989}): 
\begin{equation}
|G(z,x,y)|\leq C_\theta|z|^{-\frac{N}{2}}
e^{-\frac{|x-y|^2}{C_\theta |z|}}, 
\quad \forall\, z\in \varSigma_{\theta},\,\,\forall\, x,y\in\varOmega ,
\quad\forall\,\theta\in(0,\pi/2), 
\label{GKernelE0}
\end{equation}
where the constant $C_\theta$ depends only on $\theta$. 
Then Cauchy's integral formula says that 
\begin{equation}
\partial_t^kG(t,x,y)
=\frac{k!}{2\pi i}
\int_{|z-t|=\frac{t}{2}}\,\, \frac{G(z,x,y)}{(z-t)^{k+1}}\d z ,
\end{equation}
which further yields the following Gaussian pointwise estimate for the time derivatives of Green's function (cf. \cite[Appendix B with $\alpha=\beta=0$]{Geissert2006}): 
\begin{align}
&|\partial_t^kG(t,x,x_0)|\leq \frac{C_k}{t^{k+N/2}} e^{-\frac{|x-x_0|^2}{C_kt}}, &&
\forall\, x,x_0\in\Omega,\,\,\, \forall\, t>0, \,\, k=0,1,2,\dots \label{GausEst1}
\end{align} 
%

Let $\Gamma=\Gamma(\cdot,\cdot\, , x_0)$ 
be the regularized Green's function
of the parabolic equation, defined by 
\begin{align}\label{RGFdef}
\left\{\begin{array}{ll}
\partial_t\Gamma(\cdot,\cdot\, , x_0)-\Delta\Gamma(\cdot,\cdot\, , x_0)=0
&\mbox{in}\,\,\,(0,T)\times \Omega,\\
\Gamma(\cdot,\cdot\, , x_0) = 0
&\mbox{on}\,\,\, (0,T)\times \partial\Omega ,\\
\Gamma(0,\cdot ,x_0)=\widetilde\delta_{x_0} 
&\mbox{in}\,\,\,\Omega,
\end{array}\right.
\end{align}
and let $\Gamma_h=\Gamma_h(\cdot,\cdot,  x_0)$ 
be the finite element approximation of $\Gamma$,
defined by
\begin{align}\label{EqGammh}
\left\{\begin{array}{ll}
(\partial_t\Gamma_h(t,\cdot,x_0),v_h)+(\nabla\Gamma_h(t,\cdot,x_0),\nabla v_h)=0 ,
&\forall \, v_h\in S_h,\,\, t\in(0,T),\\[5pt]
\Gamma_h(0,\cdot,x_0)=  \delta_{h,x_0} .
\end{array}\right.
\end{align}

By using the Green's function and discrete Green's function, the solutions of \eqref{PDE1} and \eqref{FEMEq0} can be represented by 
\begin{align}
&u(t,x_0)=\int_\Omega  G(t,x,x_0)u_0(x)\d x + \int_0^t\int_\Omega G(t-s,x,x_0)f(s,x)\d x\d s ,\\
&u_h(t,x_0)=\int_\Omega  \Gamma_h(t,x,x_0)u_{0,h}(x)\d x + \int_0^t\int_\Omega \Gamma_h(t-s,x,x_0)f(s,x)\d x\d s , 
\end{align}
and
\begin{align}
&(E(t)v)(x_0)=\int_\Omega  G(t,x,x_0) v(x)\d x  , 
\qquad
(E_h(t)v_h)(x_0)=\int_\Omega  \Gamma_h(t,x,x_0)v_{h}(x)\d x .  \label{semigroup-Green-discr}
\end{align}
The regularized Green's function can be represented by 
\begin{align}\label{expr-Gamma}
&\Gamma(t,x,x_0)=\int_\Omega  G(t,y,x)\widetilde\delta_{x_0}(y)\d y=\int_\Omega  G(t,x,y)\widetilde\delta_{x_0}(y)\d y  . 
\end{align}
From the representation \eqref{expr-Gamma} one can easily derive that the regularized Green's function $\Gamma$ also satisfies the Gaussian pointwise estimate:
\begin{align}
&|\partial_t^k\Gamma(t,x,x_0)|\leq \frac{C_k}{t^{k+N/2}} e^{-\frac{|x-x_0|^2}{C_kt}}, &&
\forall\, x,x_0\in\Omega,\,\,\forall\, t>0\,\,\mbox{such that}\,\, \max(|x-x_0|,\sqrt{t})\ge 2h,  \label{GausEstGamma}
\end{align} 
with $k=0,1,2,\dots$

\subsection{Dyadic decomposition of the domain ${\cal Q}=(0,1)\times\Omega$}
\label{SecGF} 
In the proof of Theorem \ref{MainTHM1}, we need to partition the domain ${\cal Q}=(0,1)\times\Omega$ into subdomains, and present estimates of the finite element solutions in each subdomain. The following dyadic decomposition of ${\cal Q}$ was introduced in \cite{SchatzThomeeWahlbin1998} and has been used by  many authors \cite{Geissert2006,Leykekhman2004,Li2015,LiSun2016,ThomeeWahlbin2000}. The readers may pass this subsection if they are familiar with such dyadic decompositions.

For any integer $j$, we define $d_j=2^{-j}$.  For a given $x_0\in\Omega$, 
we let $J_1=1$, $J_0=0$ and $J_*$ be an integer satisfying $2^{-J_*}= C_*h$ 
with $C_*\geq 16$ to be determined later.  
If 
\begin{align}\label{h-condition}
h<1/(4C_*)  ,
\end{align}
then 
\begin{align}
2\leq J_*=\log_2[1/(C_*h)]\leq \log_2(2+1/h) .
\end{align} 
Let  
\begin{align*}
&Q_*(x_0)=\{(x,t)\in\Omega_T: \max (|x-x_0|,t^{1/2})\leq d_{J_*}\}, 
\\
&\Omega_*(x_0)=\{x\in \Omega: |x-x_0|\leq d_{J_*}\} \, . 
\end{align*} 
We define 
\begin{align*} 
&Q_j(x_0)=\{(x,t)\in \Omega_T:d_j\leq
\max (|x-x_0|,t^{1/2})\leq2d_j\}  &&\mbox{for}\,\,\, j\ge 1,
\\
&\Omega_j(x_0)=\{x\in \Omega: d_j\leq|x-x_0|\leq2d_j\}   
&&\mbox{for}\,\,\, j\ge 1,
\\
&D_j(x_0)=\{x\in \Omega:  |x-x_0|\leq2d_j\}  
&&\mbox{for}\,\,\, j\ge 1,
\end{align*} 
and 
\begin{align*} 
&Q_0(x_0)= 
{\cal Q}\big\backslash\big( \cup_{j=1}^{J_*}Q_{j}(x_0)\cup Q_*(x_0)\big) , 
\\
&\Omega_0(x_0) 
=\Omega\big\backslash\big( \cup_{j=1}^{J_*}\Omega_{j}(x_0)\cup \Omega_*(x_0)\big).
\end{align*}
For $j<0$, we simply define
$Q_{j}(x_0)=\Omega_{j}(x_0)=\emptyset$.
For all integer $j\ge 0$, we define
\begin{align*}
&\Omega_j'(x_0)=\Omega_{j-1}(x_0)\cup\Omega_{j}(x_0)\cup\Omega_{j+1}(x_0),
\quad Q_j'(x_0)=Q_{j-1}(x_0)\cup Q_{j}(x_0)\cup Q_{j+1}(x_0), \\
&\Omega_j''(x_0)=\Omega_{j-2}(x_0)\cup\Omega_{j}'(x_0)\cup\Omega_{j+2}(x_0),
\quad
Q_j''(x_0)=Q_{j-2}(x_0)\cup Q_{j}'(x_0)\cup Q_{j+2}(x_0),\\
&D_j'(x_0)=D_{j-1}(x_0)\cup D_{j}(x_0) ,
\quad D_j''(x_0)=D_{j-2}(x_0)\cup D_{j}'(x_0) .
\end{align*}
Then we have
\begin{align}\label{decomposition}
&\Omega_T=\bigcup^{J_*}_{j=0}Q_j(x_0)\,\cup Q_*(x_0)
\quad\mbox{and}\quad
\Omega=\bigcup^{J_*}_{j=0}\Omega_j(x_0)\,\cup \Omega_*(x_0) .
\end{align}
We refer to $Q_*(x_0)$ as the ``innermost" set.
We shall write $\sum_{*,j}$ when the innermost set is included and
$\sum_j$ when it is not. When $x_0$ is fixed, if there is no ambiguity, 
we simply write
$Q_j=Q_j(x_0)$, $Q_j'=Q_j'(x_0)$, $Q_j''=Q_j''(x_0)$, 
$\Omega_j=\Omega_j(x_0)$, $\Omega_j'=\Omega_j'(x_0)$ and 
$\Omega_j''=\Omega_j''(x_0)$.

We shall use the notations 
\begin{align}\label{Q-norm}
&\|v\|_{k,D}=\biggl(\int_D\sum_{|\alpha|\leq
k}|\partial^\alpha v|^2\d x\biggl)^{\frac{1}{2}},\qquad
\vertiii{v}_{k,Q}=\biggl(\int_Q\sum_{|\alpha|\leq
k}|\partial^\alpha v|^2\d x\d t\biggl)^{\frac{1}{2}} ,
\end{align}
for any subdomains $D\subset\Omega$ and $Q\subset (0,1)\times\Omega $.
Throughout this paper, we denote by $C$ a generic positive constant that 
is independent of $h$, $x_0$ and $C_*$ (until $C_*$ is determined 
in Section \ref{ProofLm03}). To simplify the notations, we also denote $d_*=d_{J_*}$.

\section{Proof of Theorem \ref{MainTHM1}}\label{Proof-Theorem}
\setcounter{equation}{0}

\subsection{Estimates of the Green's function}\label{PrfGFE}

In this subsection, we prove the following local $L^2H^{1+\alpha}$ and $L^\infty H^{1+\alpha}$ estimates for the Green's function and regularized Green's function. These local estimates are needed in the proof of Theorem \ref{MainTHM1}.  
\begin{lemma}\label{GFEst1}
{\it Let $\Omega$ be a polygon in $\R^2$ or a polyhedron in $\R^3$ $($possibly nonconvex$)$. There exists $\alpha\in (\frac{1}{2},1]$ and $C>0$, 
independent of $h$ and $x_0$, such that 
the Green's function $G$ defined in 
\eqref{GFdef} and the regularized Green's function 
$\Gamma$ defined in \eqref{RGFdef}
satisfy the following estimates:
\begin{align}
&d_j^{-4-\alpha+N/2}\|\Gamma(\cdot,\cdot,x_0)\|_{L^\infty(Q_j(x_0))}
+d_j^{-4-\alpha}\vertiii{\nabla\Gamma(\cdot,\cdot,x_0)}_{L^2(Q_j(x_0))}\nn\\
&
+d_j^{-4}\vertiii{ \Gamma(\cdot,\cdot,x_0)}_{L^2H^{1+\alpha}(Q_j(x_0))}  +d_j^{-2}\vertiii{ \partial_{t} 
\Gamma(\cdot,\cdot,x_0)}_{L^2H^{1+\alpha}(Q_j(x_0))} \nn\\
& 
+\vertiii{ \partial_{tt}\Gamma(\cdot,\cdot,x_0)}_{L^2H^{1+\alpha}(Q_j(x_0))}
\leq Cd_j^{-N/2-4-\alpha}, \label{GFest01}\\[10pt]
&\|G(\cdot,\cdot,x_0) \|_{L^{\infty}H^{1+\alpha}(\cup_{k\leq
j}Q_k(x_0))}
+d_j^2\|\partial_tG(\cdot,\cdot,x_0) \|_{L^{\infty}H^{1+\alpha}(\cup_{k\leq
j}Q_k(x_0))}\leq Cd_j^{-N/2-1-\alpha} \label{GFest03} .
\end{align}
}
\end{lemma}

To prove Lemma \ref{GFEst1}, we need to use 
the following two lemmas. 
\begin{lemma}\label{RegPoiss}
{\it Let $\Omega$ be a polygon in $\R^2$ or a polyhedron in $\R^3$ $($possibly nonconvex$)$. Then there exists a
positive constant $\alpha\in(\frac{1}{2},1]$ $($depending on the domain $\Omega)$ such that the solution
of the Poisson equation
\begin{align*}
\left\{\begin{array}{ll}
\Delta \varphi=f
&\mbox{in}\,\,\,\Omega, \\
\varphi=0
&\mbox{on}\,\,\,\partial\Omega,
\end{array}\right.  
\end{align*}
satisfies 
\begin{align*} 
\|\varphi\|_{H^{1+\alpha}} 
\leq C\|f\|_{H^{-1+\alpha}} . 
\end{align*} 
}
\end{lemma}

Lemma \ref{RegPoiss} is a consequence of \cite[Theorem 18.13]{Dauge1988}, where either ``$n_x=2$ ($x$ is an edge point) and $\mu_1(\Gamma_x)>1/2$ (first eigenvalue of the Laplacian in a 2D sector)'' or ``$n_x=3$ ($x$ is a vertex) and $\mu_1(\Gamma_x)>0$ (first eigenvalue of the Laplacian in a 3D cone)". 
Such regularity also holds for the Neumann Laplacian (cf. \cite[Corollary 23.5]{Dauge1988}).

\begin{lemma}\label{FracLapl}
{\it Let $\Omega$ be a polygon in $\R^2$ or a polyhedron in $\R^3$ $($possibly nonconvex$)$. 
Then there exists $\alpha\in(\frac{1}{2},1]$ such that 
\begin{align}
\|u\|_{H^{1+\alpha}} 
\leq  C\|\nabla u\|_{L^2}^{1-\alpha} \|\Delta u\|_{L^2}^\alpha , 
\quad\forall\, \mbox{$u\in H^1_0$\, such that \,$\Delta u\in L^2$.}
\label{FracLpl2} 
\end{align}
}
\end{lemma}
\noindent{\it Proof.}$\,\,\,$ 
Let $\phi_j$, $j=1,2,\dots$,
be the orthornormal eigenfunctions of the 
Dirichlet Laplacian $-\Delta$, 
corresponding to the eigenvalues $\lambda_j>0$, 
$j=1,2,\dots$, respectively. With these notations, 
if $u=\sum_{j=1}^\infty a_j\phi_j$ then 
$$
(-\Delta)^{\alpha/2}u := \sum_{j=1}^\infty \lambda_j^{\alpha/2} a_j\phi_j ,\qquad \alpha\in [0,2] .
$$
The norm {\color{blue}$\|(-\Delta)^{\alpha/2}u\|_{L^2}:=\left(\sum_{j=1}^\infty\lambda_j^\alpha a_j^2\right)^{\frac12}$} can be viewed as the a weighted $\ell^2$ norm of the sequence $(a_1,a_2,\cdots)$. 
Since 
$\|(-\Delta)^0 v\|_{L^2}=\|v\|_{L^2}$
and 
$C^{-1}\|v\|_{H^{1}}\leq 
\|(-\Delta)^{1/2}v\|_{L^2}\leq C \|v\|_{H^{1}}$, 
the complex interpolation method (cf. \cite[Theorem 5.4.1]{BerghLofstrom1976}) yields the following equivalence of norms: 
$$C^{-1}\|v\|_{H^{1-\alpha}(\Omega)}\leq 
\|(-\Delta)^{(1-\alpha)/2}v\|_{L^2(\Omega)} 
\leq C\|v\|_{H^{1-\alpha}(\Omega)} ,
\qquad 
\forall\, v\in H^{1-\alpha}(\Omega),\,\,\,
\forall\,\, \mbox{$\alpha\in(\frac{1}{2},1] $}. 
$$
Hence, we have  
\begin{align*}
|(-\Delta u,v)|
&=|((-\Delta)^{(1+\alpha)/2} u,(-\Delta)^{(1-\alpha)/2}v)| \\ 
&\leq \|(-\Delta)^{(1+\alpha)/2} u\|_{L^2} 
\|(-\Delta)^{(1-\alpha)/2}v\|_{L^2} \\
&\leq C\|(-\Delta)^{(1+\alpha)/2} u\|_{L^2} 
\|v\|_{H^{1-\alpha}} ,
&&\forall\, v\in C_0^\infty (\Omega) ,
\end{align*}
which implies (via the duality argument)
\begin{align*}
\|\Delta u\|_{H^{-1+\alpha}}^2
&\le C\|(-\Delta)^{(1+\alpha)/2} u\|_{L^2}^2 \\
&=\sum_j \lambda_j^{1+\alpha}a_j^2 \\
&=\sum_j (\lambda_ja_j^2)^{1-\alpha}
(\lambda_j^2a_j^2)^{\alpha} \\
&\leq \bigg(\sum_j \lambda_j a_j^2 \bigg)^{1-\alpha}
\bigg(\sum_j\lambda_j^2a_j^2\bigg)^{\alpha} 
\qquad\qquad\mbox{(H\"older's inequality)} \\
&=\|(-\Delta)^{1/2} u\|_{L^2}^{2(1-\alpha)}
\|\Delta  u\|_{L^2}^{2\alpha} \\
&=\|\nabla u\|_{L^2}^{2(1-\alpha)}
\|\Delta  u\|_{L^2}^{2\alpha}  .
\end{align*}
Lemma \ref{RegPoiss} implies the existence of $\alpha\in(\frac{1}{2},1]$ such that $\|u\|_{H^{1+\alpha}}\le C\|\Delta u\|_{H^{-1+\alpha}} $ for all $u\in H^1_0$ such that $\Delta u\in L^2$. This (together with the inequality above) yields Lemma \ref{FracLapl}. \qed\medskip

\noindent{\it Proof of Lemma \ref{GFEst1}.}$\,\,$ 
To simplify the notations, we denote $Q_j=Q_j(x_0)$. 
Let $0\leq \omega_j(x,t)\leq 1$ and $0\leq \widetilde\omega_j(x,t)\leq 1$ 
be smooth cut-off functions vanishing outside $Q_j'$ and equals $1$ in $Q_{j}$, such that 
$\widetilde\omega_{j}=1$ on the support of $\omega_{j}$,
and 
\begin{align}\label{omega-dj}
|\partial_t^{k_1}\nabla^{k_2}\omega_j|+|\partial_t^{k_1}\nabla^{k_2}\widetilde\omega_j|
\leq Cd_j^{-2k_1-k_2}
\end{align} 
for all nonnegative integers $k_1$ and $k_2$.
By the definition in \eqref{DefLpX}, it suffices to prove the corresponding global estimates for the function $\omega_1G$, which equals $G$ in $Q_j$.

Consider $\omega_jG$, which is the solution of 
\begin{align}\label{omegaG1}
\partial_{t}(\omega_jG)
-\Delta (\omega_jG)
= \widetilde\omega_jG\partial_{t}\omega_j
+\widetilde\omega_jG\Delta \omega_j 
-\nabla\cdot(2
\widetilde\omega_j G\nabla \omega_j) 
\end{align}
in the domain $(0,\infty)\times \Omega$, 
with zero boundary and initial conditions.
The standard energy estimate yields (cf. \cite[Lemma 2.1 of Chapter III]{LSU1988}, with $q_1=r_1=(2N+4)/(N+4)$) 
\begin{align}\label{EngomegaG1}
&\|\omega_jG\|_{L^\infty(0,T;L^2)}
+\|\omega_jG\|_{L^2(0,T;H^1)} \nn\\
&\leq C\|\widetilde\omega_jG 
\partial_{t}\omega_j\|_{L^{(2N+4)/(N+4)}({\cal Q})}
+ C\|\widetilde\omega_jG\Delta\omega_j 
\|_{L^{(2N+4)/(N+4)}({\cal Q})}
+C\|2\widetilde\omega_jG\nabla\omega_j
\|_{L^2(0,T;L^2)} \nn\\
&\leq Cd_j^{-2}\|G\|_{L^{(2N+4)/(N+4)}(Q_j')}
+Cd_j^{-1}\|G\|_{L^2(Q_j')} \nn\\
&\leq Cd_j^{-N/2}  ,
\end{align}
where we have used the Gaussian estimate \eqref{GausEst1} in the last step. 
The last inequality implies 
\begin{align}\label{GL2H1Qj}
&\|G\|_{L^\infty L^2(Q_j)}
+\|G\|_{L^2H^1(Q_j)} \leq Cd_j^{-N/2} ,
\end{align}
and 
\begin{align}\label{aakk}
\|\nabla\cdot(2\widetilde\omega_jG\nabla \omega_j )
\|_{L^2({\cal Q})}
&\le 
C\|G\Delta \omega_j\|_{L^2({\cal Q})}
+
C\|\nabla G \cdot \nabla \omega_j\|_{L^2({\cal Q})} \nn \\
&\le 
Cd_j^{-2}\|G\|_{L^2(Q_j')}
+
Cd_j^{-1}\|\nabla G\|_{L^2(Q_j')} \nn \\
&\le 
Cd_j^{-1}\|G\|_{L^\infty L^2(Q_j')}
+
Cd_j^{-1}\|\nabla G\|_{L^2(Q_j')} \nn \\
&\leq Cd_j^{-1-N/2} ,
\end{align}
where we have used \eqref{GL2H1Qj} in the last inequality (replacing $Q_j$ by $Q_j'$).
By applying the energy estimate to \eqref{omegaG1}, we have 
\begin{align}\label{EngomegaG2}
&\|\partial_t(\omega_jG)\|_{L^2(0,T;L^2)}
+\|\Delta(\omega_jG)\|_{L^2(0,T;L^2)} \nn\\
&\leq C\|\widetilde\omega_jG 
\partial_{t}\omega_j\|_{L^{2}({\cal Q})}
+ C\|\widetilde\omega_jG\Delta \omega_j 
\|_{L^{2}({\cal Q})}
+C\|\nabla\cdot(2\widetilde\omega_jG\nabla \omega_j)
\|_{L^2({\cal Q})} \nn\\
&\leq Cd_j^{-1-N/2} ,  
\end{align}
where we have used \eqref{GausEst1}, \eqref{omega-dj} and \eqref{aakk} in the last step.

Lemma \ref{FracLapl} implies the existence of $\alpha\in(\frac{1}{2},1]$ (depending on the domain $\Omega$)  such that 
\begin{align}\label{GH1s1}
\|\omega_jG\|_{L^2(0,T;H^{1+\alpha})} 
&\le C\|\nabla(\omega_jG)\|_{L^2(0,T;L^2)}^{1-\alpha}  \|\Delta(\omega_jG)\|_{L^2(0,T;L^2)}^{\alpha}  \nn  \\ 
&\le Cd_j^{-\alpha-N/2} ,
\end{align}
where we have used \eqref{EngomegaG1} and \eqref{EngomegaG2} in the last step. 

Similarly (replacing $G$ by $\partial_tG$
and $\partial_{tt}G$ in the estimates above), one can prove 
the following estimates: 
\begin{align}
&d_j^{-\alpha} \|  \nabla (\omega_j\partial_t G)\|_{L^2(0,T;L^2)} 
+d_j^{1-\alpha} \|  \Delta (\omega_j\partial_t G)\|_{L^2(0,T;L^2)} 
+\|\omega_j \partial_t G\|_{L^2(0,T;H^{1+\alpha})} 
\leq Cd_j^{-\alpha-2-N/2} ,  \label{GtH1s1} \\
&d_j^{-\alpha} \|  \nabla (\omega_j\partial_{tt}G)\|_{L^2(0,T;L^2)} 
+d_j^{1-\alpha} \|  \Delta (\omega_j\partial_{tt} G)\|_{L^2(0,T;L^2)} 
+\|\omega_j  \partial_{tt} G \|_{L^2(0,T;H^{1+\alpha})} 
\leq Cd_j^{-\alpha-4-N/2} .  \label{GttH1s1}
\end{align}
By using \eqref{GausEst1} and \eqref{omega-dj}, the last two inequalities imply
\begin{align}
\|\partial_t(\omega_j G)\|_{L^2(0,T;H^{1+\alpha})} 
&\le \|\partial_t \omega_jG\|_{L^2(0,T;H^{1+\alpha})} 
+ \|\omega_j \partial_t  G\|_{L^2(0,T;H^{1+\alpha})} 
\leq Cd_j^{-\alpha-2-N/2} ,
\label{GH1s2}\\[5pt]
\|\partial_{tt}(\omega_j G)\|_{L^2(0,T;H^{1+\alpha})} 
&\le  \|\partial_{tt}\omega_j G\|_{L^2(0,T;H^{1+\alpha})} 
+ 2\|\partial_t \omega_j \partial_t G\|_{L^2(0,T;H^{1+\alpha})}  
+ \|\omega_j\partial_{tt}  G\|_{L^2(0,T;H^{1+\alpha})} \nn \\
&
\leq Cd_j^{-\alpha-4-N/2} .
\label{GH1s3}
\end{align}
The estimates \eqref{GL2H1Qj} and \eqref{GH1s1}-\eqref{GH1s3} imply 
\begin{align}\label{TrueGFest}
&d_j^{-4-\alpha+N/2}\|G(\cdot,\cdot,x_0)\|_{L^\infty(Q_j(x_0))}
+d_j^{-4-\alpha}\vertiii{\nabla G(\cdot,\cdot,x_0)}_{L^2(Q_j(x_0))}\nn\\
&
+d_j^{-4}\vertiii{G(\cdot,\cdot,x_0)}_{L^2H^{1+\alpha}(Q_j(x_0))}  +d_j^{-2}\vertiii{ \partial_{t} 
G(\cdot,\cdot,x_0)}_{L^2H^{1+\alpha}(Q_j(x_0))} \nn\\
& 
+\vertiii{ \partial_{tt}G(\cdot,\cdot,x_0)}_{L^2H^{1+\alpha}(Q_j(x_0))}
\leq Cd_j^{-N/2-4-\alpha} . 
\end{align}
The estimate \eqref{TrueGFest} can also be proved for the regularized Green's function $\Gamma$ 
by using the following expression:
\begin{align}
\Gamma(t,x,x_0)
=\int_{\tau^h_l} G(t,y,x)\widetilde\delta_{x_0}(y)\d y 
=\int_{\tau^h_l} G(t,x,y)\widetilde\delta_{x_0}(y)\d y ,
\end{align}
where $\tau^h_l$ is the triangle/tetrahedron containing $x_0$ (thus $\widetilde\delta_{x_0}$ is supported in $\tau^h_l$). 
For example, if $y\in\tau^h_l$ then 
$(t,x)\in Q_j(x_0)\Rightarrow (t,x)\in Q_j'(y)$, which implies 
\begin{align*}
\|\Gamma(\cdot,\cdot , x_0)\|_{L^2 H^{1+\alpha}(Q_j(x_0))} 
&=\bigg\| \int_{\tau^h_l}  G(\cdot,\cdot ,y) 
 \widetilde\delta_{x_0}(y) \d y\bigg\|_{L^2 H^{1+\alpha}(Q_j(x_0))}  \\ 
&\le \int_{\tau^h_l} \|G(\cdot,\cdot ,y)\|_{L^2 H^{1+\alpha}(Q_j(x_0))} 
|\widetilde\delta_{x_0}(y)|\d y \\
&\le \int_{\tau^h_l} \|G(\cdot,\cdot ,y)\|_{L^2 H^{1+\alpha}(Q_j'(y))} 
|\widetilde\delta_{x_0}(y)|\d y \\
&\leq \int_\Omega Cd_j^{-\alpha-N/2}|\widetilde\delta_{x_0}(y)|\d y \\
&\leq Cd_j^{-\alpha-N/2} .
\end{align*}
This completes the proof of \eqref{GFest01}.

From \eqref{EngomegaG1} and \eqref{EngomegaG2} we see that 
\begin{align}
&\|\nabla G\|_{L^2(\cup_{k\le j}Q_k')}
\le \sum_{k\le j} \|\nabla  G\|_{L^2(Q_k')} 
\le C\sum_{k\le j} d_k^{-N/2} 
\le Cd_j^{-N/2}  , \label{GH1Qkj} \\
&\|\Delta G\|_{L^2(\cup_{k\le j}Q_k')}
\le \sum_{k\le j} \|\Delta G\|_{L^2(Q_k')} 
\le C\sum_{k\le j} d_k^{-1-N/2} 
\le Cd_j^{-1-N/2}  ,   \label{GH2Qkj} 
\end{align}
Let $\chi_j$ be a smooth cut-off function which equals $1$ on $\cup_{k\le j}Q_k$ and 
equals zero outside $\cup_{k\le j}Q_k'$, satisfying $|\partial_t^l\nabla^m\chi_j|\le Cd_j^{-2l-m}$ for all nonnegative integers $l$ and $m$. Then 
$\chi_jG$ is a function defined on ${\cal Q}$ and equals $G$ on $\cup_{k\le j}Q_k$. The inequalities  
\eqref{GausEst1} and \eqref{GH1Qkj}-\eqref{GH2Qkj} imply 
\begin{align} 
\|\nabla(\chi_j G)\|_{L^2(0,T;L^2)}\le Cd_j^{-N/2},\qquad
\|\Delta(\chi_j G)\|_{L^2(0,T;L^2)}\le Cd_j^{-1-N/2} .
\end{align}
Then Lemma \ref{FracLapl} implies  
\begin{align} 
\|\chi_j G\|_{L^2(0,T;H^{1+\alpha})} 
&\le C\|\nabla(\chi_j G)\|_{L^2(0,T;L^2)}^{1-\alpha} \|\Delta(\chi_j G)\|_{L^2(0,T;L^2)}^\alpha \le Cd_j^{-\alpha-N/2} ,
\end{align}
Similarly one can prove (by using \eqref{GtH1s1}-\eqref{GttH1s1})
\begin{align} 
\|\partial_t(\chi_j G)\|_{L^2(0,T;H^{1+\alpha})} 
&\le C\|\nabla\partial_t(\chi_j G)\|_{L^2(0,T;L^2)}^{1-\alpha} \|\Delta\partial_t(\chi_j G)\|_{L^2(0,T;L^2)}^\alpha \le Cd_j^{-\alpha-2-N/2} ,\\
\|\partial_{tt}(\chi_j G)\|_{L^2(0,T;H^{1+\alpha})} 
&\le C\|\nabla\partial_{tt}(\chi_j G)\|_{L^2(0,T;L^2)}^{1-\alpha} \|\Delta\partial_{tt}(\chi_j G)\|_{L^2(0,T;L^2)}^\alpha \le Cd_j^{-\alpha-4-N/2} .
\end{align}
Hence, the interpolation between the last two inequalities yield  
\begin{align}
&\|\chi_jG\|_{L^\infty(0,T;H^{1+\alpha})} 
\le \|\chi_jG\|_{L^2(0,T;H^{1+\alpha})}^{1/2}
\|\partial_t(\chi_jG)\|_{L^2(0,T;H^{1+\alpha})} ^{1/2}
\leq Cd_j^{-\alpha-1-N/2} , \label{GH1s4} \\
&\|\partial_t(\chi_jG)\|_{L^\infty(0,T;H^{1+\alpha})} 
\le \|\partial_t(\chi_jG)\|_{L^2(0,T;H^{1+\alpha})}^{1/2}
\|\partial_{tt}(\chi_j G)\|_{L^2(0,T;H^{1+\alpha})} ^{1/2}
\leq Cd_j^{-\alpha-3-N/2} .
\end{align}
This completes the proof of \eqref{GFest03}. 
\qed

Besides Lemma \ref{GFEst1}, we also need the following lemma in the proof of Theorem \ref{MainTHM1}. 
The proof of this lemma is deferred to Section \ref{ProofLm03}.  
\begin{lemma}\label{LemGm2}
{\it
Under the assumptions of Theorem \ref{MainTHM1}, the functions $\Gamma_h(t,x,x_0)$, $\Gamma(t,x,x_0)$ and 
$F(t,x,x_0):=\Gamma_h(t,x,x_0)-\Gamma(t,x,x_0)$ satisfy 
\begin{align}
&\sup_{t\in(0,\infty)}\,(\|\Gamma_h(t,\cdot, x_0)\|_{L^1(\Omega)}  
+t\|\partial_{t}\Gamma_h(t,\cdot, x_0)\|_{L^1(\Omega)} ) \leq C , 
\label{L1Gammh}\\ 
&\sup_{t\in(0,\infty)}\,(\|\Gamma(t,\cdot, x_0)\|_{L^1(\Omega)}  
+t\|\partial_{t}\Gamma(t,\cdot, x_0)\|_{L^1(\Omega)} ) \leq C , 
\label{L1Gammh-2}\\ 
&\|\partial_tF(\cdot,\cdot ,x_0)\|_{L^1((0,\infty)\times\Omega)}  
+\|t\partial_{tt}F(\cdot,\cdot , x_0)\|_{L^1((0,\infty)\times\Omega)} \leq C ,
\label{L1Ft}\\
&
\|\partial_t\Gamma_h(t,\cdot, x_0)\|_{L^1}\leq Ce^{-\lambda_0t} , \qquad\forall\, t\ge 1, \label{L1Gammatx0}
\end{align} 
where the constants $C$ and $\lambda_0$ are independent of $h$. 
}
\end{lemma}

\subsection{Proof of (\ref{analyticity})-(\ref{MEgodic2})}\label{sec:proof2.3}


By denoting 
$$
\Gamma_h(t)=\Gamma_h(t,\cdot,  x_0),\qquad
\Gamma(t)=\Gamma(t,\cdot,  x_0)\quad\mbox{and}\quad
F(t)=\Gamma_h(t)-\Gamma(t),
$$ 
and using the Green's function representation \eqref{semigroup-Green-discr}, we have 
\begin{align*} 
\qquad\quad\,\,\,
(E_h(t)v_h)(x_0) 
=(\Gamma_h(t),v_h) 
&= (F (t), v_h) + (\Gamma(t), v_h) \\
&
=\int_0^t(\partial_sF (s), v_h)\d s+(F(0), v_h) + (\Gamma(t), v_h)
\end{align*}
and
\begin{align*}
(t\partial_tE_h(t)v_h)(x_0)
&= (t\partial_tF (t), v_h) + (t\partial_t\Gamma(t), v_h) \\
&=\int_0^t(s\partial_{ss}F (s)+\partial_{s}F (s), v_h)\d s
+ (t\partial_t\Gamma(t), v_h) ,
\end{align*}
with (cf. \eqref{reg-Delta-est}-\eqref{Detal-pointwise} and Lemma \ref{LemGm2}) 
\begin{align*} 
\|F(0)\|_{L^1}+\|\Gamma(t)\|_{L^1}+\|t\partial_t\Gamma(t)\|_{L^1} 
&\le C(\|\delta_{h,x_0}-\widetilde\delta_{x_0}\|_{L^1}+\|\Gamma(t)\|_{L^1}+\|t\partial_t\Gamma(t)\|_{L^1} )
\le C .
\end{align*}
By applying Lemma \ref{LemGm2} to the last two equations, we obtain 
\begin{align*} 
&|(E_h(t)v_h)(x_0)|
\le \big(\|\partial_tF\|_{L^1((0,\infty)\times\Omega)}   
+\|F(0)\|_{L^1} +\|\Gamma(t)\|_{L^1}\big) \|v_h\|_{L^\infty} \le C\|v_h\|_{L^\infty},\\
&|(t\partial_tE_h(t)v_h)(x_0)|
\le \big(\|t\partial_{tt}F\|_{L^1((0,\infty)\times\Omega)}   
+\|\partial_tF\|_{L^1((0,\infty)\times\Omega)}   
 +\|t\partial_t\Gamma(t) \|_{L^1}\big) \|v_h\|_{L^\infty} \le C\|v_h\|_{L^\infty} .
\end{align*}
This proves \eqref{analyticity} in the case $q=\infty$. The case $2\le q\le \infty$ follows from the two end-point cases 
$q=2$ and $q=\infty$ via interpolation, and the case $1\le q\le 2$ follows from the case $2\le q\le \infty$ via duality (the operators $E_h(t)$ and $\partial_tE_h(t)$ are self-adjoint). 
The proof of \eqref{analyticity} is complete. \medskip

In order to prove (\ref{MEgodic2}), we need to construct a symmetrically truncated Green's function 
(since the regularized Green's function $\Gamma(t,x,x_0)$ may not be symmetric with respect to $x$ and $x_0$). In fact, there exists a truncated Green's function $G_{\rm tr}^*(t,x,y)$ satisfying the following conditions (cf. \cite{Li2015,LiSun2016}): 

(1) $G_{\rm tr}^*(t,x,y)$ is symmetric with respect to $x$ and $y$, namely, $G_{\rm tr}^*(t,x,y)=G_{\rm tr}^*(t,y,x)$. 

(2) $G_{\rm tr}^{*}(\cdot,\cdot , y)= 0$ in $Q_{*}(y):=\{(t,x)\in {\cal Q}: \max(|x-y|,\sqrt{t})<d_*\}$, and $G_{\rm tr}^{*}(0,\cdot, y)\equiv 0$ in $\Omega$.  

(3) $0\leq G^*_{\rm tr}(t,x,y)\leq G(t,x,y)$ and 
$G^*_{\rm tr}(t,x,y)= G(t,x,y)$ when
$\max(|x-y|,\sqrt{t})>2d_*$,

(4) $|\partial_tG_{\rm tr}^*(t,x,y)|\leq Cd_{*}^{-N-2}$
when $\max(|x-y|,t^{1/2}) \leq 2d_*$.

Note that for the fixed triangle/tetrahedron $\tau_l^h$ and
the point $x_0\in \tau_l^h$, the function
$\widetilde \delta_{x_0}$ is supported in
$\tau_l^h\subset \Omega_*(x_0)$ with
$\int_\Omega\widetilde\delta_{x_0}(y)\d y=1$. 
By using Lemma \ref{GFEst1}, there exists $\alpha\in(\frac{1}{2},1]$
such that (with $Q_{2*}(y):=\{(t,x)\in {\cal Q}: \max(|x-y|,\sqrt{t})<2d_*\}$)
\begin{align}\label{alkdkdkd}
&\iint_{[(0,\infty)\times\Omega ]\backslash
Q_{2*}(x_0)}|\partial_t\Gamma(t,x,x_0)
-\partial_tG^*_{\rm tr}(\tau,x,x_0)|\d x\d t \nn \\
&=\iint_{[(0,\infty)\times\Omega]\backslash
Q_{2*}(x_0)}|\partial_t\Gamma(t,x,x_0)
-\partial_tG (t,x,x_0)|\d x\d t  \nn \\
&\leq \iint_{[(0,1)\times\Omega]\backslash Q_{2*}(x_0)}\biggl|\int_{\tau^h_l}
\partial_tG(t,x,y)\widetilde\delta_{x_0}(y)\d y 
- \int_{\tau^h_l} \partial_tG(t,x,x_0)\widetilde \delta_{x_0}(y)\d y\biggl|\d x\d t  \nn \\
&~~~ 
+\iint_{(1,\infty)\times\Omega}|\partial_t\Gamma(t,x,x_0)
-\partial_tG(t,x,x_0)|\d x\d t \nn \\
&\leq 
C\iint_{{\cal Q}\backslash Q_{2*}(x_0)} h^{\alpha-(N-2)/2}| \partial_tG(t,x,\cdot)|_{C^{\alpha-(N-2)/2}(\overline \tau_l^h)}\d x\d t \nn \\
&\quad + \iint_{(1,\infty)\times\Omega}|\partial_t\Gamma(t,x,x_0)
-\partial_tG(t,x,x_0)|\d x\d t \nn \\
&=: {\cal I}_1+{\cal I}_2. 
\end{align}
By using \eqref{GausEst1} and \eqref{GausEstGamma} we have 
\begin{align*}
{\cal I}_2=\iint_{(1,\infty)\times\Omega}|\partial_t\Gamma(t,x,x_0)
-\partial_tG(t,x,x_0)|\d x\d t 
\le C\iint_{(1,\infty)\times\Omega }
\frac{C}{t^{1+N/2}} e^{-\frac{|x-x_0|^2}{Ct}} \d x\d t 
\le C .
\end{align*}
For $(t,x)\in Q_j(x_0)$ and $y\in \overline\tau^h_l$, we have $(t,y)\in Q_j'(x)$, which implies that 
\begin{align*}
{\cal I}_1
&\le C\sum_{j}\iint_{Q_j(x_0)}h^{\alpha-(N-2)/2}|
\partial_tG(t,x,\cdot)|_{C^{\alpha-(N-2)/2}(\overline \tau_l^h)}\d x\d t \\
&\le C\sum_{j}\iint_{Q_j(x_0)}h^{1+\alpha-N/2}|
\partial_tG(t,\cdot,x)|_{C^{\alpha-(N-2)/2}(\overline{Q_j'(x)})}\d x\d t \\
&\leq C\sum_{j}d_j^{N+2}h^{1+\alpha-N/2}
\sup_{x\in\Omega}\|\partial_tG(\cdot , \cdot ,x)\|_{L^\infty C^{\alpha-(N-2)/2}(Q_j'(x))}  \\ 
&\leq C\sum_{j}d_j^{N+2}h^{1+\alpha-N/2}
\sup_{x\in\Omega}\|\partial_tG(\cdot , \cdot, x)\|_{L^\infty H^{1+\alpha}(Q_j'(x))}  
\qquad\mbox{(because $H^{1+\alpha}(\Omega)\hookrightarrow C^{\alpha-(N-2)/2}(\overline\Omega)$)}\\ 
&\leq C\sum_{j}d_j^{N+2}h^{1+\alpha-N/2}d_j^{-\alpha-3-N/2} \qquad\mbox{(here we use Lemma \ref{GFEst1})} \\
&\leq C\sum_{j}(h/d_j)^{1+\alpha-N/2}  \\
&\leq C ,
\end{align*}
where $\sum_j$ indicates summation over $j=0,1,\dots,J_*$ (see the notations at the end of Section \ref{sec:Green}), and the last inequality is due to the fact that $h2^{J_*}\le C$. 

Substituting the estimates of ${\cal I}_1$ and ${\cal I}_2$ into \eqref{alkdkdkd} yields 
\begin{align}
\iint_{[(0,\infty)\times\Omega ]\backslash
Q_{2*}(x_0)}|\partial_t\Gamma(t,x,x_0)
-\partial_tG^*_{\rm tr}(\tau,x,x_0)|\d x\d t 
\le C.
\end{align}
Furthermore, by using the basic energy estimate, we have 
\begin{align}
\iint_{Q_{2*}(x_0)}|\partial_t\Gamma(t,x,x_0)|\d x\d t
&\leq
Cd_{*}^{N/2+1}\|\partial_t\Gamma(\cdot,\cdot ,  x_0)\|_{L^2({\cal Q})}  \nn \\
&\leq
Cd_{*}^{N/2+1}\|\widetilde\delta_{x_0} \|_{H^1(\Omega)} \nn\\
&\leq
Cd_{*}^{N/2+1} h^{-N/2-1} 
\leq CC_*^{N/2+1}  ,
\end{align}
and (cf. Property (4) of the function $G_{\rm tr}^*(t,x,x_0)$) 
\begin{align*}
&\iint_{Q_{2*}(x_0)} |\partial_tG_{\rm tr}^*(t,x,x_0)|\d
x\d t\leq Cd_*^{-N-2} d_*^{N+2} \leq C .
\end{align*}
The last three inequalities imply 
$\iint_{(0,\infty)\times\Omega}|\partial_t\Gamma(t,x,x_0)
-\partial_tG_{\rm tr}(t,x,x_0)|\d x\d t \le  C $, which together with Lemma \ref{LemGm2} further implies
\begin{align}
&\iint_{(0,\infty)\times\Omega}|\partial_t\Gamma_{h}(t,x,x_0)
-\partial_tG_{\rm tr}(t,x,x_0)|\d x\d t \le C.
\end{align}

Since both $\Gamma_h(t,x,y)$ and $G^*_{\rm tr}(t,x,y)$
are symmetric with respect to $x$ and $y$, from the last
inequality we see that the kernel
$$K(x,y):=\int_0^\infty|\partial_t\Gamma_h(t,x, y)
-\partial_tG^*_{\rm tr}(t,x,y)|\d t$$
is symmetric with respect to $x$ and $y$, and satisfies
\begin{align*}
&\sup_{y\in\Omega}\int_\Omega K(x,y)\d
x+\sup_{x\in\Omega}\int_\Omega K(x,y)\d y\leq C .
\end{align*}
By Schur's lemma \cite[Lemma 1.4.5]{Krantz1999}, the operator $M_K$ defined by
\begin{align}\label{M_K}
M_Kv(x):=\int_\Omega K(x,y)v(y)\d y
\end{align} 
is bounded on $L^q$ for all $1\leq q\leq\infty$, i.e.
\begin{align}\label{Ms1}
\|M_Kv\|_{L^q}\leq C\|v\|_{L^q} ,
\quad \forall\,\, 1\leq q\leq \infty
.
\end{align}
Then we have
\begin{align}\label{maximal_est}
&\sup_{t>0}(|E_h(t)|\, |v|)(x_0) \nn \\
&
=\sup_{t>0}\bigg|\int_\Omega|\Gamma_h(t,x,x_0)|\, |v(x)|\d x\bigg|   \nn \\
&\leq \sup_{t>0} \int_\Omega |\Gamma_h(t,x,x_0)-G^*_{\rm tr}(t,x,x_0)|\, |v(x)|\d x 
+\sup_{t>0} \int_\Omega |G^*_{\rm tr}(t,x,x_0)|\,|v(x)|\d x    \nn \\
&\le \sup_{t>0}\int_\Omega \bigg(|\Gamma_h(0,x,x_0)-G^*_{\rm tr}(0,x,x_0)|
+\int_0^t |\partial_t(\Gamma_h(s ,x,x_0)-G^*_{\rm tr}(s ,x,x_0))|\d s  \bigg) |v(x)| \d x   \nn \\
&\quad +\sup_{t>0} \int_\Omega |G^*_{\rm tr}(t,x,x_0)|\,|v(x)|\d x    \nn \\
&= \sup_{t>0}\int_\Omega \bigg(|\delta_{h,x_0}(x)|
+\int_0^t |\partial_t(\Gamma_h(s ,x,x_0)-G^*_{\rm tr}(s ,x,x_0))|\d s  \bigg) |v(x)| \d x 
+\sup_{t>0} \int_\Omega G(t,x,x_0)|u_h(x)|\d x    \nn \\
&\le \int_\Omega Kh^{-N} e^{-\frac{|x-x_0| }{K h }} |v(x)|\d x+(M_K|v|)(x) 
+\sup_{t>0} (E(t)|v|)(x_0) ,
\end{align}
where we have used \eqref{Detal-pointwise} and \eqref{M_K} in the last step.

From \eqref{GausEst1} we know that 
$G(t,x,x_0)\leq t^{-n/2}\Phi((x-x_0)/\sqrt{t})$
with $\Phi(x):=Ce^{-|x|^2/C}$, which is 
a radially decreasing and integrable function. 
Let $\widetilde u_h$ denote the zero extension of  $u_h$
from $\Omega$ to $\R^N$. 
Then Corollary 2.1.12 of \cite{Grafakos2008} implies 
\begin{align}\label{Etabsv}
\sup_{t>0} (E(t)|v|)(x_0)
&=\sup_{t>0} \int_\Omega G(t,x,x_0)|v(x)|\d x \nn\\
&\leq \sup_{t>0} \int_{\R^N}t^{-n/2}\Phi((x-x_0)/\sqrt{t})|\widetilde v(x)|\d x \nn\\
&\leq \sup_{t>0} \int_{\R^N} t^{-n/2}\Phi(x/\sqrt{t})|\widetilde v(x+x_0)|\d x \nn\\
& \leq \|\Phi\|_{L^1(\R^N)}{\cal M} |\widetilde v| , 
\end{align}
where ${\cal M} $ is the Hardy--Littlewood maximal operator. 
Since the Hardy--Littlewood maximal operator
is strong-type $(\infty,\infty)$ and weak-type $(1,1)$ 
(\cite{Grafakos2008}, Theorem 2.1.6), it follows that 
(via real interpolation)
\begin{align}\label{MLx0}
\|{\cal M}|\widetilde v| \|_{L^q(\R^N)}
\leq \frac{Cq}{q-1}\| \widetilde v \|_{L^q(\R^N)} 
= \frac{Cq}{q-1}\|v \|_{L^q(\Omega)} ,\quad 
\forall\,\, 1<q\leq \infty . 
\end{align}
By substituting \eqref{Ms1} and \eqref{Etabsv}-\eqref{MLx0} into \eqref{maximal_est}, we obtain 
\eqref{MEgodic2}. 
\qed

\subsection{Proof of (\ref{MaxLp1}) 
for $2\leq p=q< \infty$} \label{secp=q}

In this subsection, we prove (\ref{MaxLp1}) in the simple case $2\leq p=q< \infty$. 
The general case $1<p,q<\infty$ will be proved in the next subsection based on the result of this subsection, by using the mathematical tool of singular integral operators.

Let $f_h=P_hf$ and consider the expression
\begin{align}\label{dtuh}
&\partial_tu_h(t,x_0)  \nn\\
&=\partial_t\int_0^t(E_h(t-s)f_h(s,\cdot))(x_0) \d s \nn\\
&=\int_0^t(\partial_tE_h(t-s)f_h(s,\cdot))(x_0) \d s +f_h(t,x_0)   \nn\\
&
= \int_0^t 
\int_\Omega\partial_tF (t-s,x,x_0)f_h(s,x)\d x \d s 
+\int_0^t 
\int_\Omega \partial_t\Gamma(t-s,x,x_0) f_h(s,x) \d x  \d s +f_h(t,x_0)   \nn\\
&=:{\cal M}_hf_h+{\cal K}_hf_h +f_h , 
\end{align}
where ${\cal M}_h$ and ${\cal K}_h$ 
are certain linear operators. 
By Lemma \ref{LemGm2} we have
\begin{align}
\int_0^t\int_\Omega|\partial_tF (t-s,x,x_0)| \d x\d s 
\le \int_0^\infty\int_\Omega|\partial_tF (t,x,x_0)| \d x\d t  \le C,
\end{align}
which implies 
\begin{align}
&\|{\cal M}_hf_h\|_{L^\infty(\R_+;L^\infty)}
\leq C\|f_h\|_{L^\infty(\R_+;L^\infty)} .
\end{align}
Since the classical energy estimate implies 
\begin{align}
&\|{\cal M}_hf_h\|_{L^2(\R_+;L^2)}
\leq C\|f_h\|_{L^2(\R_+;L^2)} ,
\end{align}
the interpolation of the last two inequalities yields 
\begin{align}\label{MhfhLq}
&\|{\cal M}_hf_h\|_{L^q(\R_+;L^q)}
\leq C\|f_h\|_{L^q(\R_+;L^q)} ,\quad
\forall \, 2\leq q\leq\infty .
\end{align}

It remains to prove 
\begin{align}\label{KhLqq}
&\|{\cal K}_hf_h\|_{L^q(\R_+;L^q)}
\leq C_q\|f_h\|_{L^q(\R_+;L^q)} ,\quad
\forall \, 2\leq q<\infty .
\end{align}
To this end, we express ${\cal K}_hf_h$ as 
\begin{align*} 
{\cal K}_hf_h(t,x_0)
&=\int_0^t \int_\Omega \partial_t\Gamma(t-s,x,x_0) f_h(s,x) \d x  \d s \\
&=\int_0^t \int_\Omega \int_\Omega\partial_tG(t-s,x,y)\widetilde\delta_{x_0}(y) f_h(s,x) \d y\d x  \d s \\
&=\int_\Omega\widetilde\delta_{x_0}(y)
\bigg( \int_0^t  \int_\Omega\partial_tG(t-s,x,y)f_h(s,x) \d x  \d s\bigg) \d y    .
\end{align*} 
In view of \eqref{reg-Delta-est-2},  
Schur's lemma \cite[Lemma 1.4.5]{Krantz1999} implies 
\begin{align*}
&\|{\cal K}_hf_h(\cdot,t)\|_{L^q}
\leq C\bigg\|\int_0^t  \int_\Omega
\partial_tG(t-s,x,\cdot)f_h(x,s) \d x  \d s\bigg\|_{L^q} ,\quad
\forall \, 1\leq q\leq\infty ,
\end{align*}
and so
\begin{align}
\|{\cal K}_hf_h\|_{L^q(\R_+;L^q)}
&\leq C\bigg\|\int_0^t  \int_\Omega
\partial_tG(t-s,x,y)f_h(x,s) \d x  \d s\bigg\|_{L^q_t(\R_+;L^q_y)} \nn\\
&=C\bigg\|\partial_t\int_0^t  \int_\Omega
G(t-s,x,y)f_h(x,s) \d x  \d s-f_h(y,t)\bigg\|_{L^q_t(\R_+;L^q_y)} \nn\\
&=C\|\partial_tW-f_h \|_{L^q (\R_+;L^q )} ,
\end{align}
where $W$ is the solution of the PDE problem
\begin{align} 
\left\{\begin{array}{ll}
\partial_tW-\Delta W=f_h
&\mbox{in}\,\,\,\R_+ \times\Omega,\\
W = 0
&\mbox{on}\,\,\, \R_+ \times\partial\Omega ,\\
W(0,\cdot)= 0
&\mbox{in}\,\,\, \Omega , 
\end{array}\right.
\end{align}
which possesses the following maximal $L^q$-regularity (in view of \eqref{MaxLpReg}): 
\begin{align}
\|\partial_tW \|_{L^q (\R_+;L^q )}
&\leq C_q\|f_h \|_{L^q (\R_+;L^q )} ,\quad 
\forall \, 2\leq q<\infty .
\end{align} 
The last inequality implies \eqref{KhLqq}. Then substituting 
\eqref{MhfhLq}-\eqref{KhLqq} into \eqref{dtuh} yields 
\begin{align}
&\|\partial_tu_h\|_{L^q(\R_+;L^q)}
\leq C_q\|f_h\|_{L^q(\R_+;L^q)} ,
&& \forall\,\, 2\leq  q<\infty, 
\label{MaxLp2-2}
\end{align} 
Since replacing $f_h(t,x)$ by $f_h(t,x)1_{0<t<T}$ does not affect the value of $u_h(t,x)$ for $t\in (0,T)$, the last inequality implies \eqref{MaxLp1} for $2\leq p=q<\infty$. 


\subsection{Proof of (\ref{MaxLp1}) for $1< p,q<\infty$}\label{secpneqq}

In the last subsection, 
we have proved \eqref{MaxLp1} for $2\leq p=q<\infty$
by showing that the operator 
${\cal E}_h$ defined by 
\begin{align}\label{defEhOpt}
({\cal E}_hf_h)(t,\cdot):
=\int_0^t\partial_tE_h(t-s)f_h(s,\cdot) \d s
=({\cal M}_hf_h)(t,\cdot)+({\cal K}_hf_h)(t,\cdot)
\end{align}
satisfies
\begin{align}\label{qlarger2}
&\|{\cal E}_hf_h\|_{L^q(\R_+;L^q)}\leq C_q\|f_h\|_{L^q(\R_+;L^q)},
\quad \forall\,\, 2\leq q<\infty .
\end{align}
In this subsection, we prove \eqref{MaxLp1} for all $1< p,q<\infty$ via a duality argument and the singular integral operator approach. 

In fact, by the same method, one can also prove that
the operator ${\cal E}_h'$ defined by 
\begin{align}
({\cal E}_h'f_h)(s,\cdot)=\int_s^\infty \partial_tE_h(t-s)f_h(t,\cdot) \d t
\end{align}
satisfies
\begin{align}
&\|{\cal E}_h'f_h\|_{L^q(\R_+;L^q)}\leq C_q\|f_h\|_{L^q(\R_+;L^q)},
\quad \forall\,\, 2\leq q<\infty .
\end{align}
Since ${\cal E}_h'$ is the dual of ${\cal E}_h$,
by duality we have
\begin{align}\label{qsmaller2}
&\|{\cal E}_hf_h\|_{L^q(\R_+;L^q)}
\leq C_{q}\|f_h\|_{L^q(\R_+;L^q)} ,
\quad \forall\,\, 1< q\leq 2 .
\end{align}
The two inequalities \eqref{qlarger2} and \eqref{qsmaller2} 
can be summarized as
\begin{align}\label{qgeneral}
&\|{\cal E}_hf_h\|_{L^q(\R_+;L^q)} 
\leq C_q\|f_h\|_{L^q(\R_+;L^q)},
\quad \forall\,\, 1< q<\infty .
\end{align}
Therefore, we have
\begin{align}
&\|\partial_tu_h\|_{L^q(\R_+;L^q)}
\leq C_q \|f_h\|_{L^q(\R_+;L^q)} ,
&& \forall\,\, 1< q<\infty .
\label{MaxLp2-3}
\end{align}

Overall, for any fixed $1<q<\infty$ the operator ${\cal E}_h$ is bounded on $L^q(\R_+;L^q)$, and $\{E_h(t)\}_{t>0}$ is an analytic semigroup satisfying (see Lemma \ref{LemGm2}):  
\begin{align}
&\|\partial_tE_h(t-s) \|_{{\cal L}(L^q,L^q)}\leq C(t-s)^{-1},
&&\forall\,\, t>s>0, \label{DtEhLq}\\
&\|\partial_{tt}E_h(t-s) \|_{{\cal L}(L^q,L^q)}\leq C(t-s)^{-2},
&&\forall\,\, t>s>0 .  \label{DttEhLq}
\end{align}
From \eqref{defEhOpt} and \eqref{DtEhLq}-\eqref{DttEhLq} 
we see that ${\cal E}_h$ is an operator-valued singular integral operator whose kernel $\partial_tE_h(t-s)1_{t>s}$ satisfying 
the H\"ormander conditions (cf. \cite[condition (4.6.2)]{Grafakos2008}):
\begin{align} 
&\sup_{s,s_0\in\R}\int_{|t-s_0|\geq 2|s-s_0|}\|\partial_tE_h(t-s)1_{t>s}-\partial_tE_h(t-s_0)1_{t>s_0}\|_{{\cal L}(L^q,L^q)}\,\d t
\leq C, \label{Hormandercond}\\
&\sup_{t,t_0\in\R}\int_{|t_0-s|\geq 2|t-t_0|}\|\partial_tE_h(t-s)1_{t>s}-\partial_tE_h(t_0-s)1_{t_0>s}\|_{{\cal L}(L^q,L^q)}\,\d s\leq
C. \label{Hormandercond2} \medskip
\end{align}
Under the conditions \eqref{Hormandercond}-\eqref{Hormandercond2}, the theory of singular integral operators (cf. \cite[Theorem 4.6.1]{Grafakos2008}) says that if ${\cal E}_h$ is bounded on $L^q(\R_+;  L^q)$ for some $q\in(1,\infty)$ as proved in \eqref{MaxLp2-3}, then it is bounded on $L^p(\R_+;  L^q)$ for all $p\in(1,\infty)$: \begin{align}
&\|{\cal E}_hf_h\|_{L^p(\R_+;L^q)}
\leq \max(p\, ,\, (p-1)^{-1}) \, 
C_q\|f_h\|_{L^p(\R_+;L^q)} .
\end{align} 
Since replacing $f_h(t,x)$ by $f_h(t,x)1_{0<t<T}$ does not affect the value of $u_h(t,x)$ for $t\in (0,T)$, the last inequality implies \eqref{MaxLp1} for all $1<p,q<\infty$.

\subsection{Proof of (\ref{MaxLinfty})}

Again, we consider 
\begin{align}\label{defEhOpt0}
({\cal E}_hf_h)(t,\cdot) 
=\int_0^t\partial_tE_h(t-s)f_h(s,\cdot) \d s
=\int_0^t  \int_\Omega
\partial_t\Gamma_h(t-s,x, \cdot)f_h(s,x) \d x  \d s 
\end{align}
and use the following inequality: for $t\in(0,T)$ 
\begin{align}
\|{\cal E}_hf_h(t,\cdot)\|_{L^\infty}
&\leq \bigg(\sup_{x_0\in\Omega}\int_0^t  \int_\Omega
|\partial_t\Gamma_h(t-s,x,x_0)|\d x  \d s \bigg) 
\|f_h\|_{L^\infty(0,T;L^\infty)} \nn\\
&\le \bigg(\sup_{x_0\in\Omega}\int_0^\infty  \int_\Omega
|\partial_t\Gamma_h(t,x,x_0)|\d x   \d t \bigg) 
\|f_h\|_{L^\infty(0,T;L^\infty)} \nn \\
&\le 
\bigg(\int_0^\infty\sup_{x_0\in\Omega}  \int_\Omega
|\partial_t\Gamma_h(t,x,x_0)|\d x   \d t \bigg) 
\|f_h\|_{L^\infty(0,T;L^\infty)}  ,      \\[5pt]
\|{\cal E}_hf_h(t,\cdot)\|_{L^1} \,\,
&\leq \bigg(\int_0^t\sup_{x\in\Omega}  \int_\Omega
|\partial_t\Gamma_h(t-s,x,x_0)|\d x_0  \d s \bigg) 
\|f_h\|_{L^\infty(0,T;L^1)} \nn \\
&= \bigg(\int_0^t \sup_{x_0\in\Omega}  \int_\Omega
|\partial_t\Gamma_h(t-s,x,x_0)|\d x   \d s \bigg) 
\|f_h\|_{L^\infty(0,T;L^1)} \nn\\
&\le \bigg(\int_0^\infty\sup_{x_0\in\Omega}  \int_\Omega
|\partial_t\Gamma_h(t,x,x_0)|\d x   \d t \bigg) 
\|f_h\|_{L^\infty(0,T;L^1)} ,
\end{align}
where we have used the symmetry 
$\partial_t\Gamma_h(t-s,x,x_0)=\partial_t\Gamma_h(t-s,x_0,x)$, 
due to the self-adjointness of the operator 
$E_h(t-s)$. By interpolation between $L^\infty$ and $L^1$,
we get
\begin{align}
\|{\cal E}_hf_h\|_{L^\infty(0,T;L^{q})}
&\leq \bigg(\int_0^\infty\sup_{x_0\in\Omega}  \int_\Omega
|\partial_t\Gamma_h(t,x,x_0)|\d x   \d t \bigg)
\|f_h\|_{L^\infty(0,T;L^{q})} ,\quad\forall\, 1\le q\le \infty .
\end{align}
It remains to prove
\begin{align} \label{GtL1ST}
\int_0^\infty\sup_{x_0\in\Omega} \int_\Omega
|\partial_t\Gamma_h(t,x,x_0)|\d x   \d t 
\leq C\log(2+1/h) .
\end{align}
To this end, we note that $\partial_t\Gamma_h(t,\cdot,x_0)=\Delta_h\Gamma_h(t,\cdot,x_0)=E_h(t)\Delta_hP_h\widetilde\delta_{x_0}$. By using \eqref{L1Gammh} of Lemma \ref{LemGm2} and \eqref{analyticity} (proved in Section \ref{sec:proof2.3}), we have 
\begin{align}
&\|\partial_t\Gamma_h(t,\cdot,x_0)\|_{L^1}
\le Ct^{-1} ,\\
&\|\partial_t\Gamma_h(t,\cdot,x_0)\|_{L^1} 
\le C\|\Delta_hP_h\widetilde\delta_{x_0} \|_{L^1}
\le Ch^{-2} \|P_h\widetilde\delta_{x_0} \|_{L^1}
\le Ch^{-2} . 
\end{align} 
The interpolation of the last two inequalities gives
$\displaystyle\|\partial_t\Gamma_h(t,\cdot,x_0)\|_{L^1}
\leq \frac{C}{h^{2\theta}t^{1-\theta}}$ for arbitrary $\theta\in(0,1)$, 
where the constant $C$ is independent of $\theta$.
Hence, we have 
$\displaystyle\int_0^1\sup_{x_0\in\Omega}\|\partial_t\Gamma(t,\cdot,x_0)\|_{L^1}\d t
\leq \frac{C}{\theta h^{2\theta}}$ for arbitrary $\theta\in(0,1) .$
By choosing $\theta=1/\log(2+1/h)$, we obtain 
\begin{align}
&\int_0^1\sup_{x_0\in\Omega} \int_\Omega
|\partial_t\Gamma(t,x, x_0)|\d x  \d t\leq  C\log(2+1/h) .
\end{align}
The estimate \eqref{L1Gammatx0} implies 
\begin{align}
&\int_1^\infty \sup_{x_0\in\Omega}\int_\Omega
|\partial_t\Gamma(t,x, x_0)|\d x  \d t\leq  C .
\end{align}
The last two inequalities imply 
\eqref{GtL1ST}, and 
this completes the proof of \eqref{MaxLinfty}.

%
%
%
%

The proof of Theorem \ref{MainTHM1} is complete (up to the proof of Lemma \ref{LemGm2}).
\qed

\section{Proof of Lemma \ref{LemGm2}}
\label{ProofLm03}
\setcounter{equation}{0}

In this section we prove Lemma \ref{LemGm2}, which is used in proving Theorem \ref{MainTHM1} in the last section. 
To this end, we use the following local energy error estimate for finite element solutions of parabolic equations, which extends the existing work \cite[Lemma 6.1]{SchatzThomeeWahlbin1998} and \cite[Proposition 3.2]{LiSun2016} to nonconvex polyhedra without using the superapproximation results of the Ritz projection (cf. \cite[Theorem 5.1]{SchatzThomeeWahlbin1998} and \cite[Proposition 3.1]{LiSun2016}, which only hold in convex domains).

\begin{lemma}\label{LocEEst} 
{\it 
Suppose that $\phi\in L^2(0,T;H^1_0(\Omega))\cap H^1(0,T;L^2(\Omega))$ and 
$\phi_h\in H^1(0,T; S_h)$ satisfy the equation 
\begin{align}  
(\partial_t(\phi-\phi_h),\chi_h)
+(\nabla (\phi-\phi_h),\nabla \chi_h) =0, 
\quad\forall\, \chi_h\in  S_h, \,\, \mbox{a.e.}\,\, t>0, 
\label{p32}  
\end{align} 
with $\phi(0)=0$ in $\Omega_j''$. 
Then, under the assumptions of Theorem \ref{MainTHM1},  we have 
\begin{align} 
&\vertiii{\partial_t(\phi-\phi_h)}_{Q_j} 
+ d_j^{-1}\vertiii{\phi-\phi_h}_{1,Q_j} \nn\\
& \leq 
C\epsilon^{-3}\big(I_j(\phi_{h}(0))+X_j(I_h\phi-\phi) 
+d_j^{-2}\vertiii{\phi-\phi_h}_{Q_j'}\big) \nn\\
&\quad 
+(Ch^{1/2}d_j^{-1/2}+C\epsilon^{-1}hd_j^{-1}+\epsilon) \big(\vertiii{\partial_t(\phi-\phi_h)}_{Q_j'}
+d_j^{-1}\vertiii{\phi-\phi_h}_{1,Q_j'}\big), 
\label{LocEngErr}  
\end{align} 
where
\begin{align*}
&I_j(\phi_{h}(0))=\|\phi_{h}(0)\|_{1,\Omega_j'} +d_j^{-1}\|\phi_{h}(0)\|_{\Omega_j'} \, ,\\[5pt]
&X_j(I_h\phi-\phi)=d_j\vertiii{\partial_t(I_h\phi-\phi)}_{1,Q_j'}
+\vertiii{\partial_t(I_h\phi-\phi)}_{Q_j'} \nn\\
&\qquad\qquad\qquad\,\, 
+d_j^{-1}\vertiii{I_h\phi-\phi}_{1,Q_j'}+
d_j^{-2}\vertiii{I_h\phi-\phi}_{Q_j'}  \, ,
\end{align*}
where $\epsilon\in(0,1)$ is an arbitrary positive constant, and 
the positive constant $C$ is independent of $h$, $j$ and $C_*$; the norms $\vertiii{\cdot}_{k,Q_j'}$ and $\vertiii{\cdot}_{k,\Omega_j'}$ are defined in \eqref{Q-norm}. 
}
\end{lemma}

The proof of Lemma \ref{LocEEst} is presented in Appendix A. In the rest of this section, we apply Lemma \ref{LocEEst} to prove Lemma \ref{LemGm2} by denoting $\alpha\in(\frac{1}{2},1]$ a fixed constant satisfying Lemma \ref{GFEst1}. The proof consists of three parts. The first part is concerned with estimates for $t \in (0, 1)$, where we covert the $L^1$ estimates on ${\cal Q}=(0,1)\times\Omega=Q_*\bigcup \big(\cup_{j=0}^JQ_j\big)$ into weighted $L^2$ estimates on the subdomains $Q_*$ and $Q_j$, $j=0,1,\dots,J$. 
The second part is concerned with estimates for $t\ge 1$, which is a simple consequence of the parabolic regularity. 
The third part is concerned with the proof of \eqref{L1Gammh}-\eqref{L1Gammh-2}, which are simple consequences of the results proved in the first two parts.

\medskip

{\it Part I.}$\,\,\,$ 
First, we present estimates in the domain ${\cal Q}=(0,1)\times\Omega$ with the restriction $h<1/(4C_*) $; see \eqref{h-condition}. In this case, the basic energy estimate gives 
\begin{align} 
&\|\partial_t\Gamma\|_{L^2({\cal Q})}
+\|\partial_t\Gamma_h\|_{L^2({\cal Q})}
\leq C(\|\Gamma(0)\|_{H^1} +\|\Gamma_h(0)\|_{H^1})
\leq Ch^{-1-N/2} ,  \label{EstGamma-Q1} \\
&\|\Gamma\|_{L^\infty L^2({\cal Q})} 
+\|\Gamma_h\|_{L^\infty L^2({\cal Q})} 
\leq C(\|\Gamma(0)\|_{L^2} +\|\Gamma_h(0)\|_{L^2})
\leq Ch^{-N/2} ,\\
&\|\nabla\Gamma\|_{L^2({\cal Q})} +\|\nabla\Gamma_h\|_{L^2({\cal Q})}
\leq C(\|\Gamma(0)\|_{L^2} +\|\Gamma_h(0)\|_{L^2})
\leq Ch^{-N/2} ,\\
&\|\partial_{tt}\Gamma\|_{L^2({\cal Q})}
+\|\partial_{tt}\Gamma_h\|_{L^2({\cal Q})}
\leq C(\|\Delta\Gamma(0)\|_{H^1} +\|\Delta_h\Gamma_h(0)\|_{H^1})
\leq Ch^{-3-N/2} ,  \label{EstGamma-Q4}  \\
&\|\nabla\partial_t\Gamma\|_{L^2({\cal Q})}
+\|\nabla\partial_t\Gamma_h\|_{L^2({\cal Q})}
\leq C(\|\Delta\Gamma(0)\|_{L^2} +\|\Delta_h\Gamma_h(0)\|_{L^2})
\leq Ch^{-2-N/2} , \label{EstGamma-Q5}
\end{align} 
where we have used \eqref{reg-Delta-est} and \eqref{Detal-pointwise} to estimate $\Gamma(0)$ and $\Gamma_h(0)$, respectively. Hence, we have 
\begin{align}\label{L2Gamma-Q}
\vertiii{\Gamma}_{Q_*} +\vertiii{\Gamma_h}_{Q_*} 
\le Cd_*\|\Gamma\|_{L^\infty L^2(Q_*)}+Cd_*\|\Gamma_h\|_{L^\infty L^2(Q_*)}
\le Cd_* h^{-N/2} \le CC_*h^{1-N/2} .
\end{align} 
Since the volume of $Q_j$ is $Cd_j^{2+N}$, we can decompose $\|\partial_tF\|_{L^1({\cal Q})}+\|t\partial_{tt}F\|_{L^1({\cal Q})}$ in the following way:
\begin{align}\label{Bd31K2}
&\|\partial_tF\|_{L^1({\cal Q})}+\|t\partial_{tt}F\|_{L^1({\cal Q})} \nn\\
&\leq  \|\partial_{t} F\|_{L^1(Q_*)}
+ \|t\partial_{tt} F\|_{L^1(Q_*)}
+\sum_{j}\big(\|\partial_{t} F\|_{L^1(Q_j)}
+ \|t\partial_{tt} F\|_{L^1(Q_j)}\big) \nn\\
&\leq Cd_{*}^{1+N/2}\big(\vertiii{\partial_tF}_{Q_*} 
+d_*^2\vertiii{\partial_{tt} F}_{Q_*}\big)  
+\sum_{j}Cd_j^{1+N/2} 
\big(\vertiii{\partial_tF}_{Q_j} 
+d_j^2\vertiii{\partial_{tt} F}_{Q_j}\big) \nn\\
&\leq CC_*^{3+N/2}
+ {\mathscr K} ,
\end{align} 
where we have used \eqref{EstGamma-Q1} and \eqref{EstGamma-Q4} to estimate 
$Cd_{*}^{1+N/2}\big(\vertiii{\partial_tF}_{Q_*} 
+d_*^2\vertiii{\partial_{tt} F}_{Q_*}\big)$, and introduced the notation 
\begin{align}\label{express-K}
{\mathscr K}  :&=
\sum_{j}d_j^{1+N/2}(d_j^{-1}\vertiii{F}_{1,Q_j}
+\vertiii{\partial_tF}_{Q_j}
+d_{j}\vertiii{\partial_tF}_{1,Q_j}+d_{j}^2\vertiii{\partial_{tt}F}_{Q_j})  .
\end{align} 

It remains to estimate ${\mathscr K} $. To this end, 
we set ``$\phi_h=\Gamma_h$, $\phi=\Gamma$, $\phi_h(0)=P_h\widetilde\delta_{x_0}$ and $\phi(0)=\widetilde\delta_{x_0}$'' 
and ``$\phi_h=\partial_t\Gamma_h$, $\phi=\partial_t\Gamma$, $\phi_h(0)=\Delta_hP_h\widetilde\delta_{x_0}$ and $\phi(0)=\Delta\widetilde\delta_{x_0}$''
in Lemma \ref{LocEEst}, respectively.
Then we obtain 
\begin{align}\label{F1tE}
d_{j}^{-1}\vertiii{F}_{1,Q_j}+\vertiii{\partial_tF}_{Q_j} 
&\le 
C\epsilon_1^{-3}(\widehat{I_j}+\widehat{X_j}+ d^{-2}_j\vertiii{ F }_{Q'_j} )  \\
&\quad +(Ch^{1/2}d_j^{-1/2}+C\epsilon_1^{-1}hd_j^{-1}+\epsilon_1) \big( d_j^{-1}\vertiii{F}_{1,Q_j'}
+ \vertiii{\partial_{t}F }_{Q_j'}\big) \qquad \quad\,\, \nn
\end{align}
and
\begin{align}\label{F1ttE}
d_{j}\vertiii{\partial_tF}_{1,Q_j}+d_{j}^2\vertiii{\partial_{tt}F}_{Q_j} 
&\le 
C\epsilon_2^{-3} (\overline {I_j}+\overline{X_j} + \vertiii{ \partial_tF }_{Q'_j} )  \\
&\quad 
+(Ch^{1/2}d_j^{-1/2}+C\epsilon_2^{-1}hd_j^{-1}+\epsilon_2) \big(  d_j\vertiii{\partial_{t}F}_{1,Q_j'}
+d_{j}^2\vertiii{\partial_{tt}F }_{Q_j'} \big) 
\, , \nn
\end{align}
respectively, 
where $\epsilon_1,\epsilon_2\in(0,1)$
are arbitrary positive constants. 
By using \eqref{H1+slocal} (local interpolation error estimate), \eqref{Detal-pointwise} (exponential decay of $ P_h\widetilde \delta_{x_0}$) and Lemma \ref{GFEst1} (estimates of regularized Green's function), we have 
\begin{align}
&\widehat{I_j}=\|P_h\widetilde\delta_{x_0}\|_{1,\Omega_j'}
+d_j^{-1}\|P_h\widetilde\delta_{x_0}\|_{\Omega_j'} 
\leq Ch^2d_{j}^{-3-N/2},  \label{EstIj}  \\[5pt]
&\widehat{X_j}=
d_j\vertiii{\partial_t(I_h\Gamma-\Gamma)}_{1,Q_j'}
+\vertiii{\partial_t(I_h\Gamma-\Gamma)}_{Q_j'} \nn\\
&\qquad\, 
+d_j^{-1}\vertiii{I_h\Gamma-\Gamma}_{1,Q_j'}
+ d_j^{-2}\vertiii{I_h\Gamma-\Gamma}_{Q_j'}  \nn\\
&\quad~
\leq ( d_j h^{\alpha}+ h^{1+\alpha})\vertiii{\partial_{t}\Gamma}_{L^2H^{1+\alpha}(Q_j'')}
+(d_j^{-1}h^{\alpha}+ d_j^{-2}h^{1+\alpha})\vertiii{\Gamma}_{L^2H^{1+\alpha}(Q_j'')} \nn\\
&\quad~
\leq C h^{\alpha} d_j^{-1-\alpha-N/2}   .
\end{align}
and
\begin{align}
&\overline{I_j}=d_{j}^{2}\|\Delta_hP_h\widetilde\delta_{x_0}\|_{1,\Omega_j'}
+d_j \|\Delta_hP_h\widetilde\delta_{x_0}\|_{\Omega_j'}
\leq Ch^2d_{j}^{-3-N/2},\\[5pt]
&\overline{X_j}=
d_j^3\vertiii{I_h\partial_{tt}\Gamma-\partial_{tt}\Gamma }_{1,Q_j'}
+d_j^2\vertiii{ I_h\partial_{tt}\Gamma-\partial_{tt}\Gamma}_{Q_j'} \nn\\
&\qquad\, 
+d_j \vertiii{I_h\partial_{t}\Gamma-\partial_{t}\Gamma}_{1,Q_j'}+
 \vertiii{I_h\partial_{t}\Gamma-\partial_{t}\Gamma}_{Q_j'}  \nn\\
&\quad~
\leq ( d_j^3 h^{\alpha}+ d_j^{ 2}h^{1+\alpha})\vertiii{\partial_{tt}\Gamma}_{L^2H^{1+\alpha}(Q_j'')}
+(d_j h^{\alpha}+ h^{1+\alpha})\vertiii{\partial_{t}\Gamma}_{L^2H^{1+\alpha}(Q_j'')}  \nn\\
&\quad~
\leq C h^{\alpha} d_j^{-1-\alpha-N/2}   .
\label{ovXj}
\end{align}
By choosing $\epsilon_1=\epsilon^4$ and $\epsilon_2=\epsilon$ in \eqref{F1tE}-\eqref{F1ttE}, and substituting \eqref{F1tE}-\eqref{ovXj} into the expression of ${\mathscr K}$ in \eqref{express-K}, we have
\begin{align} 
{\mathscr K} 
&=\sum_{j}d_j^{1+N/2}(d_j^{-1}\vertiii{F}_{1,Q_j}
+\vertiii{\partial_tF}_{Q_j}
+d_{j}\vertiii{\partial_tF}_{1,Q_j}
+d_{j}^2\vertiii{\partial_{tt}F}_{Q_j})  \nn\\
&\leq C_\epsilon  \sum_{j}d_j^{1+N/2}
\big(h^2d_j^{-3-N/2}+h^{\alpha}d_j^{-1-s-N/2}
+d_j^{-2}\vertiii{ F }_{Q'_j} \big) \nn\\
&\quad 
+(C_\epsilon h^{1/2}d_j^{-1/2}+C_\epsilon hd_j^{-1}+\epsilon)\sum_{j}d_j^{1+N/2}
 (d_j^{-1}\vertiii{F}_{1,Q_j'} +\vertiii{\partial_tF}_{Q_j'}) \nn\\
&\quad 
+(C_\epsilon h^{1/2}d_j^{-1/2}+C_\epsilon hd_j^{-1}+\epsilon)\sum_{j}d_j^{1+N/2}(d_{j}\vertiii{\partial_tF}_{1,Q_j'}
+d_{j}^2\vertiii{\partial_{tt}F}_{Q_j'})  
 \nn\\
&\leq C_\epsilon 
+C_\epsilon \sum_{j}d_j^{-1+N/2}\vertiii{F}_{Q_j'}  \nn\\
&\quad 
+(C_\epsilon h^{1/2}d_j^{-1/2}+C_\epsilon hd_j^{-1}+\epsilon)\sum_{j}d_j^{1+N/2}
 (d_j^{-1}\vertiii{F}_{1,Q_j'} +\vertiii{\partial_tF}_{Q_j'}) \nn\\
&\quad 
+(C_\epsilon h^{1/2}d_j^{-1/2}+C_\epsilon hd_j^{-1}+\epsilon)\sum_{j}d_j^{1+N/2}(d_{j}\vertiii{\partial_tF}_{1,Q_j'}
+d_{j}^2\vertiii{\partial_{tt}F}_{Q_j'})    .
\end{align}

Since $\vertiii{F}_{Q_j'} \le C(\vertiii{F}_{Q_{j-1}} + \vertiii{F}_{Q_{j}} + \vertiii{F}_{Q_{j+1}})$, 
we can convert the $Q_j'$-norm in the inequality above to the $Q_j$-norm:
\begin{align}\label{slkll}
{\mathscr K} &\leq C_\epsilon
+ C_\epsilon\sum_{j}d_j^{-1+N/2}\vertiii{F}_{Q_j} 
+C_\epsilon d_{*}^{-1+N/2}\vertiii{F}_{Q_*}\nn\\ 
&\quad 
+(C_\epsilon h^{1/2}d_j^{-1/2}+C_\epsilon hd_j^{-1}+\epsilon)\sum_{j}d_j^{1+N/2}
 (d_j^{-1}\vertiii{F}_{1,Q_j} +\vertiii{\partial_tF}_{Q_j}) \nn\\
&\quad 
+(C_\epsilon h^{1/2}d_j^{-1/2}+C_\epsilon hd_j^{-1}+\epsilon)\sum_{j}d_j^{1+N/2}(d_{j}\vertiii{\partial_tF}_{1,Q_j}
+d_{j}^2\vertiii{\partial_{tt}F}_{Q_j}) \nn\\
&\quad 
+(C_\epsilon h^{1/2}d_*^{-1/2}+C_\epsilon hd_*^{-1}+\epsilon) d_*^{1+N/2}
 (d_*^{-1}\vertiii{F}_{1,Q_*} +\vertiii{\partial_tF}_{Q_*}) \nn\\
&\quad 
+(C_\epsilon h^{1/2}d_*^{-1/2}+C_\epsilon hd_*^{-1}+\epsilon)
d_*^{1+N/2}(d_*\vertiii{\partial_tF}_{1,Q_*}
+d_*^2\vertiii{\partial_{tt}F}_{Q_*}) \nn\\ 
&\leq 
C_\epsilon  
+C_\epsilon C_*^{3+N/2}+ 
\sum_{j}Cd_j^{-1+N/2}\vertiii{F}_{Q_j} 
+C(C_\epsilon C_*^{-1/2}+C_\epsilon C_*^{-1}+\epsilon) {\mathscr K}. 
\end{align}
where we have used $d_j\ge C_*h$ and \eqref{EstGamma-Q1}-\eqref{EstGamma-Q5} to estimate 
$$
\mbox{$\vertiii{F}_{1,Q_*}$,\,\, $\vertiii{\partial_tF}_{Q_*}$,\,\, $\vertiii{\partial_tF}_{1,Q_*}$\,\, and\,\, 
$\vertiii{\partial_{tt}F}_{Q_*}$}
$$
and used the expression of ${\mathscr K}$ in \eqref{express-K} to bound the terms involving $Q_j$.  
By choosing $\epsilon$ small enough 
and then choosing $C_*$ large enough 
($C_*$ is still to be determined later), the last term on the right-hand side of \eqref{slkll} will be absorbed by the left-hand side. Hence, we obtain
\begin{align} 
{\cal K}
&\leq
C
+CC_*^{3+N/2}+
\sum_{j}Cd_j^{-1+N/2}\vertiii{F}_{Q_j} .
\end{align}
%
%

It remains to estimate $\vertiii{F}_{Q_j}$. To this end, we apply a duality argument below.
Let $w$ be the solution of the backward parabolic equation
$$
-\partial_tw-\Delta w=v\quad\mbox{with}~~w(T)=0,
$$
where $v$ is a function supported on $Q_j$ and
$\vertiii{v}_{Q_j}=1$. Multiplying the above equation by $F$ yields (using integration by parts, with the notations \eqref{inner-products})
\begin{align}\label{dka6}
[F,v]=(F(0),w(0))+[F_t,w]+[\nabla 
F,\nabla w] ,
\end{align}
where (since $\widetilde\delta_{x_0}=0$ on $\Omega_j''$)
\begin{align*}
(F(0),w(0))&=(P_h\widetilde\delta_{x_0}-\widetilde\delta_{x_0},w(0))\\
&=(P_h\widetilde\delta_{x_0}-\widetilde\delta_{x_0},w(0)-I_hw(0))\\
&=
(P_h\widetilde\delta_{x_0},w(0)-I_hw(0))_{\Omega_j''}
+(P_h\widetilde\delta_{x_0}-\widetilde\delta_{x_0},
w(0)-I_hw(0))_{(\Omega_j'')^c}\\
&=: {\cal I}_1+{\cal I}_2 .
\end{align*}
By using Property (P3) and \eqref{Detal-pointwise} (the exponential decay of $P_h\widetilde\delta_{x_0}$), we derive that
\begin{align}
&|{\cal I}_1|\leq
Ch\|P_h\widetilde\delta_{x_0}\|_{L^2(\Omega_j'')}\|w(0)\|_{H^1(\Omega)}  
\nn\\
&\quad\,\,\, \leq
Ch^{1-N/2}e^{-Cd_j/h}\vertiii{v}_{Q_j} 
\nn\\
&\quad\,\,\, =
C(d_j/h)^{1+N/2}e^{-Cd_j/h}h^{2}d_j^{-1-N/2}  
\nn\\
&\quad\,\,\, \leq
Ch^{2}d_j^{-1-N/2 }, \label{SF2}\\[5pt]
&|{\cal I}_2|\leq C\| \widetilde\delta_{x_0}\|_{L^{2} }
\|w(0)-I_hw(0)\|_{L^{2}((\Omega_j'')^c)}  
\nn\\
&\quad\,\,\, \leq
Ch^{1+\alpha}\| \widetilde\delta_{x_0}\|_{L^{2} } \inf_{\widetilde w} \|\widetilde w\|_{H^{1+\alpha}(\Omega)} \nn\\
&\quad\,\,\, =
Ch^{1+\alpha}\| \widetilde\delta_{x_0}\|_{L^{2} } \|w(0)\|_{H^{1+\alpha}((\Omega_j')^c)}  \nn\\
&\quad\,\,\, \leq
Ch^{1+\alpha-N/2}\|w(0)\|_{H^{1+\alpha}((\Omega_j')^c)} ,
\label{SF22}
\end{align}
where the infimum extends over all possible $\widetilde w$ extending $w(0)$ from $(\Omega_j')^c$ to $\Omega$, 
and we have used \eqref{reg-Delta-est} in the last step.  

To estimate $\|w(0)\|_{H^{1+\alpha}((\Omega_j')^c)}$, 
we let $W_j$ be a set containing $(\Omega_j')^c$ but its
distance to $\Omega_j$ is larger than $C^{-1}d_j$. 
Since  
$$
|x-y|+s ^{1/2}\geq C_1^{-1}d_j 
\quad\mbox{for $x\in W_j$ and $(s ,y)\in Q_j$}
$$
for some positive constant $C_1$, it follows that 
$(s ,x)\in \bigcup_{k\le j+\log_2C_1} Q_{k}(y)$ for $(s ,y)\in Q_j$. 
Now, if we denote $\widetilde G(\cdot,\cdot,y)$ as any extension of $G(\cdot,\cdot,y)$ from $\bigcup_{k\le j+\log_2C_1} Q_{k}(y)$ to ${\cal Q}$, then for $x\in W_j$ we have 
$$
w(0,x)=\int_0^{T}\int_{\Omega}
G(s,y,x)v(s ,y)\d y\d s  
=\iint_{Q_j} 
G(s ,x,y)v(s ,y)\d y\d s 
=\iint_{Q_j}
\widetilde G(s ,x,y)v(s ,y)\d y\d s  ,
$$
where we have used the symmetric $G(s ,y,x)=G(s ,x,y)$ and the compact support of $v$ in $Q_j$. 
Hence, we have 
\begin{align}
\|w(0,\cdot)\|_{H^{1+\alpha}(W_j)}
 \le \|w(0,\cdot)\|_{H^{1+\alpha}(\Omega)} 
&\leq C  \int_0^t\int_\Omega \|
 \widetilde G(s ,\cdot,y)\|_{H^{1+\alpha}(\Omega)} |v(s,y)| \d y\d s  \nn\\
&\leq C\sup_{y\in\Omega}\|\widetilde G(\cdot,\cdot,y)\|_{L^\infty H^{1+\alpha}(\Omega)}\|v\|_{L^{1}(Q_j)} .\
\end{align}
Since the last inequality holds for all possible $\widetilde G(\cdot,\cdot,y)$ extending $G(\cdot,\cdot,y)$ from $\bigcup_{k\le j+\log_2C_1} Q_{k}(y)$ to ${\cal Q}$, it follows that (cf. definition \eqref{DefLpX}) 
\begin{align}
\|w(0,\cdot)\|_{H^{1+\alpha}(W_j)}
&\leq C\sup_{y\in\Omega}\|G(\cdot,\cdot,y)\|_{L^\infty H^{1+\alpha}(\bigcup_{k\le j+\log_2C_1} Q_{k}(y))}\|v\|_{L^{1}(Q_j)} \nn\\
&\leq C d_j^{-1-\alpha-N/2}\|v\|_{L^{1}(Q_j)}
\qquad\qquad
\mbox{(here we use \eqref{GFest03})}
 \nn\\
&\leq Cd_j^{-\alpha} \vertiii{v}_{Q_j} 
=C d_j^{-\alpha}  . \label{SF3}
\end{align}
From (\ref{SF2})-(\ref{SF3}), we see that the 
first term on the right-hand side of \eqref{dka6} is bounded by
\begin{align}
|(F(0),w(0))| \leq  Ch^{2}d_j^{-N/2 -1}
+Ch^{1+\alpha-N/2}d_j^{-\alpha} 
 \leq  Ch^{1+\alpha-N/2}d_j^{-\alpha}  ,
\label{f69}
\end{align}
and the rest terms are bounded by (recall that $F=\Gamma_h-\Gamma$, where $\Gamma_h$ and $\Gamma$ are solutions of \eqref{RGFdef}-\eqref{EqGammh})
\begin{align}\label{sd80}
[F_t,w]+[\nabla  F,\nabla  w]  
&=[F_t,w-I_hw]+[\nabla F,\nabla (w-I_hw)]   
\nn\\
&\leq 
\sum_{*,i}\left(C\vertiii{F_t}_{Q_i} \vertiii{w-I_hw}_{Q_i}
+ C \vertiii{F}_{1,Q_i}\vertiii{w-I_hw}_{1,Q_i}\right) \nn\\
&\leq 
\sum_{*,i}(Ch^{1+\alpha}\vertiii{F_t}_{Q_i} 
+ Ch^{\alpha} \vertiii{F}_{1,Q_i})
\|w\|_{L^2H^{1+\alpha}(Q_i')} .
\end{align}
where we have used Property (P3) of Section \ref{Sec2-2} in the last step. 

To estimate $\|w\|_{L^2H^{1+\alpha}(Q_i')}$ in \eqref{sd80}, 
we consider the expression ($v$ is supported in $Q_j$)
\begin{align}\label{dllkka}
w(t,x)
=\int_0^{T}\int_{\Omega}
G(s -t,x,y)v(s ,y)1_{s >t}\,\d y\d s 
=\iint_{Q_j}
G(s -t,x,y)v(s ,y)1_{s >t}\,\d y\d s  .
\end{align}
For $i\leq j-3$ (so that $d_i>d_j$), we have 
\begin{align}
\begin{aligned}
&\mbox{if $t>4d_j^2$ then $w(t,x)=0$ (because $v$ is supported in $Q_j$);}\\ 
&\mbox{if $t\le 4d_j^2$, $(t,x)\in Q_i'$ and
$(s ,y)\in Q_j$, then $d_i/4\leq |x-y|\leq 4 d_i$ and $s -t\in(0,d_i^2)$ }\\
&\qquad\qquad\qquad\qquad\qquad\qquad\qquad\qquad
\quad\! \mbox{thus $(s -t,x)\in Q_i'(y)$}. 
\end{aligned}
\end{align}
Hence, from \eqref{dllkka} we derive  
\begin{align*}
\|w\|_{L^2H^{1+\alpha}(Q_i')}
&\leq\sup_{y}\|G(\cdot,\cdot,y)\|_{L^2H^{1+\alpha}(Q_i'(y))}
\|v\|_{L^1(Q_j)}\\
&\leq Cd_i\, \sup_{y}\|G(\cdot,\cdot,y)\|_{L^\infty H^{1+\alpha}(Q_i'(y))}
\|v\|_{L^1(Q_j)}\\
&\leq
Cd_i^{-\alpha-N/2 }d_j^{1+N/2}\vertiii{v}_{Q_j}
\qquad\mbox{(here we use \eqref{GFest03})}\\
&\leq
Cd_j^{1-\alpha}\bigg(\frac{d_j}{d_i} \bigg)^{\alpha+N/2}   .
\end{align*}
For $i\geq j+3$ ($d_i\leq d_j$, including the case $i=*$), we have 
\begin{align}
\begin{aligned}
&\mbox{if $(t,x)\in Q_i'$ and $(s ,y)\in Q_j$, then $\max(|s -t|^{1/2},|x-y|)\geq d_{j+3}$,}\\
&\qquad\qquad\qquad\qquad\qquad\qquad
\quad \mbox{thus $(s -t,x)\in \bigcup_{k\le j+3}Q_k(y)$}, \\
&\mbox{if $(t,x)\in Q_*$ and $(s ,y)\in Q_j$ with $j\le J_*-3$, then $\max(|s -t|^{1/2},|x-y|)\geq d_{j+3}$,}\\
&\qquad\qquad\qquad\qquad\qquad\qquad\qquad\qquad\qquad\,
\qquad \mbox{thus $(s -t,x)\in \bigcup_{k\le j+3}Q_k(y)$}. 
\end{aligned}
\end{align}
If $\widetilde G(\cdot,\cdot,y)$ is a function satisfying 
\begin{align} \label{Exend-G}
\widetilde G(\cdot,\cdot,y)=G(\cdot,\cdot,y)\quad\mbox{on}\quad \bigcup_{k\le j+3}Q_k(y) ,
\end{align}
then for $(t,x)\in Q_i'$ we have 
$$
w(t,x)
=\int_0^{T}\int_{\Omega}
G(s-t,x,y)v(y,s)1_{s>t}\,\d y\d s 
=\iint_{Q_j}
\widetilde G(s-t,x,y)v(s,y)1_{s>t}\,\d y\d s  .
$$
Hence, we have
\begin{align*}
\|w\|_{L^2H^{1+\alpha}(Q_i')}
&\le C d_i\|w\|_{L^\infty H^{1+\alpha}(Q_i')} \le C d_i\|w\|_{L^\infty H^{1+\alpha}({\cal Q})} \\
&\le C d_i \iint_{Q_j}
\|\widetilde G(s-t,x,y)1_{s>t}\|_{L^\infty_tH^{1+\alpha}_x({\cal Q})}  |v(s,y)|\d y\d s \\
&\le C d_i \iint_{Q_j}
\|\widetilde G(\cdot,\cdot,y)\|_{L^\infty H^{1+\alpha}({\cal Q})}  |v(s,y)|\d y\d s \\
&\le C d_i \sup_{y\in\Omega}\|\widetilde G(\cdot,\cdot,y)\|_{L^\infty H^{1+\alpha}({\cal Q})} 
\|v\|_{L^1(Q_j)}\\ 
&\leq Cd_i \sup_{y\in\Omega}\|\widetilde G(\cdot,\cdot,y)\|_{L^\infty H^{1+\alpha}({\cal Q})} 
\vertiii{v}_{Q_j}d_j^{1+N/2}  . 
\end{align*} 
In view of the definition \eqref{DefLpX}, by taking infimum over all the possible choices of $\widetilde G(\cdot,\cdot,y)$ satisfying \eqref{Exend-G}, we have 
\begin{align*}
\|w\|_{L^2H^{1+\alpha}(Q_i')}
&\leq Cd_i \sup_{y\in\Omega}\|G(\cdot,\cdot,y)\|_{L^\infty H^{1+\alpha}(\cup_{k\le j+3}Q_k(y))} 
\vertiii{v}_{Q_j}d_j^{1+N/2}  .\\ 
&\leq Cd_i d_j^{-1-\alpha-N/2}d_j^{N/2+1}
\qquad\qquad\qquad \mbox{(here we use \eqref{GFest03})}\\[5pt] 
&=Cd_i^{1-\alpha}\bigg(\frac{d_i}{d_j} \bigg)^{\alpha} . 
\end{align*} 

For $|i-j|\leq 2$, applying the standard energy estimate yields
$$ 
\|w\|_{L^2H^{1+\alpha}(Q_i')}
\le  
\|w\|_{L^2H^{1+\alpha}({\cal Q})}
\le 
C\|v\|_{L^2H^{\alpha-1}({\cal Q})}
\le 
C\|v\|_{L^2L^{\frac{N}{1+N/2-\alpha}}({\cal Q})}
\le 
Cd_i^{1-\alpha}\vertiii{v}_{{\cal Q}}=Cd_i^{1-\alpha} ,
$$ 
where we have used the Sobolev embedding $L^{\frac{N}{1+N/2-\alpha}}\hookrightarrow H^{\alpha-1}$ and 
the H\"older's inequality 
$$\|v\|_{L^{\frac{N}{1+N/2-\alpha}}}\le Cd_j^{1-\alpha} \|v\|_{L^2}
\quad\mbox{(this requires the volume of the support of $v$ to be bounded by $d_j^N$)} .
$$ 
 
Combining the three cases above (corresponding to $i\le j-3$, $i\ge j+3$ and $|i-j|\le 2$), we have  
\begin{align}\label{wW2pQi} 
\|w\|_{L^2H^{1+\alpha}(Q_i')} 
\leq 
Cd_i^{1-\alpha} 
\bigg(\frac{\min (d_i ,d_j )}{\max (d_i ,d_j )}\bigg)^{\alpha} . 
\end{align} 
Substituting \eqref{f69}-\eqref{sd80} and \eqref{wW2pQi} into (\ref{dka6}) yields 
\begin{align}\label{FL2Qj} 
\vertiii{F}_{Q_j} 
&\leq Ch^{1+\alpha-N/2}d_j^{-\alpha} 
+C\sum_{*,i}(h^{1+\alpha}\vertiii{F_t}_{Q_i} 
+ h^{\alpha} \vertiii{F}_{1,Q_i})
d_i^{1-\alpha} \bigg(\frac{\min (d_i ,d_j )}{\max (d_i ,d_j )}\bigg)^{\alpha} .
\end{align}
Since $\alpha>1/2$, it follows that 
\begin{align}\label{lalnu}
\sum_jd_j^{N/2-1}
\bigg(\frac{\min (d_i ,d_j )}{\max (d_i ,d_j )}\bigg)^{\alpha}
\leq Cd_i^{N/2-1} .
\end{align}
Hence, we have 
\begin{align*}
{\mathscr K} &\leq 
C+CC_*^{3+N/2}+
C\sum_{j}d_j^{-1+N/2}\vertiii{F}_{Q_j} \nn\\
&\leq C+CC_*^{3+N/2}
+C\sum_j \left(\frac{h}{d_j}\right)^{1+\alpha-N/2} 
\qquad\mbox{here we substitute \eqref{FL2Qj} }\\
&\quad
+C\sum_j d_j^{-1+N/2}
  \sum_{*,i}(h^{1+\alpha}\vertiii{F_t}_{Q_i} 
+ h^{\alpha} \vertiii{F}_{1,Q_i})
d_i^{1-\alpha} \bigg(\frac{\min (d_i ,d_j )}{\max (d_i ,d_j )}\bigg)^{\alpha}\\
&\leq
C+CC_*^{3+N/2} +C 
\qquad\mbox{here we exchange the order of summation} \\
&\quad 
+C\sum_{*,i}(h^{1+\alpha}\vertiii{F_t}_{Q_i} 
+ h^{\alpha} \vertiii{F}_{1,Q_i})d_i^{1-\alpha}
\sum_jd_j^{N/2-1}
\bigg(\frac{\min (d_i ,d_j )}{\max (d_i ,d_j )}\bigg)^{\alpha}\\
&\leq
C+CC_*^{3+N/2}    
+C\sum_{*,i}(h^{1+\alpha}\vertiii{F_t}_{Q_i} 
+ h^{\alpha} \vertiii{F}_{1,Q_i})d_i^{N/2-\alpha} 
\qquad\mbox{(here we use \eqref{lalnu})}\\
&=
C+CC_*^{3+N/2}    
+C\sum_{*,i}d_i^{1+N/2}\left(\vertiii{F_t}_{Q_i} 
\bigg(\frac{h}{d_i}\bigg)^{1+\alpha} 
+ d_i^{-1}\vertiii{F}_{1,Q_i}\bigg(\frac{h}{d_i}\bigg)^{\alpha}\right) \\
&\leq 
C+CC_*^{3+N/2}    
+Cd_{*}^{1+N/2}
\left(\vertiii{F_t}_{Q_*} + d_j^{-1}\vertiii{F}_{1,Q_*} \right) \\
&\quad
+C\sum_{i}d_i^{1+N/2}\left(\vertiii{F_t}_{Q_i} 
+ d_j^{-1}\vertiii{F}_{1,Q_i}\right)\bigg(\frac{h}{d_i}\bigg)^{\alpha} \\
&\leq  C+C C^{3+N/2 }_* +\frac{C {\mathscr K} }{C_*^{\alpha}}  .
\end{align*}
By choosing $C_*$ to be large enough ($C_*$ is determined now), the term $\displaystyle\frac{C {\mathscr K} }{C_*^{\alpha}} $ will be absorbed by the left-hand side of the inequality above. In this case, the inequality above implies 
\begin{align}
{\mathscr K} \le C .
\end{align}  
Substituting the last inequality into \eqref{Bd31K2} yields 
\begin{align}\label{L1FtQ1}
\|\partial_tF\|_{L^1({\cal Q})}+\|t\partial_{tt}F\|_{L^1({\cal Q})} 
&\leq C .
\end{align} 

{\it Part II.}$\,\,\,$ 
Second, we present estimates for $(t,x)\in(1,\infty)\times\Omega$. 
For $t>1$, 
we differentiate \eqref{EqGammh} with respect to $t$ and
integrate the resulting equation against
$\partial_t\Gamma_h$. Then we get
\begin{align*}
&\frac{\d}{\d t}\|\partial_t\Gamma_h(t,\cdot, x_0)\|_{L^2}^2
+\lambda_0\|\partial_t\Gamma_h(t,\cdot, x_0)\|_{L^2}^2  \\
&\leq \frac{\d}{\d t}\|\partial_t\Gamma_h(t,\cdot, x_0)\|_{L^2}^2
+(\nabla\partial_t\Gamma_h(t,\cdot, x_0),
\nabla\partial_t\Gamma_h(t,\cdot, x_0)) \\
&= 0 ,
\end{align*}
for $t\geq 1$, where
$\lambda_0>0$ is the smallest eigenvalue of the operator $-\Delta$. 
From the last inequality we derive the exponential decay of $\partial_t\Gamma_h$ with respect to $t$
\begin{align*}
\|\partial_t\Gamma_h(t,\cdot, x_0)\|_{L^2}^2\leq
e^{-\lambda_0(t-1)}\|\partial_t\Gamma_h(1,\cdot, x_0)\|_{L^2}^2
\leq Ce^{-\lambda_0 (t-1)} ,
\end{align*}
where the inequality $\|\partial_t\Gamma_h(1,\cdot, x_0)\|_{L^2}\leq C$
can be proved by a simple energy estimate (omitted here). 
Similarly, we also have 
\begin{align*}
\|\partial_{tt}\Gamma_h(t,\cdot, x_0)\|_{L^2}^2
+\|\partial_{t}\Gamma(t,\cdot, x_0)\|_{L^2}^2
+\|\partial_{tt}\Gamma(t,\cdot, x_0)\|_{L^2}^2
\leq Ce^{-\lambda_0(t-1)} \quad\mbox{for}\,\,\, t\geq 1. 
\end{align*}
The estimate \eqref{L1FtQ1} and the last two inequalities imply \eqref{L1Ft}-\eqref{L1Gammatx0} in the case $h<h_*:=1/(4C_*)$. 

For $h\geq h_*$, some basic energy estimates would yield 
\begin{align*}
&\|\partial_{t}\Gamma_h(t,\cdot, x_0)\|_{L^2}^2
+\|\partial_{tt}\Gamma_h(t,\cdot, x_0)\|_{L^2}^2
+\|\partial_{t}\Gamma(t,\cdot, x_0)\|_{L^2}^2
+\|\partial_{tt}\Gamma(t,\cdot, x_0)\|_{L^2}^2 \\
&\leq Ce^{-\lambda_0t}
(\|\partial_{t}\Gamma_h(0,\cdot, x_0)\|_{L^2}^2
+\|\partial_{tt}\Gamma_h(0,\cdot, x_0)\|_{L^2}^2
+\|\partial_{t}\Gamma(0,\cdot, x_0)\|_{L^2}^2
+\|\partial_{tt}\Gamma(0,\cdot, x_0)\|_{L^2}^2)  \\
&= Ce^{-\lambda_0t}
(\|\Delta_hP_h\widetilde\delta_{x_0} \|_{L^2}^2
+\|\Delta_h^2P_h\widetilde\delta_{x_0} \|_{L^2}^2
+\|\Delta\widetilde\delta_{x_0}\|_{L^2}^2
+\|\Delta^2\widetilde\delta_{x_0} \|_{L^2}^2)  \\
&\leq Ce^{-\lambda_0t} (h_*^{-4-N}+h_*^{-8-N}) 
\end{align*}
for arbitrary $t>0$.  This implies \eqref{L1Ft}-\eqref{L1Gammatx0} in the case $h\ge h_*$. 
\medskip

{\it Part III.}$\,\,\,$ 
Finally, we note that \eqref{L1Gammh-2} is a simple consequence of  \eqref{reg-Delta-est}, \eqref{GausEst1} and \eqref{expr-Gamma}, while \eqref{L1Gammh} is a consequence of \eqref{L1Gammh-2} and the following inequalities: 
\begin{align}
&\|\Gamma_h(t,\cdot, x_0)\|_{L^1(\Omega)}\le \|\Gamma(t,\cdot, x_0)\|_{L^1(\Omega)}+\|F(t,\cdot, x_0)\|_{L^1(\Omega)},\\
&\|t\partial_t\Gamma_h(t,\cdot, x_0)\|_{L^1(\Omega)}\le \|t\partial_t\Gamma(t,\cdot, x_0)\|_{L^1(\Omega)}+\|t\partial_t F(t,\cdot, x_0)\|_{L^1(\Omega)},
\end{align} 
with 
\begin{align}
\|F(t,\cdot, x_0)\|_{L^1(\Omega)}
&\le \|F(0,\cdot, x_0)\|_{L^1(\Omega)}
+\bigg\|\int_0^t\partial_s  F(s ,\cdot, x_0) \d s \bigg\|_{L^1(\Omega)} \nn\\
&\le \|P_h\widetilde\delta_{x_0} 
-\widetilde\delta_{x_0}\|_{L^1(\Omega)}
+\|\partial_tF\|_{L^1({\cal Q})}  
\leq C ,\\[5pt]
\|t\partial_t F(t,\cdot, x_0)\|_{L^1(\Omega)}
&\le \bigg\|\int_0^t \Big(s \partial_{s s }F(s ,\cdot, x_0)
+\partial_{s }F_h(s ,\cdot, x_0)\Big) \d s \bigg\|_{L^1(\Omega)} \nn\\
&\leq \|t\partial_{tt}F\|_{L^1({\cal Q})}  +\|\partial_tF\|_{L^1({\cal Q})}   
\le C ,
\end{align} 
where we have used \eqref{L1Ft} in the last two inequalities, which was proved in Part I and Part II.\medskip   

The proof of Lemma \ref{LemGm2} is complete.
\qed

\section{Conclusion}
The analyticity and maximal $L^p$-regularity of finite element solutions of the heat equation are proved in general polygons and polyhedra, possibly nonconvex. The $L^\infty$-stability of the finite element parabolic projection has been reduced to the $L^\infty$-stability of the Ritz projection. Such $L^\infty$-stability of the Dirichlet Ritz projection is currently known in general polygons \cite{Schatz1980} and convex polyhedra \cite{LeykekhmanVexler2016}, but still remains open in nonconvex polyhedra. The $L^\infty$-stability of the Neumann Ritz projection remains open in both nonconvex polygons and nonconvex polyhedra. This article focuses on the Lagrange finite element method. Extension of the results to other numerical methods, such as finite volume methods and discontinuous Galerkin methods, are interesting and nontrivial. Such extension may need more precise $W^{s,p}$-approximation properties of local elliptic projectors onto finite element spaces (e.g., see \cite{DiPietroDroniou2017}).

\section*{Appendix A: Proof of Lemma \ref{LocEEst}}
\label{LocEngPrf}

\renewcommand{\thelemma}{A.\arabic{lemma}}
\renewcommand{\theequation}{A.\arabic{equation}}
\setcounter{equation}{0}\setcounter{lemma}{0}

In this subsection, we prove Lemma \ref{LocEEst}, which is used in the last section in proving Lemma \ref{LemGm2}. Before we prove Lemma \ref{LocEEst}, we present a local energy estimate for finite element solutions of parabolic equations based on the decomposition $Q_j=[(0,d_j^2) \times \Omega_j]
\cup [(d_j^2,4d_j^2)\times D_j]$. 
\begin{lemma}\label{lemlocDEng0}
{\it
Suppose that $\phi_h(t)\in S_h$, $t\in(0,T)$, satisfies 
\begin{align*} 
&\big(\partial_t\phi_h ,\chi_h\big)
+\big(\nabla \phi_h ,\nabla\chi_h\big) =0, 
\qquad\mbox{for}\,\,\,\,\chi_h\in S_h^0(\Omega_j''),
\,\,\, t\in(0,d_j^2)  ,\\
&\big(\partial_t \phi_h ,\chi_h\big)
+\big(\nabla \phi_h ,\nabla\chi_h\big)
=0   ,\qquad\mbox{for}\,\,\,\,\chi_h\in S_h^0(D_j'') ,
\,\,\, t\in(d_j^2/4,2d_j^2) .
\end{align*}
Then we have 
\begin{align}\label{locDEng0}
&\|\partial_t\phi_h\|_{L^2(Q_j)} +d_j^{-1}\|\nabla \phi_h\|_{L^2(Q_j)} \nn\\
&\leq  
 (Ch^{1/2}d_j^{-1/2}+C\epsilon^{-1} hd_j^{-1}+\epsilon)\big(
\|\partial_t\phi_h\|_{L^2(Q_{j}')}   
+ d_j^{-1}\|\nabla\phi_h\|_{L^2(Q_{j}')}\big) \nn\\
&\quad 
+C\epsilon^{-1}\left(d_j^{-2} \|\phi_h\|_{L^2(Q_{j}')}
+d_j^{-1}\|\phi_h(0)\|_{L^2(\Omega_j')} 
+\|\nabla \phi_h(0)\|_{L^2(\Omega_j')}  \right)  ,
\end{align} 
where the constant $C$ is independent of $h$, $j$ and $C_*$. 
}
\end{lemma} 

\noindent{\it Proof of Lemma \ref{lemlocDEng0}.}$\,\,\,$
We shall present estimates in the 
two subdomains $(0,d_j^2) \times \Omega_j$
and $(d_j^2,4d_j^2)\times D_j$, separately.

First, we present estimates in $(0,d_j^2) \times \Omega_j$. 
To this end, we let $\omega$ be a smooth cut-off function 
which equals $1$ on $\Omega_j$ 
and equals $0$ outside $\Omega_j'$,
and let $\widetilde\omega$ be a smooth cut-off function 
which equals $1$ on $\Omega_j'$ 
and equals $0$ outside $\Omega_j''$, such that   

%
(1) ${\rm dist}\big({\rm supp}(\omega)\cap \Omega,
\Omega\backslash\Omega_j'\big)\geq 
d_j/8\geq 2 h$ 
and ${\rm dist}\big({\rm supp}(\widetilde\omega)\cap \Omega,
\Omega\backslash\Omega_j''\big)
\geq d_j/8\geq 2 h$,

(2) $ |\partial^\alpha\omega|
+|\partial^\alpha\widetilde\omega|\leq C_\alpha d_j^{-|\alpha|}$
for any multi-index $\alpha$. 
 
By Property (P3) of Section \ref{Sec2-2}, the function 
$v_h:=I_h(\widetilde\omega\phi_h)\in S_h^0(\Omega_j'')$
satisfies $v_h=\phi_h$ on $\Omega_j'$ and 
\begin{align} 
&\|v_h\|_{L^2(\Omega)}
\leq C\|\phi_h\|_{L^2(\Omega_j'')}  , \label{vhphih1}\\
&\|\nabla v_h\|_{L^2(\Omega)}
\leq C\|\nabla \phi_h\|_{L^2(\Omega_j'')}
+ Cd_j^{-1}\| \phi_h\|_{L^2(\Omega_j'')} ,  \label{vhphih2} \\[5pt]
&
\big(\partial_tv_h,\chi_h\big)
+\big(\nabla v_h,\nabla\chi_h\big)
=0  ,\quad\forall\, \chi_h\in S_h^0(\Omega_j') \, ,
\quad\forall\,\, t\in (0,d_j^2) .
\end{align}
Property (P3) of Section \ref{Sec2-2} also implies that $I_h(\omega^2v_h)\in S_h^0(\Omega_j')$
such that
\begin{align*}
&\|\omega^2v_h-I_h(\omega^2v_h)\|_{L^2}
+h\|\nabla (\omega^2v_h-I_h(\omega^2v_h))\|_{L^2}
\leq Chd_j^{-1} \|v_h\|_{L^2}. 
\end{align*}
Since $\omega$ and $\widetilde\omega$ are time-independent, it follows that
\begin{align*}
&\frac{1}{2}\frac{\d}{\d t}\|\omega v_h\|^2
+(\omega^2 \nabla v_h,\nabla v_h)\\
&=\big[\big(\partial_tv_h,\omega^2v_h\big)
+\big(\nabla v_h,\nabla(\omega^2v_h)\big)\big]  
-\big(2 v_h\omega\nabla\omega , \nabla v_h\big)  \\
&=\big[\big(\partial_tv_h,\omega^2v_h-I_h(\omega^2v_h)h\big)
+\big(\nabla v_h,\nabla(\omega^2v_h -I_h(\omega^2v_h))\big)\big]  
-\big(2v_h\nabla\omega , \omega \nabla v_h\big)  \\
&\leq 
\big[ C\|\partial_tv_h\|_{L^2}\|v_h\|_{L^2}hd_j^{-1}
+ C\|\nabla v_h\|_{L^2}\|v_h\|_{L^2}d_j^{-1}  \big]
+Cd_j^{-1}\|v_h\|_{L^2}\|\omega\nabla v_h\|_{L^2} \\
&\leq 
C\|\partial_tv_h\|_{L^2}^2h^2 
+ \epsilon^4\|\nabla v_h\|_{L^2}^2
+C\epsilon^{-4}\|v_h\|_{L^2}^2d_j^{-2} ,
\qquad\forall\,\,\epsilon\in(0,1) .
\end{align*} 
By using \eqref{vhphih1}-\eqref{vhphih2}, integrating the last inequality from $0$ to $d^2_j$ yields
\begin{align}\label{locEng1}
&\|\phi_h\|_{L^\infty(0,d_j^2; L^2(\Omega_j))} 
+\|\nabla \phi_h\|_{L^2(0,d_j^2; L^2(\Omega_j))} \nn\\
&\leq C\|\phi_h(0)\|_{L^2(\Omega_j'')}
+C\|\partial_t\phi_h\|_{L^2(0,d_j^2; L^2(\Omega_j''))}h \nn \\
&\quad
+\epsilon^2\|\nabla \phi_h\|_{L^2(0,d_j^2; L^2(\Omega_j''))}  
+C\epsilon^{-2}\|\phi_h\|_{L^2(0,d_j^2; L^2(\Omega_j''))} d_j^{-1} . 
\end{align}

Furthermore, we have
\begin{align*} 
&\|\omega^2\partial_tv_h\|_{L^2}^2
+\frac{1}{2}\frac{\d}{\d t}\big(\omega^4 \nabla v_h,\nabla v_h\big)\\
&=\big(\partial_tv_h,\omega^4\partial_tv_h \big)
+\big(\nabla v_h,\nabla(\omega^4\partial_tv_h ) \big)   -\big(4\partial_tv_h\omega^3\nabla \omega,\nabla v_h \big) \\
&=\big(\partial_tv_h,\omega^4\partial_tv_h-I_h(\omega^4\partial_tv_h)\big)
+\big(\nabla v_h,\nabla(\omega^4\partial_tv_h -I_h(\omega^4\partial_tv_h))\, \big)   -\big(4\partial_tv_h\omega^3\nabla \omega,\nabla v_h \big) \\
&\le C\|\partial_tv_h\|_{L^2}^2hd_j^{-1}
+C\|\nabla v_h\|_{L^2}\|\partial_tv_h\|_{L^2}d_j^{-1}  
+ C\|\omega^2\partial_tv_h\|_{L^2}
\|\nabla v_h\|_{L^2} d_j^{-1} \\
&\le (Chd_j^{-1}+\epsilon^2)\|\partial_tv_h\|_{L^2}^2
+\epsilon^2\|\omega^2\partial_tv_h\|_{L^2}^2
+ C\epsilon^{-2}\|\nabla v_h\|_{L^2}^2 d_j^{-2} , 
\qquad\forall\,\,\epsilon\in(0,1/2) , 
\end{align*} 
which reduces to 
\begin{align*} 
\|\omega^2\partial_tv_h\|_{L^2(0,d_j^2;L^2)}^2 
&\leq (Ch d_j^{-1}+\epsilon^2)\|\partial_tv_h\|_{L^2(0,d_j^2;L^2)}^2 
+ C\epsilon^{-2}\|\nabla v_h\|_{L^2(0,d_j^2;L^2)}^2 d_j^{-2} 
+ C\|\nabla v_h(0)\|_{L^2}^2 .
\end{align*}
By using \eqref{vhphih1}-\eqref{vhphih2}, 
the last inequality further implies 
\begin{align}\label{locEng2}
\|\partial_t\phi_h\|_{L^2(0,d_j^2;L^2(\Omega_j))}  
&\leq C(\|\nabla\phi_h(0)\|_{L^2(\Omega_j'')} 
+d_j^{-1}\|\phi_h(0)\|_{L^2(\Omega_j'')}) \nn\\
&\quad 
+  (Ch^{1/2}d_j^{-1/2}+\epsilon)\|\partial_t\phi_h\|_{L^2(0,d_j^2;L^2(\Omega_j''))} 
\nn\\
&\quad
+ C\epsilon^{-1}( \|\phi_h\|_{L^2(0,d_j^2;L^2(\Omega_j''))}  d_j^{-2}
+\|\nabla \phi_h\|_{L^2(0,d_j^2;L^2(\Omega_j''))}  d_j^{-1} )  \nn\\
&\leq C(\|\nabla\phi_h(0)\|_{L^2(\Omega_j'')}
+d_j^{-1}\|\phi_h(0)\|_{L^2(\Omega_j'')})  \nn\\
&\quad 
+ (Ch^{1/2}d_j^{-1/2}+\epsilon)\|\partial_t\phi_h\|_{L^2(0,d_j^2;L^2(\Omega_j''))} 
\nn\\
&\quad
+ C\epsilon^{-1}(\|\phi_h\|_{L^\infty(0,d_j^2; L^2(\Omega_j''))}  d_j^{-1}
+\|\nabla \phi_h\|_{L^2(0,d_j^2;L^2(\Omega_j''))}  d_j^{-1} )  .
\end{align}
With an obvious change of domains (replacing $\Omega_j$ by $\Omega_j'$ on the left-hand side of \eqref{locEng1} and replacing $\Omega_j''$ by $\Omega_j'$ on the right-hand side of \eqref{locEng2}), 
the two estimates \eqref{locEng1} and \eqref{locEng2} imply  
\begin{align}\label{Ineq-A7}
&\|\phi_h\|_{L^\infty(0,d_j^2;L^2(\Omega_j'))} 
+\|\nabla \phi_h\|_{L^2(0,d_j^2;L^2(\Omega_j'))} \nn\\
&\leq \hat C\|\phi_h(0)\|_{L^2(\Omega_j'')}
+\hat C\|\partial_t\phi_h\|_{L^2(0,d_j^2; L^2(\Omega_j''))}h \nn \\
&\quad
+\epsilon^2\|\nabla \phi_h\|_{L^2(0,d_j^2; L^2(\Omega_j''))}  
+\hat C\epsilon^{-2}\|\phi_h\|_{L^2(0,d_j^2; L^2(\Omega_j''))} d_j^{-1} . 
\end{align} 
and 
\begin{align} \label{Ineq-A8}
\|\partial_t\phi_h\|_{L^2(0,d_j^2;L^2(\Omega_j))}  
&\leq \hat C(\|\nabla\phi_h(0)\|_{L^2(\Omega_j')}
+d_j^{-1}\|\phi_h(0)\|_{L^2(\Omega_j')})  \nn\\
&\quad 
+ (\hat Ch^{1/2}d_j^{-1/2}+\epsilon)\|\partial_t\phi_h\|_{L^2(0,d_j^2;L^2(\Omega_j'))}  
\nn\\
&\quad 
+ \hat C\epsilon^{-1}d_j^{-1}  (\|\phi_h\|_{L^\infty(0,d_j^2;L^2(\Omega_j'))}  
+\|\nabla \phi_h\|_{L^2(0,d_j^2;L^2(\Omega_j'))}  ) ,
\end{align} 
where $\hat C\ge 1$ is some positive constant and $\epsilon\in(0,1)$ can be arbitrary. 
Then $2\hat C\epsilon^{-1}d_j^{-1}\times$\eqref{Ineq-A7}+\eqref{Ineq-A8} yields (the last term in \eqref{Ineq-A8} can be absorbed by left-hand side of $2\hat C\epsilon^{-1}d_j^{-1}\times$\eqref{Ineq-A7}) 
\begin{align} \label{phihtxOj}
&\|\partial_t\phi_h\|_{L^2(0,d_j^2;L^2(\Omega_j))}  
+d_j^{-1}\|\nabla \phi_h\|_{L^2(0,d_j^2;L^2(\Omega_j))} \nn\\
&\leq (Ch^{1/2}d_j^{-1/2} +C\epsilon^{-1}hd_j^{-1}+\epsilon)
(\|\partial_t\phi_h\|_{L^2(0,d_j^2;L^2(\Omega_j''))}  
+d_j^{-1}\|\nabla \phi_h\|_{L^2(0,d_j^2;L^2(\Omega_j''))}) \nn\\
&\quad 
+ C\epsilon^{-3}\|\phi_h\|_{L^2(0,d_j^2;L^2(\Omega_j''))} d_j^{-2}
+ C\epsilon^{-1}(\|\nabla\phi_h(0)\|_{L^2(\Omega_j'')}
+d_j^{-1}\|\phi_h(0)\|_{L^2(\Omega_j'')})   .
\end{align}

Second, we present estimates in 
$(d_j^2,4d_j^2)\times D_j$.  
We re-define $\omega(x,t):=\omega_1(x)\omega_2(t)$ 
and $\widetilde\omega(x,t):=
\widetilde\omega_1(x) \widetilde\omega_2(t)$ 
such that 

(1) $\omega_1=1$ in $D_j$ and  $\omega_1=0$ outside $D_j'$, 
$\widetilde\omega_1=1$ in $D_j'$ 
and $\widetilde\omega_1=0$ outside $D_j''$;

(2) ${\rm dist}({\rm supp}(\omega_1)\cap\Omega,
\Omega\backslash D_j')\geq d_j/4\geq\kappa h$
and 
${\rm dist}({\rm supp}(\widetilde\omega_1)\cap\Omega,
\Omega\backslash D_j'')\geq d_j/4\geq\kappa h$;

(3) $\omega_2=1$ for $t\in(d_j^2,4d_j^2)$
and $\omega_2=0$ for $t\in(0,d_j^2/2)$;

(4) $\widetilde\omega_2=1$ for $t\in(d_j^2/4,4d_j^2)$ and
$\widetilde\omega_2=0$ for $t\in(0,d_j^2/8)$; 

(5) $|\partial^\alpha\omega_1|
+|\partial^\alpha\widetilde\omega_1|\leq Cd_j^{-|\alpha|}$
for any multi-index $\alpha$;

(6) $|\partial_t^k\omega_2|+|\partial_t^k\widetilde\omega_2|\leq Cd_j^{-2k}$
for any nonnegative integer $k$. 

Then the function 
$v_h:=I_h(\widetilde\omega\phi_h)\in S_h^0(D_j'')$
satisfies $v_h=\widetilde\omega_2\phi_h$ on $D_j'$ and 
\begin{align} 
&\|v_h\|_{L^2(\Omega)}
\leq C\|\widetilde\omega_2\phi_h\|_{L^2(D_j'')}  , \label{vhphih1-2}\\
&\|\nabla v_h\|_{L^2(\Omega)}
\leq C\|\widetilde\omega_2\nabla \phi_h\|_{L^2(D_j'')}
+ Cd_j^{-1}\| \widetilde\omega_2\phi_h\|_{L^2(D_j'')} ,  \label{vhphih2-2} \\ 
&
\big(\partial_tv_h,\chi_h\big)
+\big(\nabla v_h,\nabla\chi_h\big)
=0  ,\quad\forall\, \chi_h\in S_h^0(D_j') \quad\forall\,\,
t\in(d_j^2/4,4d_j^2) .
\end{align} 
According to (P3) of Section \ref{Sec2-2},  
the function $\chi_h=I_h(\omega^2v_h)\in S_h^0(D_j')$ satisfies
\begin{align*}
&\|\omega^2v_h-\chi_h\|_{L^2}
+h\|\nabla (\omega^2v_h- \chi_h)\|_{L^2}
\leq Chd_j^{-1} \|v_h\|_{L^2}. 
\end{align*}
Therefore,  we have 
\begin{align*}
&\frac{1}{2}\frac{\d}{\d t}\|\omega v_h\|^2
+(\omega^2 \nabla v_h,\nabla v_h)\\
&=\big[\big(\partial_t v_h,\omega^2 v_h\big)
+\big(\nabla v_h,\nabla(\omega^2 v_h)\big)\big]  
+\big(\partial_t\omega v_h, v_h\big)
-\big(2 v_h\omega\nabla \omega , \nabla v_h\big)  \\
&=\big[\big(\partial_t v_h, \omega^2 v_h-\chi_h\big)
+\big(\nabla v_h,\nabla( \omega^2 v_h -\chi_h)\big)\big]  
+\big(\omega\partial_t\omega v_h,v_h\big)
-\big(2 v_h\nabla \omega ,\omega \nabla v_h\big)  \\
&\leq 
\big[ C\|\partial_t v_h\|_{L^2(\Omega)}\| v_h\|_{L^2(\Omega)}hd_j^{-1}
+ C\|\nabla v_h\|_{L^2(\Omega)}\| v_h\|_{L^2(\Omega)}d_j^{-1}  \big] \\
&\quad +Cd_j^{-2}\| v_h\|_{L^2(\Omega)}^2
+Cd_j^{-1}\| v_h\|_{L^2(\Omega)}\|\omega\nabla v_h\|_{L^2(\Omega)} \\
&\leq 
C\|\partial_tv_h\|_{L^2(\Omega)}^2h^2 
+ \epsilon^4\|\nabla v_h\|_{L^2(\Omega)}^2
+ \epsilon^4\|\omega\nabla v_h\|_{L^2(\Omega)}^2
+C\epsilon^{-4}\|v_h\|_{L^2(\Omega)}^2d_j^{-2} ,
\quad\forall\,\,\epsilon\in(0,1/2) .
\end{align*} 
Integrating the last inequality in time for $t\in(d_j^2/2,4d_j^2)$, we obtain 
\begin{align}\label{phihDj}
&\|\phi_h\|_{L^\infty(d_j^2,4d_j^2; L^2(D_j))} 
+\|\nabla \phi_h\|_{L^2(d_j^2,4d_j^2; L^2(D_j))} \\
&\leq C\|\partial_t\phi_h\|_{L^2(d_j^2/2,4d_j^2; L^2(D_j''))}h
+\epsilon^2\|\nabla \phi_h\|_{L^2(d_j^2/2,4d_j^2; L^2(D_j''))}  
+C\epsilon^{-2}\|\phi_h\|_{L^2(d_j^2/2,4d_j^2; L^2(D_j''))} d_j^{-1} . 
\nn
\end{align}
Furthermore, we have
\begin{align*} 
&\|\omega^2\partial_tv_h\|_{L^2}^2
+\frac{1}{2}\frac{\d}{\d t}\big(\omega^4 \nabla v_h,\nabla v_h\big)\\
&=\big[\big(\partial_tv_h,\omega^4\partial_tv_h-\eta_h\big)
+\big(\nabla v_h,\nabla(\omega^4\partial_tv_h -\eta_h)\, \big)\big]  \\
&\quad +\big(2\omega^3\partial_t\omega \nabla v_h,\nabla v_h\big) 
-\big(4\partial_tv_h\omega^3\nabla \omega,\nabla v_h \big) \\
&\leq [C\|\partial_tv_h\|_{L^2(\Omega)}^2hd_j^{-1}
+C\|\nabla v_h\|_{L^2(\Omega)}\|\partial_tv_h\|_{L^2(\Omega)}d_j^{-1} \\
&\quad + C\|\nabla v_h\|_{L^2(\Omega)}^2 d_j^{-2} 
+ C\|\omega^2\partial_tv_h\|_{L^2(\Omega)}
\|\nabla v_h\|_{L^2(\Omega)} d_j^{-1} \\
&\leq (Chd_j^{-1}+\epsilon)\|\partial_tv_h\|_{L^2(\Omega)}^2
+\epsilon^2\|\omega^2\partial_tv_h\|_{L^2(\Omega)}^2
+ C\epsilon^{-2}\|\nabla v_h\|_{L^2(\Omega)}^2 d_j^{-2} , 
\quad\forall\,\,\epsilon\in(0,1/2) , 
\end{align*} 
which implies (by integrating the last inequality in time for $t\in(d_j^2/2, 4d_j^2)$)
\begin{align*} 
\|\omega^2\partial_tv_h\|_{L^2(d_j^2,4d_j^2;L^2(\Omega))}^2 
\leq (Chd_j^{-1}+\epsilon^2) \|\partial_tv_h\|_{L^2(d_j^2/2,4d_j^2;L^2(\Omega))}^2
+ C\epsilon^{-2}\|\nabla v_h\|_{L^2(d_j^2/2,4d_j^2;L^2(\Omega))}^2 d_j^{-2}  .
\end{align*}
By using \eqref{vhphih1-2}-\eqref{vhphih2-2}, 
the last inequality further implies 
\begin{align}\label{phihtDj}
&\|\partial_t\phi_h\|_{L^2(d_j^2,4d_j^2;L^2(D_j))}   \nn\\ 
&\leq  
(Ch^{1/2}d_j^{-1/2}+\epsilon) \|\partial_t\phi_h\|_{L^2(d_j^2/2,4d_j^2;L^2(D_j''))} \nn\\
&\quad 
+ C\epsilon^{-1}( \|\phi_h\|_{L^2(d_j^2/2,4d_j^2;L^2(D_j''))}  d_j^{-2}
+\|\nabla \phi_h\|_{L^2(d_j^2/2,4d_j^2;L^2(D_j''))}  d_j^{-1} )  \nn\\
&\leq  
(Ch^{1/2}d_j^{-1/2}+\epsilon) \|\partial_t\phi_h\|_{L^2(d_j^2/2,4d_j^2;L^2(D_j''))} \nn\\
&\quad 
+ C\epsilon^{-1}( \|\phi_h\|_{L^\infty(d_j^2/2,4d_j^2;L^2(D_j''))}  d_j^{-1}
+\|\nabla \phi_h\|_{L^2(d_j^2/2,4d_j^2;L^2(D_j''))}  d_j^{-1} )  .
\end{align}
With an obvious change of domains, \eqref{phihDj} and \eqref{phihtDj} imply 
\begin{align}
&\|\phi_h\|_{L^\infty(d_j^2/2,4d_j^2; L^2(D_j'))} 
+\|\nabla \phi_h\|_{L^2(d_j^2/2,4d_j^2; L^2(D_j'))} \\
&\leq C\|\partial_t\phi_h\|_{L^2(d_j^2/4,4d_j^2; L^2(D_j''))}h
+\epsilon^2\|\nabla \phi_h\|_{L^2(d_j^2/4,4d_j^2; L^2(D_j''))}  
+C\epsilon^{-2}\|\phi_h\|_{L^2(d_j^2/4,4d_j^2; L^2(D_j''))} d_j^{-1} . 
\nn
\end{align}
and
\begin{align}
\|\partial_t\phi_h\|_{L^2(d_j^2,4d_j^2;L^2(D_j))}   
&\leq  
(Ch^{1/2}d_j^{-1/2}+\epsilon) \|\partial_t\phi_h\|_{L^2(d_j^2/2,4d_j^2;L^2(D_j'))}\\
&\quad 
+ C\epsilon^{-1}( \|\phi_h\|_{L^\infty(d_j^2/2,4d_j^2;L^2(D_j'))}  d_j^{-1}
+\|\nabla \phi_h\|_{L^2(d_j^2/2,4d_j^2;L^2(D_j'))}  d_j^{-1} )  .
\nn
\end{align}
respectively. The last two inequalities further imply
\begin{align}\label{phihtxDj}
&\|\partial_t\phi_h\|_{L^2(d_j^2,4d_j^2;L^2(D_j))}  
+d_j^{-1}\|\nabla \phi_h\|_{L^2(d_j^2,4d_j^2;L^2(D_j))}  \nn\\ 
&\leq  
(Ch^{1/2}d_j^{-1/2}+C\epsilon^{-1}hd_j^{-1}+\epsilon) 
(\|\partial_t\phi_h\|_{L^2(d_j^2/4,4d_j^2;L^2(D_j''))}  
+d_j^{-1}\|\nabla \phi_h\|_{L^2(d_j^2/4,4d_j^2;L^2(D_j''))})\nn\\
&\quad 
+C\epsilon^{-3}\|\phi_h\|_{L^2(d_j^2/4,4d_j^2; L^2(D_j''))} d_j^{-2}  .
\end{align}

Finally, combining \eqref{phihtxOj} and \eqref{phihtxDj} yields 
\begin{align}
&\|\partial_t \phi_h\|_{L^2(Q_j)} +d_j^{-1}\|\nabla  \phi_h\|_{L^2(Q_j)} \nn\\
&\leq  (Ch^{1/2}d_j^{-1/2}+C\epsilon^{-1}hd_j^{-1}+\epsilon) 
\big( \|\partial_t \phi_h\|_{L^2(Q_{j}'' )}  
+ d_j^{-1}\|\nabla \phi_h\|_{L^2(Q_{j}'' )}\big)
+ C\epsilon^{-3}d_j^{-2} \|  \phi_h\|_{L^2(Q_{j}'')}  \nn\\
&\quad
+ C\epsilon^{-1}(\|\nabla  \phi_h(0)\|_{L^2(\Omega_j'')} 
+d_j^{-1}\| \phi_h(0)\|_{L^2(\Omega_j'')} )     .
\end{align}
Replacing $\Omega_j''$ by $\Omega_j'$ and replacing $Q_j''$ by $Q_j'$ in the last inequality, we obtain \eqref{locDEng0} and complete the proof of Lemma \ref{lemlocDEng0}.\qed 
\bigskip

\noindent{\it Proof of Lemma \ref{LocEEst}.}$\,\,\,$
Let $\widetilde\omega(t,x)$ be a smooth cut-off function 
which equals $1$ in $Q_j'$ 
and vanishes outside $Q_j''$, 
and let $\widetilde  \phi=\widetilde\omega  \phi$. 
Then $\widetilde\phi=\phi$ in $Q_j'$, which implies that 
\begin{align*}
&\big(\partial_t(\widetilde \phi-\phi_h),\chi_h \big)
+\big(\nabla (\widetilde \phi-\phi_h),\nabla\chi_h \big)
=0  ,\qquad\mbox{for}\,\,\,\,\chi_h \in S_h^0(\Omega_j ') ,\,\,\, t\in(0,d_j^2) ,\\
& \big(\partial_t(\widetilde \phi-\phi_h),\chi_h \big)
+\big(\nabla (\widetilde \phi-\phi_h),\nabla\chi_h \big)
=0  ,\qquad\mbox{for}\,\,\,\,\chi_h \in S_h^0(D_j ') ,\,\,\, t\in(d_j^2/4,4d_j^2) .
\end{align*}
Let $\widetilde \phi_h\in S_h$ be the solution of
\begin{align}\label{wdtzh}
\big(\partial_t(\widetilde\phi-\widetilde \phi_h),\chi_h \big)
+\big(\nabla (\widetilde\phi-\widetilde \phi_h),\nabla\chi_h \big)
=0   ,\qquad\forall\, \chi_h \in S_h  ,
\end{align}
with $\widetilde\phi_h(0)=\widetilde \phi(0)=0$ so that 
\begin{align}
&\big(\partial_t(\widetilde \phi_h-\phi_h),\chi_h \big)
+\big(\nabla (\widetilde \phi_h-\phi_h),\nabla\chi_h \big)
=0   ,\qquad\mbox{for}\,\,\,\,\chi_h \in S_h^0(\Omega_j'),\,\,\, t\in(0,d_j^2)  ,
\label{wdtzh2}\\
&\big(\partial_t(\widetilde \phi_h-\phi_h),\chi_h \big)
+\big(\nabla (\widetilde \phi_h- \phi_h),\nabla\chi_h \big)
=0   ,\qquad\mbox{for}\,\,\,\,\chi_h \in S_h^0(D_j') ,\,\,\, t\in(d_j^2/4,4d_j^2) .
\label{wdtzh3}
\end{align}
We shall estimate $\widetilde\phi-\widetilde\phi_h$ and $\widetilde\phi_h-\phi_h$ separately. 

The basic global energy estimates of \eqref{wdtzh} are 
(substituting $\chi_h=P_h\widetilde\phi-\widetilde\phi_h$ and $\chi_h=\partial_t(P_h\widetilde\phi-\widetilde\phi_h)$, respectively) 
\begin{align*} 
&\|\nabla (\widetilde\phi-\widetilde \phi_h) \|_{L^2({\cal Q})}^2
+\|\widetilde\phi-\widetilde \phi_h\|_{L^\infty(0,T;L^2(\Omega))}^2 
\leq C\|\partial_t (\widetilde \phi-\widetilde \phi_h)\|_{L^2({\cal Q})}
\| \widetilde \phi \|_{L^2({\cal Q})}
+C\| \nabla\widetilde\phi \|_{L^2({\cal Q})}^2,\\
&\|\partial_t (\widetilde\phi-\widetilde \phi_h) \|_{L^2({\cal Q})}^2
\leq C\|\partial_t  \widetilde \phi \|_{L^2({\cal Q})}^2
+C\|\nabla (\widetilde \phi-\widetilde \phi_h) \|_{L^2({\cal Q})}
\| \nabla\partial_t\widetilde \phi \|_{L^2({\cal Q})} ,
\end{align*}
which imply 
\begin{align*} 
&\|\partial_t(\widetilde \phi-\widetilde\phi_h)\|_{L^2({\cal Q})}^2
+d_j^{-2}\|\nabla \widetilde \phi-\widetilde\phi_h)\|_{L^2({\cal Q})}^2
+d_j^{-2}\|\widetilde \phi-\widetilde \phi_h\|_{L^\infty(0,T;L^2(\Omega))}^2 \\
&\leq Cd_j^{-2}\|\partial_t (\widetilde\phi-\widetilde\phi_h)\|_{L^2({\cal Q})}
\| \widetilde\phi \|_{L^2({\cal Q})}
+Cd_j^{-2}\| \nabla\widetilde\phi \|_{L^2({\cal Q})}^2,\\
&\quad
+C\|\partial_t  \widetilde\phi \|_{L^2({\cal Q})}^2
+C\|\nabla (\widetilde \phi-\widetilde \phi_h)\|_{L^2({\cal Q})}
\| \nabla\partial_t\widetilde \phi \|_{L^2({\cal Q})} \\
&\leq \frac{1}{2}\|\partial_t (\widetilde\phi-\widetilde\phi_h)\|_{L^2({\cal Q})}^2
+Cd_j^{-4}\| \widetilde \phi \|_{L^2({\cal Q})}^2
+Cd_j^{-2}\| \nabla\widetilde \phi \|_{L^2({\cal Q})}^2,\\
&\quad
+C\|\partial_t  \widetilde \phi \|_{L^2({\cal Q})}^2
+\frac{1}{2}d_j^{-2}\|\nabla (\widetilde \phi-\widetilde \phi_h)\|_{L^2({\cal Q})}^2
+Cd_j^{2}\| \nabla\partial_t\widetilde\phi \|_{L^2({\cal Q})}^2 .
\end{align*}
The first and fifth terms on the right-hand side above can be absorbed by the left-hand side, and 
the last inequality further reduces to  
\begin{align} \label{global-tilde-phi}
&\|\partial_t(\widetilde\phi-\widetilde\phi_h)\|_{L^2({\cal Q})} 
+d_j^{-1}\|\nabla (\widetilde\phi-\widetilde\phi_h)\|_{L^2({\cal Q})}
+d_j^{-1}\|\widetilde \phi-\widetilde \phi_h\|_{L^\infty(0,T;L^2(\Omega))} \nn\\
&\leq Cd_j^{-2}\| \widetilde\phi \|_{L^2({\cal Q})} 
+Cd_j^{-1}\| \nabla\widetilde\phi \|_{L^2({\cal Q})}  
+C\|\partial_t  \widetilde\phi \|_{L^2({\cal Q})} 
+Cd_j \| \nabla\partial_t\widetilde \phi \|_{L^2({\cal Q})} \nn\\
&\le 
Cd_j^{-2}\| \phi \|_{L^2(Q_j'')} 
+Cd_j^{-1}\| \nabla\phi \|_{L^2(Q_j'')}  
+C\|\partial_t\phi \|_{L^2(Q_j'')} 
+Cd_j \| \nabla\partial_t \phi \|_{L^2(Q_j'')}. 
\end{align}
By applying Lemma \ref{lemlocDEng0} to \eqref{wdtzh2}-\eqref{wdtzh3}, 
we obtain 
\begin{align}
&\|\partial_t(\widetilde\phi_h-\phi_h)\|_{L^2(Q_j)} 
+ d_j^{-1}\|\nabla(\widetilde \phi_h-\phi_h)\|_{L^2(Q_j)} \nn\\
&\leq  C\epsilon^{-1}\left(\|\nabla (\widetilde \phi_h-\phi_h)(0)\|_{L^2(\Omega_j')} 
+d_j^{-1}\|(\widetilde \phi_h-\phi_h)(0)\|_{L^2(\Omega_j')} \right) 
+ C\epsilon^{-3}d_j^{-2} \| \widetilde \phi_h-\phi_h\|_{L^2(Q_{j}')}\nn\\
&\quad
+ (Ch^{1/2}d_j^{-1/2}+C\epsilon^{-1}hd_j^{-1}+\epsilon)\big(
\|\partial_t(\widetilde \phi_h-\phi_h)\|_{L^2(Q_{j}')}  
+ d_j^{-1}\|\nabla(\widetilde \phi_h-\phi_h)\|_{L^2(Q_{j}')}\big) \nn\\
&= C\epsilon^{-1} \left(\|\nabla \phi_h(0)\|_{L^2(\Omega_j')} 
+d_j^{-1}\|\phi_h(0)\|_{L^2(\Omega_j')} \right) 
+ C\epsilon^{-3}d_j^{-2} \| \widetilde \phi_h-\phi_h\|_{L^2(Q_{j}')}\nn\\
&\quad
+ (Ch^{1/2}d_j^{-1/2} +C\epsilon^{-1}hd_j^{-1} +\epsilon)\big(
\|\partial_t(\widetilde \phi_h-\phi_h)\|_{L^2(Q_{j}')}  
+ d_j^{-1}\|\nabla(\widetilde \phi_h-\phi_h)\|_{L^2(Q_{j}')}\big)
\end{align}
where we have used the identity $\widetilde \phi_h(0)=0$ in the last step. 
Splitting $\widetilde \phi_h-\phi_h$ into $(\widetilde \phi-\phi_h)+(\widetilde \phi_h- \widetilde \phi)$ in the right-hand side of the last inequality yields 
\begin{align}
&\|\partial_t(\widetilde\phi_h-\phi_h)\|_{L^2(Q_j)} 
+ d_j^{-1}\|\nabla(\widetilde \phi_h-\phi_h)\|_{L^2(Q_j)} \nn\\
&\le C\epsilon^{-1} \left(\|\nabla \phi_h(0)\|_{L^2(\Omega_j')} 
+d_j^{-1}\|\phi_h(0)\|_{L^2(\Omega_j')} \right)  \nn\\
&\quad
+ C\epsilon^{-3}d_j^{-2} \| \widetilde \phi- \phi_h\|_{L^2(Q_{j}')} \nn\\
&\quad 
+ (Ch^{1/2}d_j^{-1/2} +C\epsilon^{-1}hd_j^{-1} +\epsilon)\big(
\|\partial_t(\widetilde \phi -\phi_h)\|_{L^2(Q_{j}')}  
+ d_j^{-1}\|\nabla(\widetilde \phi -\phi_h)\|_{L^2(Q_{j}')}\big) \nn\\
&\quad
+ C\epsilon^{-3}d_j^{-2} \| \widetilde \phi_h- \widetilde \phi\|_{L^2(Q_{j}')}\nn\\
&\quad
+ (Ch^{1/2}d_j^{-1/2} +C\epsilon^{-1}hd_j^{-1} +\epsilon)\big(
\|\partial_t(\widetilde \phi_h-\widetilde\phi)\|_{L^2(Q_{j}')}  
+ d_j^{-1}\|\nabla(\widetilde \phi_h-\widetilde\phi )\|_{L^2(Q_{j}')}\big) \nn\\ 
&\le C\epsilon^{-1} \left(\|\nabla \phi_h(0)\|_{L^2(\Omega_j')} 
+d_j^{-1}\|\phi_h(0)\|_{L^2(\Omega_j')} \right) 
+ C\epsilon^{-3}d_j^{-2} \| \phi -\phi_h\|_{L^2(Q_{j}')}\nn\\
&\quad
+ (Ch^{1/2}d_j^{-1/2} +C\epsilon^{-1}hd_j^{-1} +\epsilon)\big(
\|\partial_t(\phi -\phi_h)\|_{L^2(Q_{j}')}  
+ d_j^{-1}\|\nabla(\phi -\phi_h)\|_{L^2(Q_{j}')}\big) \nn \\
&\quad
+  \bigg(Ch^{1/2}d_j^{-1/2} +C\epsilon^{-1}hd_j^{-1} +C\epsilon +C\epsilon^{-3}\bigg)  d_j^{-2}\| \phi \|_{L^2(Q_j'')} \nn \\
&\quad
+  \bigg(Ch^{1/2}d_j^{-1/2} +C\epsilon^{-1}hd_j^{-1} +C\epsilon +C\epsilon^{-3}\bigg) \bigg(d_j^{-1}\| \nabla\phi \|_{L^2(Q_j'')}  
+\|\partial_t\phi \|_{L^2(Q_j'')} 
+d_j \| \nabla\partial_t \phi \|_{L^2(Q_j'')}\bigg) \nn \\
\end{align}
where we have used the identity $\widetilde\phi=\phi$ on $Q_j'$ and \eqref{global-tilde-phi} in the last step. 
Since 
$$\bigg(Ch^{1/2}d_j^{-1/2} +C\epsilon^{-1}hd_j^{-1} +C\epsilon +C\epsilon^{-3}\bigg)
\le C\epsilon^{-3} ,$$ 
the last inequality reduces to 
\begin{align}\label{local-tilde-phi}
&\|\partial_t(\widetilde\phi_h-\phi_h)\|_{L^2(Q_j)} 
+ d_j^{-1}\|\nabla(\widetilde \phi_h-\phi_h)\|_{L^2(Q_j)} \nn\\
&\le C\epsilon^{-1} \left(\|\nabla \phi_h(0)\|_{L^2(\Omega_j')} 
+d_j^{-1}\|\phi_h(0)\|_{L^2(\Omega_j')} \right) 
+ C\epsilon^{-3}d_j^{-2} \| \phi -\phi_h\|_{L^2(Q_{j}')}\nn\\
&\quad
+ (Ch^{1/2}d_j^{-1/2} +C\epsilon^{-1}hd_j^{-1} +\epsilon)\big(
\|\partial_t(\phi -\phi_h)\|_{L^2(Q_{j}')}  
+ d_j^{-1}\|\nabla(\phi -\phi_h)\|_{L^2(Q_{j}')}\big)  \\
&\quad
+ C\epsilon^{-3}\big( d_j^{-2}\| \phi \|_{L^2(Q_j'')} +d_j^{-1}\| \nabla\phi \|_{L^2(Q_j'')}  
+\|\partial_t\phi \|_{L^2(Q_j'')} 
+d_j \| \nabla\partial_t \phi \|_{L^2(Q_j'')}\big) \nn
,
\end{align}

The estimates \eqref{global-tilde-phi} and \eqref{local-tilde-phi} imply  
\begin{align*}
&\|\partial_t(\phi-\phi_h)\|_{L^2(Q_j)} 
+ d_j^{-1}\|\nabla(\phi-\phi_h)\|_{L^2(Q_j)} \nn\\
&=\|\partial_t(\widetilde\phi-\phi_h)\|_{L^2(Q_j)} 
+ d_j^{-1}\|\nabla(\widetilde \phi-\phi_h)\|_{L^2(Q_j)} \nn\\
&\le \|\partial_t(\widetilde\phi_h-\phi_h)\|_{L^2(Q_j)} 
+ d_j^{-1}\|\nabla(\widetilde \phi_h-\phi_h)\|_{L^2(Q_j)}  
\qquad\mbox{(here we use triangle inequality)}\nn\\
&\quad 
+\|\partial_t(\widetilde\phi-\widetilde\phi_h)\|_{L^2(Q_j)} 
+ d_j^{-1}\|\nabla(\widetilde \phi-\widetilde\phi_h)\|_{L^2(Q_j)}\nn\\
&\leq  C\epsilon^{-1}\left(\|\nabla \phi_h(0)\|_{L^2(\Omega_j'')} 
+d_j^{-1}\|\phi_h(0)\|_{L^2(\Omega_j'')} \right) 
+ C\epsilon^{-3}d_j^{-2} \| \phi-\phi_h\|_{L^2(Q_{j}'')}\nn\\
&\quad
+ (Ch^{1/2}d_j^{-1/2}+C\epsilon^{-1}hd_j^{-1}+\epsilon)\big(
\|\partial_t(\phi-\phi_h)\|_{L^2(Q_{j}'')}  
+ d_j^{-1}\|\nabla(\phi-\phi_h)\|_{L^2(Q_{j}'')}\big) \nn\\
&\quad
+ C\epsilon^{-3}\big( d_j^{-2}\| \phi \|_{L^2(Q_j'')} +d_j^{-1}\| \nabla\phi \|_{L^2(Q_j'')}  
+\|\partial_t\phi \|_{L^2(Q_j'')} 
+d_j \| \nabla\partial_t \phi \|_{L^2(Q_j'')}\big)
\end{align*}
Replacing $\phi$ by $\phi-I_h\phi$, replacing $\Omega_j''$ by $\Omega_j'$, and replacing $Q_j''$ by $Q_j'$ in the last inequality, we obtain \eqref{LocEngErr} and complete the proof of Lemma \ref{LocEEst}.
\qed

\section*{Appendix B: Property (P3) and the operator $I_h$}
\label{AppB}

\renewcommand{\thelemma}{B.\arabic{lemma}}
\renewcommand{\theequation}{B.\arabic{equation}}
\setcounter{equation}{0}\setcounter{lemma}{0}

Let $\Phi_i$ be the basis function of the finite element space $S_h$ corresponding to the finite element nodes $x_i\in\Omega$, $i=1\dots,M$. In other words, we have $\Phi_j(x_i)=\delta_{ij}$ (the Kronecker symbol).  
Let $\tau_i$ denote the union of triangles (or tetrahedra in $\R^3$) whose closure contain the node $x_i$. 
For any function $v\in L^2(\Omega)$, we denote by 
$P^{(i)}_hv$ the local $L^2$ projection onto $S_h(\tau_i)$ (the space of finite element functions defined on the region $\tau_i$). The operator $I_h:L^2(\Omega)\rightarrow S_h$ is defined as (in the spirit of Cl\'ement's interpolation operator, cf. \cite{Clement1975}) 
\begin{align}\label{Def-Ih}
(I_hv)(x)=\sum_{i=1}^M(P^{(i)}_hv)(x_i)\Phi_i(x),\quad \mbox{for}\,\,\, x\in\Omega ,
\end{align}
which equals zero on the boundary $\partial\Omega$ (as every $\Phi_i$ equals zero on $\partial\Omega$). 

Now we prove that the operator $I_h$ defined in \eqref{Def-Ih} satisfies property (P3) of Section \ref{Sec:notation}. To this end, we let $S_h'$ be the finite element space subject to the same mesh as $S_h$, with the same order of finite elements, but not necessarily zero on the boundary $\partial\Omega$. We denote by $x_j'$, $j=1,\dots,m$, the finite element nodes on the boundary $\partial\Omega$, and we denote by $\Phi_j'$ the basis function corresponding to the node $x_j'$. The notation $\tau_j'$ will denote the union of triangles (or tetrahedra in $\R^3$) whose closure contain $x_j'$.  With these notations, the space $S_h'$ is spanned by the basis functions $\Phi_i$, $i=1,\dots,M$, and $\Phi_j'$, $j=1,\dots,m$. 
We define an auxiliary operator $\widetilde I_h:H^1(\Omega)\rightarrow S_h'$ by setting 
\begin{align}\label{Def-tilde-Ih} 
(\widetilde I_h\phi)(x)=(I_h\phi)(x)
+\sum_{j=1}^m(\widetilde P^{(j)}_h\phi)(x_j')\Phi_j'(x),\quad \mbox{for}\,\,\, x\in\Omega ,
\end{align}
where $\widetilde P^{(j)}_h\phi$ is the $L^2$ projection of $\phi|_{\partial\Omega}$ (trace of $\phi$ on the boundary) onto $S_h(\partial\Omega\cap\overline\tau_j')$ (the space of finite element functions on $\partial\Omega\cap\overline\tau_i'$, a piece of the boundary).
The definition \eqref{Def-tilde-Ih} implies  
\begin{align}
I_h\phi=\widetilde I_h\phi,\qquad\forall\, \phi\in H^1_0(\Omega) .
\end{align}
Hence, in order to prove property (P3)-(1), we only need to prove the corresponding error estimate for the operator $\widetilde I_h$. 

In fact, the definition \eqref{Def-tilde-Ih} guarantees the following local stability: 
\begin{align*}
&\|\widetilde I_h\phi\|_{L^2(\tau_j')}+h^{\frac12}\|\widetilde I_h\phi\|_{L^2(\partial\Omega\cap \tau_j')} 
+ h \|\nabla\widetilde I_h\phi\|_{L^2(\tau_j')}  
\le C(\|\phi\|_{L^2(\widetilde\tau_j')}+h^{\frac12}\|\phi\|_{L^2(\partial\Omega\cap\widetilde\tau_j')}),\\
&\|\widetilde I_h\phi\|_{L^2(\tau_i)} + h \|\nabla\widetilde I_h\phi\|_{L^2(\tau_j)}  
\le C\|\phi\|_{L^2(\widetilde\tau_i)},
\end{align*}
where $\widetilde\tau_j‘$ is the union of triangles (or tetrahedra in $\R^3$) whose closure intersect the closure of $\tau_j'$ (boundary triangle/tetrahedron), and 
$\widetilde\tau_i$ is the union of triangles (or tetrahedra in $\R^3$) whose closure intersect the closure of $\tau_i$ (interior triangle/tetrahedron). 
Let $P_h'$ denote the $L^2$ projection from $L^2(\Omega)$ onto the finite element space $S_h'$ . Then substituting $\phi=v-P_h'v$ into the two inequalities above yields 
\begin{align*}
&\|v-\widetilde I_hv\|_{L^2(\tau_j')}+h^{\frac12}\|v-\widetilde I_hv\|_{L^2(\partial\Omega\cap \tau_j')} 
 +h \|\nabla(v-\widetilde I_hv)\|_{L^2(\tau_j')}  \,\nonumber \\ 
&\le C(\|v-P_h'v\|_{L^2(\widetilde\tau_j')}+h^{\frac12}\|v-P_h'v\|_{L^2(\partial\Omega\cap\widetilde\tau_j')}) 
\end{align*}
and 
\begin{align*}
&\|v-\widetilde I_hv\|_{L^2(\tau_i)} + h \|\nabla(v-\widetilde I_hv)\|_{L^2(\tau_j)}  
\le C\|v-P_h'v\|_{L^2(\widetilde\tau_i)}. \quad 
\end{align*}
Summing up the two inequalities above for $i=1,\dots,M$ and $j=1,\dots,m$ yields 
\begin{align*}
&\|v-\widetilde I_hv\|_{L^2(\Omega)}+h^{\frac12}\|v-\widetilde I_hv\|_{L^2(\partial\Omega)} 
+h \|\nabla(v-\widetilde I_hv)\|_{L^2(\Omega)}  \\
&\le C(\|v-P_h'v\|_{L^2(\Omega)}+h^{\frac12}\|v-P_h'v\|_{L^2(\partial\Omega)}) \\ 
&\le C(\|v-P_h'v\|_{L^2(\Omega)}+h^{\frac12}\|v-P_h'v\|_{L^2(\Omega)}^{\frac12}\|v-P_h'v\|_{H^1(\Omega)}^{\frac12})
&&\mbox{(interpolation inequality)}\\
&\le C(\|v-P_h'v\|_{L^2(\Omega)}+h\|v-P_h'v\|_{H^1(\Omega)}). 
&&\mbox{(H\"older's inequality)}\\
&\le 
\left\{\begin{array}{ll}
Ch\|v\|_{H^1(\Omega)}
&\mbox{if} \,\,v\in  H^1(\Omega),\\[3pt]
Ch^2\|v\|_{H^2(\Omega)}
&\mbox{if} \,\,v\in  H^2(\Omega),
\end{array}
\right.
\end{align*}
where the last inequality is the basic estimate of the $L^2$ projection $P_h':L^2(\Omega)\rightarrow S_h'$ (without imposing boundary condition). 
By the complex interpolation method, we have
\begin{align*}
&\|v-\widetilde I_hv\|_{L^2(\Omega)}+h^{\frac12}\|v-\widetilde I_hv\|_{L^2(\partial\Omega)} +h \|\nabla(v-\widetilde I_hv)\|_{L^2(\Omega)} 
\le 
Ch^{1+\alpha}\|v\|_{H^{1+\alpha}(\Omega)} .
\end{align*}
Hence, for $v\in H^{1+\alpha}(\Omega)\cap H^1_0(\Omega)$, we have $I_hv=\widetilde I_hv=0$ on $\partial\Omega$ and 
\begin{align}\label{v-Ihv-intp}
&\|v-I_hv\|_{L^2(\Omega)} +h \|\nabla(v- I_hv)\|_{L^2(\Omega)}  \nonumber \\
&=\|v-\widetilde I_hv\|_{L^2(\Omega)}+h^{\frac12}\|v-\widetilde I_hv\|_{L^2(\partial\Omega)} +h \|\nabla(v- \widetilde I_hv)\|_{L^2(\Omega)}    \nonumber \\
&\le 
Ch^{1+\alpha}\|v\|_{H^{1+\alpha}(\Omega)} .
\end{align}
This proves property (P3)-(1) in Section \ref{Sec:notation}. 
The other properties in (P3) are simple consequences of the definition of the operator $I_h$.   
\qed

\section*{Appendix C: Proof of \eqref{Linfty-Deltah-1}}\label{Sec:AppendC}

\renewcommand{\thelemma}{C.\arabic{lemma}}
\renewcommand{\theequation}{C.\arabic{equation}}
\setcounter{equation}{0}\setcounter{lemma}{0}

The proof of \eqref{Linfty-Deltah-1} requires some properties of the finite element space described in Section \ref{Sec2-2}. 

It suffices to prove that the solution $w_h\in S_h$ of the finite element equation 
\begin{align}\label{Eq-wh-fh}
\Delta_hw_h=f_h 
\end{align}
satisfies 
\begin{align}\label{Est-wh-fh}
\|w_h\|_{L^\infty} \le C \|f_h\|_{L^\infty} .
\end{align}
To this end, we define $w\in H^1_0(\Omega)$ as the solution of the following PDE problem:
\begin{align}
\left\{
\begin{array}{ll}
\Delta w =f_h &\mbox{in}\,\,\,\Omega, \\
w=0               &\mbox{on}\,\,\,\partial\Omega . 
\end{array}
\right.
\end{align}
Then $w_h$ is the Ritz projection of $w$, and the following standard $H^1$-norm error estimate holds for some $\alpha\in \left(\frac{1}{2},1\right)$: 
\begin{align}\label{CInhwwh}
\|I_hw-w_h\|_{H^1}
\le C\|I_hw-w\|_{H^{1}}
\le Ch^{\alpha}\|w\|_{H^{1+\alpha}}
\le Ch^{\alpha}\|f_h\|_{H^{-1+\alpha}} ,
\end{align}
where the first inequality above is due to $H^1$-stability of the Ritz projection, the second inequality due to \eqref{v-Ihv-intp}, and the last inequality due to Lemma \ref{RegPoiss}. 
Consequently, we have  
\begin{align}
\|I_hw-w_h\|_{L^\infty}
&\le Ch^{-\frac{1}{2}}\|I_hw-w_h\|_{H^1} &&\mbox{(inverse inequality)} \nn\\
&\le Ch^{\alpha-\frac{1}{2}}\|f_h\|_{H^{-1+\alpha}} &&\mbox{(use \eqref{CInhwwh})}\nn\\
&\le Ch^{\alpha-\frac{1}{2}}\|f_h\|_{L^2} &&\mbox{($L^2\hookrightarrow H^{-1+\alpha}$)} \nn\\ 
&\le Ch^{\alpha-\frac{1}{2}}\|f_h\|_{L^\infty} , \label{CInhwwh1} \\ 
\|w-I_hw\|_{L^\infty}
&\le \|w-P_h w\|_{L^\infty}+\|P_h w-I_hw\|_{L^\infty} &&\mbox{(triangle inequality)} \nn\\ 
&\le Ch^{1+\alpha-\frac{N}{2}} \|w\|_{C^{1+\alpha-\frac{N}{2}}}
       +Ch^{-\frac{N}{2}}\|P_h w-I_hw\|_{L^2} &&\mbox{(inverse inequality)} \nn\\
&\le Ch^{1+\alpha-\frac{N}{2}} \|w\|_{C^{1+\alpha-\frac{N}{2}}}
       +Ch^{-\frac{N}{2}}\|w-I_hw\|_{L^2} &&\mbox{($L^2$-stability of $P_h$)} \nn\\
&\le Ch^{1+\alpha-\frac{N}{2}} \|w\|_{H^{1+\alpha}}
       +Ch^{1+\alpha-\frac{N}{2}}\|w\|_{H^{1+\alpha}} &&\mbox{($H^{1+\alpha}\hookrightarrow C^{1+\alpha-\frac{N}{2}}$ and \eqref{v-Ihv-intp})}  \nn\\
&\le Ch^{\alpha-\frac{1}{2}} \|f_h\|_{H^{-1+\alpha}} &&\mbox{(use Lemma \ref{RegPoiss} and $N=2,3$)} \nn\\
&\le Ch^{\alpha-\frac{1}{2}}\|f_h\|_{L^2} &&\mbox{($L^2\hookrightarrow H^{-1+\alpha}$ for $\alpha\in\left(\frac12,1\right)$)} \nn\\ 
&\le Ch^{\alpha-\frac{1}{2}}\|f_h\|_{L^\infty} . \label{CInhwwh2}
\end{align}
Consequently, the triangle inequality implies
\begin{align}
\|w_h\|_{L^\infty}
&\le \|w\|_{L^\infty}+\|w-I_hw\|_{L^\infty}+ \|I_hw-w_h\|_{L^\infty} \nn\\
&\le \|w\|_{L^\infty}+ Ch^{\alpha-\frac{1}{2}}\|f_h\|_{L^\infty}  + Ch^{\alpha-\frac{1}{2}}\|f_h\|_{L^\infty}  
&&\mbox{(use \eqref{CInhwwh1} and \eqref{CInhwwh2})} \nn\\
&\le C\|w\|_{H^{1+\alpha}}+ Ch^{\alpha-\frac{1}{2}}\|f_h\|_{L^\infty}  + Ch^{\alpha-\frac{1}{2}}\|f_h\|_{L^\infty}  
&&\mbox{($H^{1+\alpha}\hookrightarrow L^\infty$ for $\alpha\in\left(\frac12,1\right)$, $N=2,3$)}  \nn\\
&\le C \|f_h\|_{H^{-1+\alpha}}  + Ch^{\alpha-\frac{1}{2}}\|f_h\|_{L^\infty} 
&&\mbox{(use Lemma \ref{RegPoiss})} \nn\\
&\le C\|f_h\|_{L^2} + Ch^{\alpha-\frac{1}{2}}\|f_h\|_{L^\infty}  
&&\mbox{($L^2\hookrightarrow H^{-1+\alpha}$ for $\alpha\in\left(\frac12,1\right)$)} \nn\\ 
&\le C\|f_h\|_{L^\infty} .
\end{align}
This proves \eqref{Est-wh-fh}. The proof of \eqref{Linfty-Deltah-1} is complete. \qed

{\small
\bibliographystyle{abbrv}
\bibliography{max-reg}

\begin{thebibliography}{10}

\bibitem{AdamsFournier2003}
R.~A. Adams and J.~J.~F. Fournier.
\newblock {\em {Sobolev Spaces}}.
\newblock Academic Press, Amsterdam, second edition, 2003.

\bibitem{AkrivisLi2017}
G.~Akrivis and B.~Li.
\newblock Maximum norm analysis of implicit--explicit backward difference
  formulas for nonlinear parabolic equations.
\newblock {\em IMA J. Numer. Anal.}, 2017, DOI: 10.1093/imanum/drx008.

\bibitem{AkrivisLiLubich2016}
G.~Akrivis, B.~Li, and C.~Lubich.
\newblock {Combining maximal regularity and energy estimates for time
  discretizations of quasilinear parabolic equations}.
\newblock {\em Math. Comp.}, (86):1527--1552, 2017.

\bibitem{Amann1995}
H.~Amann.
\newblock {\em {Linear and Quasilinear Parabolic Problems. Volume I: Abstract
  linear theory.}}
\newblock Birkh\"auser Boston Inc., Boston, MA, 1995.

\bibitem{BerghLofstrom1976}
J.~Bergh and J.~L\"ofstr\"om.
\newblock {\em {Interpolation Spaces: An Introduction}}.
\newblock Springer-Verlag Berlin Heidelberg, Printed in Germany, 1976.

\bibitem{Clement1975}
P.~Cl\'ement.
\newblock {Approximation by finite element functions using local
  regularization}.
\newblock {\em Revue française d'automatique, informatique, recherche
  op\'erationnelle. Analyse numérique}, 9:77--84, 1975.

\bibitem{ClementPruss1992}
P.~Cl\'ement and J.~Pr\"uss.
\newblock {Global existence for a semilinear prabolic Volterra equation}.
\newblock {\em Mathematische Zeitschrift}, 209:17--26, 1992.

\bibitem{Crouzeix2003}
M.~Crouzeix.
\newblock {Contractivity and analyticity in $l^p$ of some approximation of the
  heat equation}.
\newblock {\em Numerical Algorithms}, 33:193--201, 2003.

\bibitem{Dauge1988}
M.~Dauge.
\newblock {\em {Elliptic Boundary Value Problems in Corner Domains}}.
\newblock Springer-Verlag Berlin Heidelberg, 1988.

\bibitem{Davis1989}
E.~B. Davies.
\newblock {\em {Heat Kernels and Spectral Theory}}.
\newblock Cambridge University Press, Cambridge, 1989.

\bibitem{DLM2009}
A.~Demlow, O.~Lakkis, and C.~Makridakis.
\newblock {A posteriori error estimates in the maximum norm for parabolic
  problem}.
\newblock {\em SIAM J. Numer. Anal.}, 47:2157--2176,, 2009.

\bibitem{DLSW2012}
A.~Demlow, D.~Leykekhman, A.~H. Schatz, and L.~B. Wahlbin.
\newblock {Best approximation property in the norm for finite element methods
  on graded meshes}.
\newblock {\em Math. Comp.}, 81:743--764, 2012.

\bibitem{DiPietroDroniou2017}
D.~A. Di~Pietro and J.~Droniou.
\newblock {$W^{s,p}$-approximation properties of elliptic projectors on
  polynomial spaces, with application to the error analysis of a Hybrid
  High-Order discretisation of Leray--Lions problems}.
\newblock {\em Math. Models Methods Appl. Sci.}, 27:879--908, 2017.

\bibitem{Geissert2006}
M.~Geissert.
\newblock {Discrete maximal $L^p$ regularity for finite element operators}.
\newblock {\em SIAM J. Numer. Anal.}, 44:677--698, 2006.

\bibitem{Geissert2007}
M.~Geissert.
\newblock {Applications of discrete maximal $L^p$ regularity for finite element
  operators}.
\newblock {\em Numer. Math.}, 108:121--149, 2007.

\bibitem{Grafakos2008}
L.~Grafakos.
\newblock {\em {Classical Fourier Analysis}}.
\newblock Springer Science + Business Media, LLC, second edition, 2008.

\bibitem{GLRS2009}
J.~Guzm\'an, D.~Leykekhman, J.~Rossmann, and A.~H. Schatz.
\newblock {H\"older estimates for Green’s functions on convex polyhedral
  domains and their applications to finite element methods}.
\newblock {\em Numer. Math.}, 112:221--243,, 2009.

\bibitem{Hansbo2002}
A.~Hansbo.
\newblock {Strong stability and non-smooth data error estimates for
  discretizations of linear parabolic problems}.
\newblock {\em BIT Numer. Math.}, 42:351--379, 2002.

\bibitem{DouglasDupontWahlbin1974}
T.~D. J.~Douglas, Jr. and L.~Wahlbin.
\newblock {The stability in $L^q$ of the $L^2$-projection into finite element
  function spaces}.
\newblock {\em Numer. Math.}, 23:193--197, 1975.

\bibitem{KaltonoLancien2000}
N.~J. Kalton and G.~Lancien.
\newblock {A solution to the problem of $L^p$-maximal regularity}.
\newblock {\em Math. Z.}, 235:559--568, 2000.

\bibitem{KemmochiSaito}
T.~Kemmochi and N.~Saito.
\newblock {Discrete maximal regularity and the finite element method for
  parabolic equations}.
\newblock http://arXiv.org/abs/1602.06864.

\bibitem{KovacsLiLubich2016}
B.~Kov\'acs, B.~Li, and C.~Lubich.
\newblock A-stable time discretizations preserve maximal parabolic regularity.
\newblock {\em SIAM J. Numer. Anal.}, 54:3600--3624, 2016.

\bibitem{Krantz1999}
S.~Krantz.
\newblock {\em {A Panorama of Harmonic Analysis}}.
\newblock American History Series. Mathematical Association of America, 1999.

\bibitem{KunstmanLiLubich2016}
P.~C. Kunstmann, B.~Li, and C.~Lubich.
\newblock {Runge-Kutta time discretization of nonlinear parabolic equations
  studied via discrete maximal parabolic regularity}.
\newblock 2016, arXiv:1606.03692.

\bibitem{KunstmannWeis2004}
P.~C. Kunstmann and L.~Weis.
\newblock {\em {Maximal $L^p$-regularity for parabolic equations, Fourier
  multiplier theorems and $H^\infty$-functional calculus}}, pages 65--311.
\newblock Functional Analytic Methods for Evolution Equations, edited by M.
  Lannelli, R. Nagel, and S. Piazzera. Springer Berlin Heidelberg, Berlin,
  Heidelberg, 2004.

\bibitem{LSU1988}
O.~Ladyzhenskai͡a, V.~Solonnikov, and N.~Ural'tseva.
\newblock {\em {Linear and Quasi-linear Equations of Parabolic Type}}.
\newblock American Mathematical Society, Translations of Mathematical
  Monographs.

\bibitem{Leykekhman2004}
D.~Leykekhman.
\newblock {Pointwise localized error estimates for parabolic finite element
  equations}.
\newblock {\em Numer. Math.}, 96:583--600, 2004.

\bibitem{LeykekhmanVexler2016}
D.~Leykekhman and B.~Vexler.
\newblock {Finite element pointwise results on convex polyhedral domains}.
\newblock {\em SIAM J. Numer. Anal.}, 54:561--587, 2016.

\bibitem{LeykekhmanVexler2016-2}
D.~Leykekhman and B.~Vexler.
\newblock {Pointwise best approximation results for Galerkin finite element
  solutions of parabolic problems}.
\newblock {\em SIAM J. Numer. Anal.}, 54:1365--1384, 2016.

\bibitem{LeykekhmanVexler2016-NM}
D.~Leykekhman and B.~Vexler.
\newblock {Discrete maximal parabolic regularity for Galerkin finite element
  methods}.
\newblock {\em Numer. Math.}, 135:923--952, 2017.

\bibitem{Li2015}
B.~Li.
\newblock Maximum-norm stability and maximal ${L}^p$ regularity of {FEM}s for
  parabolic equations with {L}ipschitz continuous coefficients.
\newblock {\em Numer. Math.}, 131:489--516, 2015.

\bibitem{LiSun2015-regularity}
B.~Li and W.~Sun.
\newblock Regularity of the diffusion-dispersion tensor and error analysis of
  {FEM}s for a porous media flow.
\newblock {\em SIAM J. Numer. Anal.}, 53:1418--1437, 2015.

\bibitem{LiSun2016}
B.~Li and W.~Sun.
\newblock Maximal ${L}^p$ analysis of finite element solutions for parabolic
  equations with nonsmooth coefficients in convex polyhedra.
\newblock {\em Math. Comp.}, 86:1071--1102, 2017.

\bibitem{LiSun2017}
B.~Li and W.~Sun.
\newblock {Maximal regularity of fully discrete finite element solutions of
  parabolic equations}.
\newblock {\em SIAM J. Numer. Anal.}, 55:521--542, 2017.

\bibitem{Lions1996}
P.~L. Lions.
\newblock {\em {Mathematical Topics in Fluid Mechanics. Volume 1.
  Incompressible Models.}}
\newblock Clarendon Press, Oxford, 1996.

\bibitem{NitscheWheeler1982}
J.~A. Nitsche and M.~F. Wheeler.
\newblock {$L_\infty$-boundedness of the finite element Galerkin operator for
  parabolic problems}.
\newblock {\em Numer. Funct. Anal. Optim.}, 4:325--353, 1982.

\bibitem{OdenReddy2012}
J.~T. Oden and J.~N. Reddy.
\newblock {\em An Introduction to the Mathematical Theory of Finite Elements}.
\newblock Dover Books on Engineering. Dover Publications, 2012.

\bibitem{Ouhabaz1995}
E.-M. Ouhabaz.
\newblock {Gaussian estimates and holomorphy of semigroups}.
\newblock {\em Proceedings of the American Mathematical Society},
  123:1465--1474, 1995.

\bibitem{Palencia1996}
C.~Palencia.
\newblock {Maximum norm analysis of completely discrete finite element methods
  for parabolic problems}.
\newblock {\em SIAM J. Numer. Anal.}, 33:1654--1668,, 1996.

\bibitem{Rannacher1991}
R.~Rannacher.
\newblock {\em {$L^\infty$-stability estimates and asymptotic error expansion
  for parabolic finite element equations}}.
\newblock Extrapolation and Defect Correction, Bonner Mathematische Schriften
  228. University of Bonn, pp. 74--94, 1991.

\bibitem{RannacherScott1982}
R.~Rannacher and R.~Scott.
\newblock {Some optimal error estimates for piecewise linear finite element
  approximations}.
\newblock {\em Math. Comp.}, 38:437--445,, 1982.

\bibitem{Schatz1980}
A.~H. Schatz.
\newblock {A weak discrete maximum principle and stability of the finite
  element method in $L_\infty$ on plane polygonal domains. I.}
\newblock {\em Math. Comp.}, 31:77--91,, 1980.

\bibitem{SchatzThomeeWahlbin1980}
A.~H. Schatz, V.~Thom\'{e}e, and L.~B. Wahlbin.
\newblock {Maximum norm stability and error estimates in parabolic finite
  element equations}.
\newblock {\em Comm. Pure Appl. Math.}, 33:265--304,, 1980.

\bibitem{SchatzThomeeWahlbin1998}
A.~H. Schatz, V.~Thom\'{e}e, and L.~B. Wahlbin.
\newblock {Stability, analyticity, and almost best approximation in maximum
  norm for parabolic finite element equations}.
\newblock {\em Comm. Pure Appl. Math.}, 51:1349--1385, 1998.

\bibitem{SchatzWahlbin1995}
A.~H. Schatz and L.~B. Wahlbin.
\newblock {Interior maximum-norm estimates for finite element methods II}.
\newblock {\em Math. Comp.}, 64:907--928, 1995.

\bibitem{Thomee2006}
V.~Thom\'{e}e.
\newblock {\em {Galerkin Finite Element Methods for Parabolic Problems}}.
\newblock Springer-Verlag, New York, second edition, 2006.

\bibitem{ThomeeWahlbin2000}
V.~Thom\'{e}e and L.~B. Wahlbin.
\newblock {Stability and analyticity in maximum-norm for simplicial Lagrange
  finite element semidiscretizations of parabolic equations with Dirichlet
  boundary conditions}.
\newblock {\em Numer. Math.}, 87:373--389, 2000.

\bibitem{Wahlbin1991}
L.~B. Wahlbin.
\newblock {\em {Local behavior in finite element methods}}.
\newblock ``Handbook of Numerical Analysis II: Finite Element Methods'', pp.
  353-522, edited by P. G. Ciarlet and J. L. Lions. North-Holland, Amsterdam,
  1991.

\bibitem{Weis2001-2}
L.~Weis.
\newblock {Operator-valued Fourier multiplier theorems and maximal
  $L_p$-regularity}.
\newblock {\em Math. Ann.}, 319.

\bibitem{Weis2001-1}
L.~Weis.
\newblock {\em {A new approach to maximal $L_p$-regularity}}.
\newblock Lecture Notes in Pure and Applied Mathematics 215: Evolution
  Equations and Their Applications in Physical and Life Sciences, pp. 195--214,
  edited by G. Lumer. Marcel Dekker, 2001.

\bibitem{Yosida1980}
K.~Yosida.
\newblock {\em {Functional Analysis}}.
\newblock Springer-Verlag, New York, sixth edition, 1980.

\end{thebibliography}

}
\end{document}